\documentclass[10pt]{article}
\usepackage{amsmath}
\usepackage{latexsym}
\usepackage{amssymb}
\usepackage{amsfonts}
\usepackage{amsthm}

\newcommand{\Q}{{\mathbb{Q}}}

\usepackage{rotating}

\title{BIANCHI IDENTITIES \\
FOR THE RIEMANN AND WEYL TENSORS }
\author{J.-F. Pommaret \\ CERMICS, Ecole des Ponts ParisTech,\\ 6/8 Av. Blaise Pascal, 77455 Marne-la-Vall\'ee Cedex 02, France \\
E-mail: jean-francois.pommaret@wanadoo.fr, pommaret@cermics.enpc.fr \\
URL: http://cermics.enpc.fr/$\sim$pommaret/home.html }
\date{  }
\begin{document}
\maketitle

\noindent
{\bf ABSTRACT}  \\

The purpose of this paper is to revisit the Bianchi identities existing for the Riemann and Weyl tensors in the combined framework of the {\it formal
 theory} of systems of partial differential equations (Spencer cohomology, differential systems, formal integrability) and {\it Algebraic Analysis} (homological algebra, differential modules, duality). In particular, we prove that the $n^2(n^2-1)(n-2)/24$ generating Bianchi identities for the Riemann tensor are {\it first order} and can be easily described by means of the Spencer cohomology of the first order Killing symbol in arbitrary dimension $n\geq 2$. Similarly, the $n(n^2-1)(n+2)(n-4)/24$ generating Bianchi identities for the Weyl tensor are {\it first order} and can be easily described by means of the Spencer cohomology of the first order conformal Killing symbol in arbitrary dimension $n\geq 5$. As a {\it most surprising result}, the $9$ generating Bianchi identities for the Weyl tensor are of {\it second order} in dimension $n=4$ while the analogue of the Weyl tensor has $5$ components of {\it third order} in the metric with $3$ {\it first order} generating Bianchi identities in dimension $n=3$. The above results, {\it which could not be obtained otherwise}, are valid for any non-degenerate metric of constant riemannian curvature and do not depend on any conformal factor. They are checked in an Appendix produced by Alban Quadrat (INRIA, Lille) by means of computer algebra. We finally explain why the work of Lanczos and followers is not coherent with these results and must therefore be also revisited.\\


\noindent
{\bf 1) INTRODUCTION} \\

The language of {\it differential modules} has been recently introduced in applications as a way to understand the {\it structural properties} of systems of partial differential equations and the Poincar\'{e} duality between geometry and physics by using {\it adjoint operators} or variational calculus with differential constraints ([2],[23],[38]). In order to explain briefly the ideas of Lanczos as a way to justify the title of this paper, let us revisit briefly the foundation of n-dimensional elasticity theory as it can be found today in any textbook. If $x=(x^1,...,x^n)$ is a point in space and $\xi=({\xi}^1(x),...,{\xi}^n(x))$ is the displacement vector, lowering the indices by means of the Euclidean metric, we may introduce the "small" deformation tensor $\epsilon = ({\epsilon}_{ij}={\epsilon}_{ji}=(1/2)({\partial}_i{\xi}_j+{\partial}_j{\xi}_i))$ with $n(n+1)/2$ (independent) {\it components} $({\epsilon}_{i\leq j})$. If we study a part of a deformed body by means of a variational principle, we may introduce the local density of free energy $\varphi (\epsilon)=\varphi({\epsilon}_{ij}{\mid} i\leq j)$ and vary the total free energy $\Phi={\int}\varphi (\epsilon)dx$ with $dx=dx^1\wedge... \wedge dx^n$ by introducing ${\sigma}^{ij}=\partial \varphi/\partial {\epsilon}_{ij}$ for $i\leq j$ and "{\it deciding}" to define the stress tensor $\sigma$ by a symmetric matrix with ${\sigma}^{ij}={\sigma}^{ji}$ in a purely artificial way within such a variational principle. Indeed, the usual Cauchy Tetrahedron device (1828) assumes that each element of a boundary surface is acted on by a surface density of force ${\vec{\sigma}}$ with a linear dependence $\vec{\sigma}=({\sigma}^{ir}(x)n_r)$ on the outward normal unit vector $\vec{n}=(n_r)$ and does not make any assumption on the stress tensor. It is only by an equilibrium of forces and couples, namely the well known {\it phenomenological static torsor equilibrium}, that one can {\it prove} the symmetry of $\sigma$. However, if we assume this symmetry, we may now consider the summation $\delta \Phi={\int}{\sigma}^{ij}\delta{\epsilon}_{ij}dx={\int}{\sigma}^{ir}{\partial}_r\delta {\xi}_idx$. An integration by parts and a change of sign produce the integral ${\int}({\partial}_r{\sigma}^{ir})\delta {\xi}_idx$ leading to the stress equations ${\partial}_r{\sigma}^{ir}=0$. This classical approach to elasticity theory, based on invariant theory with respect to the group of rigid motions, cannot therefore describe equilibrium of torsors by means of a variational principle where the proper torsor concept is totally lacking. It is however widely used through the technique of " {\it finite elements} " where it can also be applied to electromagnetism (EM) with similar quadratic (piezoelectricity) or cubic (photoelasticity) lagrangian integrals. In this situation, the $4$-potential $A$ of EM is used in place of $\xi$ while the EM field $dA=F=(\vec{B}, \vec{E})$ is used in place of $\epsilon$.  \\

However, there exists another equivalent procedure dealing with a {\it variational calculus with constraint}. Indeed, as we shall see later on, the deformation tensor is not any symmetric tensor as it must satisfy $n^2(n^2-1)/12$ {\it Riemann} compatibility conditions (CC), that is the only condition ${\partial}_{22}{\epsilon}_{11}+{\partial}_{11}{\epsilon}_{22}-2{\partial}_{12}{\epsilon}_{12}=0$ when $n=2$. In this case, introducing the {\it Lagrange multiplier } $\lambda$, {\it we have to vary the new integral} $\int[{\varphi}(\epsilon) + \lambda ({\partial}_{22}{\epsilon}_{11}+{\partial}_{11}{\epsilon}_{22}-2{\partial}_{12}{\epsilon}_{12})]dx$ {\it for an arbitrary} $\epsilon$. Setting $\lambda= - \phi$, a double integration by parts now provides the parametrization ${\sigma}^{11}={\partial}_{22}\phi,{\sigma}^{12}={\sigma}^{21}=-{\partial}_{12}\phi, {\sigma}^{22}={\partial}_{11}\phi$ of the stress equations by means of the Airy function $\phi$ and the {\it formal adjoint} of the Riemann CC ([1],[26]). The same variational calculus with constraint may thus also be used in order to avoid the introduction of the EM potential $A$ by using the Maxwell equations $dF=0$ in place of the Riemann CC for $\epsilon$ but, in all these situations, {\it we have to eliminate the Lagrange multipliers or use them as potentials}.\\

In arbitrary dimension, the above compatibility conditions are nothing else but the linearized Riemann tensor in Riemannian geometry, a crucial mathematical tool in the theory of general relativity and a good reason for studying the work of Cornelius Lanczos (1893-1974) as it can be found in ([14],[15]) or in a few modern references ([5],[6],[7],[18],[36]). The starting point of Lanczos has been to take EM as a model in order to introduce a Lagrangian that should be quadratic in the Riemann tensor $({\rho}^k_{l,ij}\Rightarrow {\rho}_{ij}={\rho}^r_{i,rj}={\rho}_{ji}\Rightarrow \rho={\omega}^{ij}{\rho}_{ij})$ while considering it independently of its expression through the second order derivatives of a metric $({\omega}_{ij})$ with inverse $({\omega}^{ij})$ or the first order derivatives of the corresponding Christoffel symbols $({\gamma}^k_{ij})$. According to the previous paragraph, the corresponding variational calculus {\it must} involve PD constraints made by the Bianchi identities and {\it the new lagrangian to vary must therefore contain as many Lagrange multipliers as the number of generating Bianchi identities} that can be written under the form:  
\[   {\nabla}_r{\rho}^k_{l,ij}+{\nabla}_i{\rho}^k_{l,jr}+{\nabla}_j{\rho}^k_{l,ri}=0 \Rightarrow   {\nabla}_r{\rho}^r_{l,ij}={\nabla}_i{\rho}_{lj}-{\nabla}_j{\rho}_{li}       \]
Meanwhile, Lanczos and followers have been looking for a kind of "{\it parametrization} " by using the corresponding "{\it Lanczos potential} ", {\it exactly like the Lagrange multiplier has been used as an Airy potential for the stress equations}. However, we shall prove that the definition of a {\it Riemann candidate} cannot be done without the knowledge of the Spencer cohomology. Moreover, we have pointed out the existence of well known couplings between elasticity and electromagnetism, namely piezoelectricity and photoelasticity, which are showing that, in the respective Lagrangians, the EM field is on equal footing with the deformation tensor and {\it not} with the Riemann tensor. The {\it shift by one step backwards} that must be used in the physical interpretation of the differential sequences involved cannot therefore be avoided. Meanwhile, the {\it ordinary derivatives} ${\partial}_i$ can be used in place of the {\it covariant derivatives} ${\nabla}_i$ when dealing with the linearized framework as the Christoffel symbols vanish when Euclidean or Minkowskian metrics are used. \\

The next tentative of Lanczos has been to extend his approach to the Weyl tensor:  \\
\[  {\tau}^k_{l,ij} = {\rho}^k_{l,ij} - \frac{1}{(n-2)}({\delta}^k_i{\rho}_{lj} - {\delta}^k_j{\rho}_{li} +{\omega}_{lj}{\omega}^{ks}{\rho}_{si} - {\omega}_{li}{\omega}^{ks}{\rho}_{sj}) + \frac{1}{(n-1)(n-2)}({\delta}^k_i{\omega}_{lj} - {\delta}^k_j{\omega}_{li})\rho  \]
The main problem is now that the Spencer cohomology of the symbols of the conformal Killing equations, in particular the $2$-acyclicity, will be {\it absolutely 
needed} in order to provide the Weyl tensor and its relation with the Riemann tensor. It will follow that the CC for the Weyl tensor may not be first order contrary to the CC for the Riemann tensor made by the Bianchi identities, another reason for justifying the shift by one step already quoted. In order to provide an idea of the difficulty involved, let us define the following tensors:  \\
\[     Schouten= ({\sigma}_{ij}={\rho}_{ij} - \frac{1}{2(n-1)}{\omega}_{ij}\rho) \Rightarrow Cotton=({\sigma}_{k,ij}={\nabla}_i{\sigma}_{kj} - {\nabla}_j{\sigma}_{ki})  \]
An elementary but tedious computation allows to prove the formula:  \\
\[              {\nabla}_r{\tau}^r_{k,ij}=\frac{(n-3)}{(n-2)}{\sigma}_{k,ij}   \]
Then, of course, {\it if Einstein equations in vacuum are valid}, the Schouten and Cotton tensors vanish but the left member is by no way a differential identity for the Weyl tensor and {\it great care must be taken when mixing up mathematics with physics}.  \\

The author thanks Prof. Lars Andersson (Einstein Institute, Postdam) for having suggested him to study the {\it Lanczos potential} within this new framework and Alban Quadrat (INRIA, Lille), a specialist of control theory and computer algebra, for having spent time checking directly in an Appendix the many striking results contained in this paper.\\

\noindent
{\bf 2) HOMOLOGICAL ALGEBRA}  \\

We now need a few definitions and
results from homological algebra ([3],[9],[17],[23],[37]). In the following two classical theorems, $A,B,C,D,K,L,M,N,Q,R,S,T$ will be modules over a ring $A$ or 
vector spaces over a field $k$ and the linear maps are making the diagrams commutative. We start recalling the well known Cramer's rule for linear systems through
the exactness of the ker/coker sequence for modules. When $f:M\rightarrow N$ is a linear map (homomorphism), we introduce the so-called ker/coker long exact sequence:  \\
\[0\longrightarrow ker(f) \longrightarrow M
\stackrel{f}{\longrightarrow} N \longrightarrow coker(f) 
\longrightarrow 0 \]
In the case of vector spaces over a field $k$, we successively have 
$rk(f)=dim(im(f))$, $dim(ker(f))=dim(M)-rk(f)$ and $dim(coker(f))=dim(N)-rk(f)$ is the proper number of compatibility conditions. We obtain by substraction:
\[ dim(ker(f))-dim(M)+dim(N)-dim(coker(f))=0 \]
In the case of modules, we may replace the dimension by
the rank with $rk_A(M)=r$ when $F\simeq A^r$ is the greatest free submodule of $M$ and obtain the same relations because of the additive property of
the rank ([23],[24],[33]). The following theorems will be crucially used through the whole paper:\\

\noindent
{\bf SNAKE THEOREM  2.1}: When one has the following commutative diagram
resulting from the two central vertical short exact sequences by
exhibiting the three corresponding horizontal ker/coker exact sequences:
\[
\begin{array}{ccccccccccc}
 & & 0 & & 0 & & 0 & & & & \\
 & &\downarrow & & \downarrow & & \downarrow & & & & \\
0&\longrightarrow &K&\longrightarrow &A&\longrightarrow
&A'&\longrightarrow &Q&\longrightarrow &0\\
 & &\downarrow & &\;\;\;\downarrow \! f&
&\;\;\;\;\downarrow \! f'& &\downarrow & & \\
0&\longrightarrow &L&\longrightarrow &B&\longrightarrow
&B'&\longrightarrow &R& \longrightarrow &0 \\
 & &\downarrow & &\;\;\;\downarrow \! g & &\;\;\;\;
\downarrow \! g'& & \downarrow & & \\
0 & \longrightarrow &M& \longrightarrow &C&
\longrightarrow &C'& \longrightarrow &S& \longrightarrow
&0 \\
 & & & & \downarrow & & \downarrow & & \downarrow & & \\
 & & & &0& &0& &0& &
\end{array}
\]
then there exists a connecting map $M \longrightarrow Q$ both with a long
exact sequence:
\[0 \longrightarrow K \longrightarrow L \longrightarrow M
\longrightarrow Q \longrightarrow R \longrightarrow S
\longrightarrow 0.\]\\

\noindent
{\it Proof}: We construct the connecting map by using the
following succession of elements:
\[
\begin{array}{ccccccc}
 & & a & \cdots & a'& \longrightarrow &q \\
 & & \vdots & & \downarrow & & \\
 & & b & \longrightarrow &b' & & \\
 & & \downarrow & & \vdots & & \\
m& \longrightarrow &c& \cdots &0& & 
\end{array}
\]
Indeed, starting with $m\in M$, we may identify it with $c\in C$ in the
kernel of the next horizontal map. As $g$ is an epimorphism, we may find
$b\in B$ such that $c=g(b)$ and apply the next horizontal map to get
$b'\in B'$ in the kernel of $g'$ by the commutativity of the lower
square. Accordingly, there is a unique $a'\in A'$ such that
$b'=f'(a')$ and we may finally project $a'$ to $q\in Q$. The map is
well defined because, if we take another lift for $c$ in $B$, it will
differ from $b$ by the image under $f$ of a certain $a\in A$ having zero
image in $Q$ by composition. The remaining of the proof is similar. The above explicit procedure is called " {\it chase} " and will not be
repeated. \\
\hspace*{12cm}  Q.E.D. \\

We may now introduce {\it cohomology theory} through the following 
definition:\\

\noindent
{\bf DEFINITION  2.2}: If one has a sequence $L\stackrel{f}{\longrightarrow} M \stackrel{g}{\longrightarrow} N $, that is if $g\circ f=0$, 
then one may introduce the submodules $coboundary=B=im(f)\subseteq ker(g)=cocycle=Z\subseteq M$ and define 
the cohomology at $M$ to be the quotient $H=Z/B $. The sequence is said to be {\it exact} at $M$ if $im(f)=ker(g)$.  \\

\noindent
{\bf COHOMOLOGY THEOREM  2.3}: The following commutative diagram where the two central 
vertical sequences are long exact sequences and the horizontal lines are
ker/coker exact sequences:
\[        \begin{array}{ccccccccccccc}
 & &0& &0& &0& & & & & & \\
 & &\downarrow & & \downarrow & & \downarrow & & & & & &
\\
0&\longrightarrow &K&\longrightarrow &A&\longrightarrow
&A'&\longrightarrow &Q&\longrightarrow &0& & \\
 & &\downarrow & &\;\;\;\downarrow \! f &
&\;\;\;\;\downarrow \! f' & &\downarrow & & & & \\
0&\longrightarrow &L&\longrightarrow&B&\longrightarrow
&B'&\longrightarrow &R&\longrightarrow &0& & \\
\cdots &\cdots &\downarrow &\cdots &\;\;\;\downarrow \!
g &\cdots &\;\;\;\;\downarrow \! g'&\cdots
&\downarrow &\cdots &\cdots &\cdots &cut \\
0&\longrightarrow &M&\longrightarrow &C&\longrightarrow
&C'&\longrightarrow &S&\longrightarrow &0& & \\
 & &\downarrow & &\;\;\;\downarrow \! h &
&\;\;\;\;\downarrow \! h'& &\downarrow & & & & \\
0&\longrightarrow &N&\longrightarrow &D&\longrightarrow
&D'&\longrightarrow &T&\longrightarrow &0& & \\
 & & & &\downarrow & &\downarrow & &\downarrow & & & & \\
 & & & &0& &0& &0 & & & & 
\end{array}   \]
induces an isomorphism between the cohomology at $M$ in the left vertical
column and the kernel of the morphism $Q\rightarrow R$ in the right
vertical column.\\

\noindent
{\it Proof}: Let us ``cut'' the preceding diagram along the dotted line. We obtain the following two
commutative and exact diagrams with $im(g)=ker(h), im(g')=ker(h')$:
\[   \begin{array}{ccccccccccc}
 & &0& &0& &0 & & & \\
 & &\downarrow & &\downarrow & &\downarrow & & & & \\
 0&\longrightarrow &K&\longrightarrow &A&\longrightarrow
&A'&\longrightarrow &Q&\longrightarrow &0 \\
 & &\downarrow & &\;\;\;\downarrow \! f &
&\;\;\;\;\downarrow \! f' & &\downarrow & & \\
0&\longrightarrow &L&\longrightarrow &B&\longrightarrow
&B'&\longrightarrow &R&\longrightarrow &0 \\
 & &\downarrow & &\;\;\;\downarrow \! g &
&\;\;\;\;\downarrow \! g' & & & & \\
0&\longrightarrow & cocycle &\longrightarrow &im (g)
&\longrightarrow &im (g') & & & & \\
 & & & &\downarrow & &\downarrow & & & & \\
 & & & &0& &0& & & & 
\end{array}    \]

\[   \begin{array}{ccccccc}
 & &0& &0& &0 \\
 & &\downarrow & &\downarrow & &\downarrow \\
0&\longrightarrow & cocycle &\longrightarrow & ker (h) &\longrightarrow &ker (h')\\
 & &\downarrow & &\downarrow & &\downarrow \\
0&\longrightarrow &M&\longrightarrow &C&\longrightarrow
&C' \\
 & &\downarrow & &\;\;\;\downarrow \! h &
&\;\;\;\;\downarrow \! h' \\
0&\longrightarrow &N&\longrightarrow &D&\longrightarrow
&D' \\
 & & & &\downarrow & &\downarrow \\
 & & & &0& &0 
\end{array}   \]
Using the snake theorem, we successively obtain the following long exact sequences: \\
\[   \begin{array}{ccccc}
\Longrightarrow &\exists &\qquad &0\longrightarrow K
\longrightarrow L \stackrel{g}{\longrightarrow}
cocycle \longrightarrow Q \longrightarrow R &\qquad \\
\Longrightarrow &\exists &\qquad & 0 \longrightarrow
coboundary \longrightarrow cocycle \longrightarrow ker \,(Q\longrightarrow R) 
\longrightarrow 0 &\qquad \\
\Longrightarrow & &\qquad & cohomology \hspace{1mm}  \simeq
ker \,(Q \longrightarrow R) & 
\end{array}   \]
\hspace*{12cm}   Q.E.D.  \\

We finally quote for a later use:\\

\noindent
{\bf PROPOSITION  2.4}: If one has a short exact sequence:
\[0\longrightarrow
M'\stackrel{f}{\longrightarrow}M\stackrel{g}{\longrightarrow}M''
\longrightarrow 0  \]
then the following conditions are equivalent:\\
$\bullet$ There exists an epimorphism $u:M\rightarrow M'$ such that $u\circ f=id_{M'}$ ({\it left inverse} of $f$).\\
$\bullet$ There exists a monomorphism $v:M''\rightarrow M$ such that $g\circ v=id_{M''}$ ({\it right inverse} of $g$).\\

\noindent
{\bf DEFINITION  2.5}: In the above situation, we say that the short exact
sequence {\it splits}. The relation $f\circ u+v\circ g=id_M$ provides an isomorphism $(u,g):M\rightarrow M'\oplus M''$ with inverse $f+v:M'\oplus M"\rightarrow M$.
The short exact sequence $0 \rightarrow \mathbb{Z} \rightarrow  \mathbb{Q} \rightarrow \mathbb{Q}/\mathbb{Z} \rightarrow 0$ cannot split over $\mathbb{Z}$. \\

\noindent
{\bf 3) DIFFERENTIAL SYSTEMS}  \\

If $E$ is a vector bundle over the base manifold $X$ with projection $\pi$ and local coordinates $(x,y)=(x^i,y^k)$ projecting onto $x=(x^i)$ for $i=1,...,n$ and $k=1,...,m$, identifying a map with its graph, a (local) section $f:U\subset X \rightarrow E$ is such that $\pi\circ f =id$ on $U$ and we write $y^k=f^k(x)$ or simply $y=f(x)$. For any change of local coordinates $(x,y)\rightarrow (\bar{x}=\varphi(x),\bar{y}=A(x)y)$ on $E$, the change of section is $y=f(x)\rightarrow \bar{y}=\bar{f}(\bar{x})$ such that ${\bar{f}}^l(\varphi(x)\equiv A^l_k(x)f^k(x)$. The new vector bundle $E^*$ obtained by changing the {\it transition matrix} $A$ to its inverse $A^{-1}$ is called the {\it dual vector bundle} of $E$. In particular, let $T$ be the tangent vector bundle of vector fields on $X$, $T^*$ be the cotangent vector bundle of 1-forms on $X$ and $S_qT^*$ be the vector bundle of symmetric q-covariant tensors on $X$. Differentiating with respect to $x^i$ and using new coordinates $y^k_i$ in place of ${\partial}_if^k(x)$, we obtain ${\bar{y}}^l_r{\partial}_i{\varphi}^r(x)=A^l_k(x)y^k_i+{\partial}_iA^l_k(x)y^k$. Introducing a multi-index $\mu=({\mu}_1,...,{\mu}_n)$ with length $\mid \mu \mid={\mu}_1+...+{\mu}_n$ and prolonging the procedure up to order $q$, we may construct in this way, by patching coordinates, a vector bundle $J_q(E)$ over $X$, called the {\it jet bundle of order} $q$ with local coordinates $(x,y_q)=(x^i,y^k_{\mu})$ with $0\leq \mid\mu\mid \leq q$ and $y^k_0=y^k$. We have therefore epimorphisms ${\pi}^{q+r}_q:J_{q+r}(E)\rightarrow J_q(E), \forall q,r\geq 0$ and the short exact sequences $0 \rightarrow S_qT^*\otimes E \rightarrow J_q(E) \stackrel{{\pi}^q_{q-1}}{\longrightarrow} J_{q-1}(E) \rightarrow 0$. For a later use, we shall set $\mu+1_i=({\mu}_1,...,{\mu}_{i-1},{\mu}_i+1,{\mu}_{i+1},...,{\mu}_n)$ and define the operator $j_q:E \rightarrow J_q(E):f \rightarrow j_q(f)$ on sections by the local formula $j_q(f):(x)\rightarrow({\partial}_{\mu}f^k(x)\mid 0\leq \mid\mu\mid \leq q,k=1,...,m)$. Moreover, a jet coordinate $y^k_{\mu}$ is said to be of {\it class} $i$ if ${\mu}_1=...={\mu}_{i-1}=0, {\mu}_i\neq 0$. We finally inroduce the {\it Spencer operator} $D:J_{q+1}(E) \rightarrow T^*\otimes J_q(E):f_{q+1} \rightarrow j_1(f_q)-f_{q+1}$ with $(Df_{q+1})^k_{\mu,i}={\partial}_if^k_{\mu}-f^k_{\mu+1_i}$.\\

\noindent
{\bf DEFINITION  3.1}:  A {\it system} of PD equations of order $q$ on $E$ is a vector subbundle $R_q\subset J_q(E)$ locally defined by a constant rank system of linear equations for the jets of order $q$ of the form $ a^{\tau\mu}_k(x)y^k_{\mu}=0$. Its {\it first prolongation} $R_{q+1}\subset J_{q+1}(E)$ will be defined by the equations $ a^{\tau\mu}_k(x)y^k_{\mu}=0, a^{\tau\mu}_k(x)y^k_{\mu+1_i}+{\partial}_ia^{\tau\mu}_k(x)y^k_{\mu}=0$ which may not provide a system of constant rank as can easily be seen for $xy_x-y=0 \Rightarrow xy_{xx}=0$ where the rank drops at $x=0$.\\

The next definition of {\it formal integrability} will be crucial for our purpose.\\

\noindent
{\bf DEFINITION  3.2}: A system $R_q$ is said to be {\it formally integrable} if the $R_{q+r}$ are vector bundles $\forall r\geq 0$ (regularity condition) and no new equation of order $q+r$ can be obtained by prolonging the given PD equations more than $r$ times, $\forall r\geq 0$ or, equivalently, we have induced epimorphisms ${\pi}^{q+r+1}_{q+r}:R_{q+r+1} \rightarrow R_{q+r}, \forall r\geq 0$ allowing to compute " {\it step by step} " formal power series solutions.\\

A formal test first sketched by C. Riquier in 1910, has been improved by M. Janet in 1920 ([10],[19]) and by E. Cartan in 1945 ([4]), finally rediscovered in 1965, totally independently, by B. Buchberger who introduced Gr\"{o}bner bases, using the name of his thesis advisor ([25]). However all these tentatives have been largely superseded and achieved in an intrinsic way, again totally independently of the previous approaches, by D.C. Spencer in 1965 ([19],[22],[39]). \\

\noindent
{\bf DEFINITION  3.3}: The family $g_{q+r}$ of vector spaces over $X$ defined by the purely linear equations $ a^{\tau\mu}_k(x)v^k_{\mu+\nu}=0$ for $ \mid\mu\mid= q, \mid\nu\mid =r $ is called the {\it symbol} at order $q+r$ and only depends on $g_q$.\\

The following procedure, {\it where one may have to change linearly the independent variables if necessary}, is the key towards the next definition which is intrinsic even though it must be checked in a particular coordinate system called $\delta$-{\it regular} (See [19],[22],[23] and [39] for more details):\\

\noindent
$\bullet$ {\it Equations of class} $n$: Solve the maximum number ${\beta}^n_q$ of equations with respect to the jets of order $q$ and class $n$. Then call $(x^1,...,x^n)$ {\it multiplicative variables}.\\
\[  - - - - - - - - - - - - - - - -  \]
$\bullet$ {\it Equations of class} $i$: Solve the maximum number of {\it remaining} equations with respect to the jets of order $q$ and class $i$. Then call $(x^1,...,x^i)$ {\it multiplicative variables} and $(x^{i+1},...,x^n)$ {\it non-multiplicative variables}.\\
\[ - - - - - - - - - - - - - - - - - \]
$\bullet$ {\it Remaining equations equations of order} $\leq q-1$: Call $(x^1,...,x^n)$ {\it non-multiplicative variables}.\\

\noindent
{\bf DEFINITION  3.4}: The above multiplicative and non-multiplicative variables can be visualized respectively by integers and dots in the corresponding {\it Janet board}. A system of PD equations is said to be {\it involutive} if its first prolongation can be achieved by prolonging its equations only with respect to the corresponding multiplicative variables. The following numbers are called {\it characters}:  \\
\[ {\alpha}^i_q=m(q+n-i-1)!/((q-1)!(n-i)!)-{\beta}^i_q , \hspace{3mm} \forall 1\leq i \leq n\hspace{3mm} \Rightarrow \hspace{3mm}{\alpha}^1_q\geq ... \geq {\alpha}^n_q \]
For an involutive system, $(y^{{\beta}^n_q +1},...,y^m)$ can be given arbitrarily.  \\

For an involutive system of order $q$ in the above {\it solved form}, we shall use to denote by $y_{pri}$ the {\it principal jet coordinates}, namely the leading terms of the solved equations in the sense of involution, and any formal derivative of a principal jet coordinate is again a principal jet coordinate. The remaining jet coordinates will be called {\it  parametric jet coordinates} and denoted by $y_{par}$.\\

\noindent
{\bf PROPOSITION  3.5}: Using the Janet board and the definition of involutivity, we get:  \\
\[   dim(g_{q+r})={\sum}_{i=1}^n\frac{(r+i-1)!}{r!(i-1)!}{\alpha}^i_q   \Rightarrow  dim(R_{q+r})=dim(R_{q-1})+{\sum}_{i=1}^n\frac{(r+i)!}{r!i!}{\alpha}^i_q  \]   \\

Let now ${\wedge}^sT^*$ be the vector bundle of s-forms on $X$ with usual bases $\{dx^I=dx^{i_1}\wedge ... \wedge dx^{i_s}\}$ where we have set $I=(i_1< ... <i_s)$. . Moreover, if  $\xi,\eta\in T$ are two vector fields on $X$, we may define their {\it bracket} $[\xi,\eta]\in T$ by the local formula $([\xi,\eta])^i(x)={\xi}^r(x){\partial}_r{\eta}^i(x)-{\eta}^s(x){\partial}_s{\xi}^i(x)$ leading to the {\it Jacobi identity} $[\xi,[\eta,\zeta]]+[\eta,[\zeta,\xi]]+[\zeta,[\xi,\eta]]=0, \forall \xi,\eta,\zeta \in T$. We may finally introduce the {\it exterior derivative} $d:{\wedge}^rT^*\rightarrow {\wedge}^{r+1}T^*:\omega={\omega}_Idx^I \rightarrow d\omega={\partial}_i{\omega}_Idx^i\wedge dx^I$ with $d^2=d\circ d\equiv 0$ in the {\it Poincar\'{e} sequence}:\\
\[  {\wedge}^0T^* \stackrel{d}{\longrightarrow} {\wedge}^1T^* \stackrel{d}{\longrightarrow} {\wedge}^2T^* \stackrel{d}{\longrightarrow} ... \stackrel{d}{\longrightarrow} {\wedge}^nT^* \longrightarrow 0  \]

In a purely algebraic setting, one has ([19],[22],[23],[24],[39]):  \\

\noindent
{\bf PROPOSITION  3.6}: There exists a map $\delta:{\wedge}^sT^*\otimes S_{q+1}T^*\otimes E\rightarrow {\wedge}^{s+1}T^*\otimes S_qT^*\otimes E$ which restricts to $\delta:{\wedge}^sT^*\otimes g_{q+1}\rightarrow {\wedge}^{s+1}T^*\otimes g_q$ and ${\delta}^2=\delta\circ\delta=0$.\\

{\it Proof}: Let us introduce the family of s-forms $\omega=\{ {\omega}^k_{\mu}=v^k_{\mu,I}dx^I \}$ and set $(\delta\omega)^k_{\mu}=dx^i\wedge{\omega}^k_{\mu+1_i}$. We obtain at once $({\delta}^2\omega)^k_{\mu}=dx^i\wedge dx^j\wedge{\omega}^k_{\mu+1_i+1_j}=0$.\\
\hspace*{12cm} Q.E.D.  \\

The kernel of each $\delta$ in the first case is equal to the image of the preceding $\delta$ but this may no longer be true in the restricted case and we set (See [22], p 85-88 for more details):   \\

\noindent
{\bf DEFINITION  3.7}: We denote by $B^s_{q+r}(g_q)\subseteq Z^s_{q+r}(g_q)$ and $H^s_{q+r}(g_q)=Z^s_{q+r}(g_q)/B^s_{q+r}(g_q)$ respectively the coboundary space, cocycle space and cohomology space at ${\wedge}^sT^*\otimes g_{q+r}$ of the restricted $\delta$-sequence which only depend on $g_q$ and may not be vector bundles. The symbol $g_q$ is said to be {\it s-acyclic} if $H^1_{q+r}=...=H^s_{q+r}=0, \forall r\geq 0$, {\it involutive} if it is n-acyclic and {\it finite type} if $g_{q+r}=0$ becomes trivially involutive for r large enough. For a later use, we notice that a symbol $g_q$ is involutive {\it and} of finite type if and only if $g_q=0$. Finally, $S_qT^*\otimes E$ is involutive $\forall q\geq 0$ if we set $S_0T^*\otimes E=E$. We shall prove later on that {\it any} symbol $g_q$ is $1$-acyclc. \\

\noindent
{\bf CRITERION THEOREM  3.8}: If ${\pi}^{q+1}_q:R_{q+1} \rightarrow R_q$ is an epimorphism of vector bundles and $g_q$ is $2$-acyclic (involutive), then $R_q$ is formally integrable (involutive).  \\

\noindent
{\bf EXAMPLE  3.9}: The system $R_2$ defined by the three PD equations\\
\[     {\Phi}^3 \equiv y_{33}=0, \hspace{5mm} {\Phi}^2 \equiv y_{23}-y_{11}=0,\hspace{5mm}{\Phi}^1 \equiv  y_{22}=0  \]
is homogeneous and thus automatically formally integrable but $g_2$ and $g_3$ are not involutive though finite type because $g_4=0$ and the sequence $0 \rightarrow {\wedge}^3T^*\otimes g_3 \rightarrow 0$ is not exact. Elementary computations of ranks of matrices shows that the $\delta$-map:\\
\[    0\rightarrow  {\wedge}^2T^*\otimes g_3  \stackrel{\delta}{\longrightarrow} {\wedge}^3T^*\otimes g_2 \rightarrow 0  \]
is a $3\times 3$ isomorphism and thus $g_3$ is 2-acyclic with $dim(g_3)=1$, a {\it crucial intrinsic} property totally absent from any "old" work and quite more easy to handle than its Koszul dual. We invite the reader to treat similarly the system $y_{33}-y_{11}=0, y_{23}=0, y_{22}-y_{11}=0$ and compare.\\

The main use of involution is to construct differential sequences that are made up by successive {\it compatibility conditions} (CC) of order one. In particular, when $R_q$ is involutive, the differential operator ${\cal{D}}:E\stackrel{j_q}{\rightarrow} J_q(E)\stackrel{\Phi}{\rightarrow} J_q(E)/R_q=F_0$ of order $q$ with space of solutions $\Theta\subset E$ is said to be {\it involutive} and one has the canonical {\it linear Janet sequence} ([22], p 144):\\
\[  0 \longrightarrow  \Theta \longrightarrow E \stackrel{\cal{D}}{\longrightarrow} F_0 \stackrel{{\cal{D}}_1}{\longrightarrow}F_1 \stackrel{{\cal{D}}_2}{\longrightarrow} ... \stackrel{{\cal{D}}_n}{\longrightarrow} F_n \longrightarrow 0   \]
where each other operator is first order involutive and generates the CC of the preceding one with the {\it Janet bundles} $F_r={\wedge}^rT^*\otimes J_q(E)/({\wedge}^rT^*\otimes R_q+\delta ({\wedge}^{r-1}T^*\otimes S_{q+1}T^*\otimes E))$. As the Janet sequence can be "cut at any place", that is can also be constructed anew from any intermediate operator, {\it the numbering of the Janet bundles has nothing to do with that of the Poincar\'{e} sequence for the exterior derivative}, contrary to what many physicists  still believe ($n=3$ with ${\cal{D}}=div$ provides the simplest example). Moreover, the fiber dimension of the Janet bundles can be computed at once inductively from the board of multiplicative and non-multiplicative variables that can be exhibited for $\cal{D}$ by working out the board for ${\cal{D}}_1$ and so on. For this, the number of rows of this new board is the number of dots appearing in the initial board while the number $nb(i)$ of dots in the column $i$ just indicates the number of CC of class $i$ for $i=1, ... ,n$ with $nb(i) < nb(j), \forall i<j$. \\

\noindent
{\bf MAIN THEOREM  3.10}: When $R_q\subset J_q(E)$ is {\it not} involutive but formally integrable and {\it its symbol} $g_q$ {\it becomes} $2$-{\it acyclic after exactly} $s$ {\it prolongations}, the generating CC are of order $s+1$ (See [22], Example 6, p 120 and previous Example). \\

\noindent
{\it Proof}: We may introduce the canonical epimorphism $\Phi=J_q(E) \rightarrow J_q(E)/R_q=F_0$ and denote by ${\cal{D}}=\Phi \circ j_q:E \rightarrow F_0$ the corresponding differential operator. As before, we may write formally ${\Phi}^{\tau}(x,y_q)\equiv a^{\tau}_{\mu,k}(x)y^k_{\mu} = z^{\tau}$, obtain $d_i{\Phi}^{\tau}\equiv a^{\tau}_{\mu,k}(x)y^k_{\mu+1_i} + {\partial}_ia^{\tau}_{\mu,k}(x)y^k_{\mu}=z^{\tau}_i$ for the first prolongation ${\rho}_1(\Phi):J_{q+1}(E) \rightarrow j_1(F_0)$ and so on with ${\rho}_r(\Phi):J_{q+r}(E) \rightarrow J_r(F_0)$ defined by $ d_{\nu}{\Phi}^{\tau} =z^{\tau}_{\nu}$ with $0 \leq \mid \nu \mid \leq r$. Setting $B_r=im ({\rho}_r(\Phi))\subseteq J_r(F_0)$, we may introduce the canonical epimorphism $\Psi:J_{s+1}(F_0) \rightarrow J_{s+1}(F_0)/B_{s+1}=F_1$. Taking into account the formal integrability of $R_q$ ({\it care}), we obtain by composition of jets the following commutative {\it prolongation diagrams} $\forall r\geq 1$:  \\

\[ \begin{array}{rcccccccl}
   & 0                    &  & 0                      &  & 0                       & & 0                    &\\
   & \downarrow &  &  \downarrow  &  &  \downarrow  &  & \downarrow &  \\
0    \rightarrow & g_{q+r+s+1} & \rightarrow & S_{q+r+s+1}T^*\otimes E & \stackrel{{\sigma}_{r+s+1}(\Phi)}{\longrightarrow} & S_{r+s+1}T^*\otimes F_0 & 
\stackrel{{\sigma}_r(\Psi)}{\longrightarrow} & S_rT^*\otimes F_1 &   \\  
   & \downarrow &  &  \downarrow  &  &  \downarrow  &  & \downarrow &   \\
0    \rightarrow & R_{q+r+s+1} & \rightarrow & J_{q+r+s+1}(E) & \stackrel{{\rho}_{r+s+1}(\Phi)}{\longrightarrow} & J_{r+s+1}( F_0) & 
\stackrel{{\rho}_r(\Psi)}{\longrightarrow} & J_r(F_1) &   \\   
 & \downarrow &  &  \downarrow  &  &  \downarrow  &  & \downarrow &  \\
0    \rightarrow & R_{q+r+s} & \rightarrow & J_{q+r+s}(E) & \stackrel{{\rho}_{r+s}(\Phi)}{\longrightarrow} & J_{r+s}( F_0) & 
\stackrel{{\rho}_{r-1}(\Psi)}{\longrightarrow} & J_{r-1}(F_1) &   \\   
& \downarrow &  &  \downarrow  &  &  \downarrow  &  & \downarrow &   \\   
& 0                    &  & 0                      &  & 0                       & & 0                    &        
\end{array}  \]
\noindent
and the only thing we know is that the bottom sequence is exact for $r=1$ by construction and that the upper induced sequence is exact when $r=0$ when 
${\sigma}_0(\Psi)=\sigma(\Psi)$ is the restriction of $\Psi$ to $S_{s+1}T^*\otimes F_0\subset J_{s+1}(F_0)$ after a chase in the following commutative diagram:  \\
\[ \begin{array}{rcccccccl}
   & 0                    &  & 0                      &  & 0                       & & 0                    &\\
   & \downarrow &  &  \downarrow  &  &  \downarrow  &  & \downarrow &  \\
0    \rightarrow & g_{q+s+1} & \rightarrow & S_{q+s+1}T^*\otimes E & \stackrel{{\sigma}_{s+1}(\Phi)}{\longrightarrow} & S_{s+1}T^*\otimes F_0 & 
\stackrel{{\sigma}(\Psi)}{\longrightarrow} &  F_1 &   \\  
   & \downarrow &  &  \downarrow  &  &  \downarrow  &  & \downarrow &   \\
0    \rightarrow & R_{q+s+1} & \rightarrow & J_{q+s+1}(E) & \stackrel{{\rho}_{s+1}(\Phi)}{\longrightarrow} & J_{s+1}( F_0) & 
\stackrel{\Psi}{\longrightarrow} & F_1 &   \\   
 & \downarrow &  &  \downarrow  &  &  \downarrow  &  & &  \\
0    \rightarrow & R_{q+s} & \rightarrow & J_{q+s}(E) & \stackrel{{\rho}_s(\Phi)}{\longrightarrow} & J_s( F_0) &  & &   \\   
& \downarrow &  &  \downarrow  &  &  \downarrow  &  &  &   \\   
& 0                    &  & 0                      &  & 0                       & &                     &        
\end{array}  \]
because $R_q$ is formally integrable ({\it care}). Appying now the $\delta$-maps to the upper row of the previous prolongation diagram and proceeding by induction, starting from $r=1$, we shall prove that the upper row is exact. Indeed, setting $h_{r+s}=im({\sigma}_{r+s}(\Phi))\subseteq S_{r+s}T^*\otimes F_0$, we may cut the full commutative diagram thus obtained as in the proof of the previous "cohomology theorem" into the following two commutative diagrams:  \\
\[ \begin{array}{rcccccl}
   & 0                    &  & 0                      &  & 0                       & \\
  & \downarrow &  &  \downarrow  &  &  \downarrow  &  \\
0    \rightarrow & g_{q+r+s+1} & \rightarrow & S_{q+r+s+1}T^*\otimes E & \rightarrow & h_{r+s+1} &  \rightarrow 0  \\  
   & \downarrow &  &  \downarrow  &    & \downarrow &  \\                                                                                   
0    \rightarrow & T^*\otimes g_{q+r+s} & \rightarrow & T^*\otimes S_{q+r+s}T^*\otimes E & \rightarrow & T^*\otimes h_{r+s} & \rightarrow 0 \\  
   & \downarrow &  &  \downarrow  &  & \downarrow &  \\
0    \rightarrow & {\wedge}^2T^*\otimes g_{q+r+s-1} & \rightarrow & {\wedge}^2T^*\otimes S_{q+r+s-1}T^*\otimes E&\rightarrow&
 {\wedge}^2T^*\otimes h_{r+s-1} &\rightarrow 0   \\  
& \downarrow &  &  \downarrow  &  &  &     \\
0    \rightarrow & {\wedge}^3T^*\otimes S_{q+r+s-2}T^*\otimes E & = & {\wedge}^3T^*\otimes S_{q+r+s-2}T^*\otimes E & & 
\end{array}  \]

\[  \begin{array}{rcccccl}
 & 0 & & 0 & & 0 &   \\
 & \downarrow & & \downarrow & & \downarrow &  \\
0 \rightarrow & h_{r+s+1} & \rightarrow &S_{r+s+1}T^*\otimes F_0 & \rightarrow & S_rT^*\otimes F_1 &   \\  
 & \downarrow &  &\downarrow     &  & \downarrow &  \\
0 \rightarrow & T^*\otimes h_{r+s} &\rightarrow & T^* \otimes S_{r+s}T^*\otimes F_0 & \rightarrow & T^*\otimes S_{r-1}T^*\otimes F_1 & \\
&  \downarrow  &  & \downarrow  & &  &  \\
0 \rightarrow & {\wedge}^2T^*\otimes h_{r+s-1} & \rightarrow &{\wedge}^2T^*\otimes S_{r+s-1}T^*\otimes F_0 &  & &     
\end{array}  \]
An easy chase is showing that $g_q$ is always $1$-acyclic and that we have an induced monomorphism $0\rightarrow h_{r+s+1}\rightarrow T^*\otimes h_{r+s}$. The crucial result that no classical approach could provide is that, whenever $g_{q+s}$ is $2$-acyclic, then the full right column of the first diagram is also exact or, equivalently, $h_{r+s+1}\subseteq S_{r+s+1}T^*\otimes F_0$ is the $r$-prolongation of the symbol $h_{s+1}\subseteq S_{s+1}T^*\otimes F_0$. Using finally the second diagram, it follows by induction and a chase that the upper row is exact whenever the central row is exact, a result achieving the first part of the proof.  \\
We may also use an inductive chase in the full diagram, showing directly that the cohomology at $S_{r+s+1}T^* \otimes F_0$ of the upper sequence:  \\
\[   0    \rightarrow g_{q+r+s+1}  \rightarrow  S_{q+r+s+1}T^*\otimes E  \stackrel{{\sigma}_{r+s+1}(\Phi)}{\longrightarrow}  S_{r+s+1}T^*\otimes F_0 
\stackrel{{\sigma}_r(\Psi)}{\longrightarrow}  S_rT^*\otimes F_1   \]  
is the same as the $\delta$-cohomology of the left column at ${\wedge}^2T^*\otimes g_{q+r+s-1}$ because all the other vertical $\delta$-sequences are exact.  \\
Finally, starting from the $ker/coker$ long exact sequence allowing to define $\Psi$ and ending with $F_1$ while taking into account that the upper symbol row of the prolongation diagram is exact, we deduce by an inductive chase that the central row is also exact. It follows that $B_{r+s+1}$ is the $r$-prolongation of $B_{s+1} $ which 
is  formally integrable because a chase shows that $B_{r+s+1}$ projects onto $B_{r+s}$, $\forall r\geq 1$.The case of an involutive symbol can be studied similarly by choosing $s=0$ and explains why all the CC operators met in the Janet sequence are first order involutive operators.  \\

\hspace*{12cm}  Q.E.D.  \\

 As we shall see through explicit examples, in particular the conformal Killing system, there is no rule in general in order to decide about the minimum number $s'\geq 0$ such that $h_{s+s'+1}$ becomes $2$-acyclic in order to repeat the above procedure. However, replacing $r$ by $s'+2$ and chasing in the first of the last two diagrams, we have:  \\

\noindent
{\bf COROLLARY  3.11}: The symbol $h_{s+s'+1}$ {\it becomes} $2$-acyclic whenever the symbol $g_{q+s+s'}$ {\it becomes} $3$-acyclic.\\

\noindent
{\bf DEFINITION  3.12}: More generally, a differential sequence is said to be {\it formally exact} if each operator generates the CC of the operator preceding it.\\

\noindent
{\bf EXAMPLE  3.13}: ([16],$\S 38$, p 40) The second order system $y_{11}=0, y_{13}-y_2=0$ is neither formally integrable nor involutive. Indeed, we get $d_3y_{11}-d_1(y_{13}-y_2)=y_{12}$ and $d_3y_{12}-d_2(y_{13}-y_2)=y_{22}$, that is to say {\it each first and second} prolongation does bring a new second order PD equation. Considering the new system $y_{22}=0, y_{12}=0, y_{13}-y_2=0, y_{11}=0$, the (evident !) permutation of coordinates $(1,2,3)\rightarrow (3,2,1)$ provides the following involutive second order system with one equation of class $3$, $2$ equations of class $2$ and $1$ equation of clas $1$:   \\
\[  \left\{  \begin{array}{lcl}
{\Phi}^4 \equiv y_{33}  & = & 0  \\
{\Phi}^3 \equiv y_{23} & = & 0  \\
{\Phi}^2 \equiv y_{22} & = &  0  \\
{\Phi}^1 \equiv y_{13}-y_2 & = &  0 
\end{array}
\right. \fbox{$\begin{array}{lll}
1 & 2 & 3 \\
1 & 2 & \bullet \\
1 & 2 & \bullet \\
1 & \bullet & \bullet
\end{array}$}  \]
We have ${\alpha}^3_2=0,{\alpha}^2_2=0,{\alpha}^1_2=2$ and we get therefore the (formally exact) Janet sequence:    
\[    0 \longrightarrow  \Theta \longrightarrow 1 \longrightarrow 4 \longrightarrow 4 \longrightarrow 1  \longrightarrow  0    \]
However, keeping only ${\Phi}^1$ and ${\Phi}^4$ while using the fact that $d_{33}$ commutes with $d_{13}-d_2$, we get the formally exact sequence $ 0 \rightarrow  \Theta \rightarrow 1 \rightarrow 2 \rightarrow 1 \rightarrow 0$ {\it which is not a Janet sequence} though the Euler-Poincar\'{e} characteristics vanishes in both cases with $1-4+4-1=1-2+1=0$ ([22], p 159 and [23]). \\

\noindent
{\bf EXAMPLE  3.14}: Coming back to Example 3.9 while intoducing the three second order operators $P=d_{22},Q=d_{23}-d_{11}, R=d_{33}$ which are commuting between themselves, we have now $q=2, s=1$ and we obtain the second order CC:  \\
\[    {\Psi}^3\equiv P{\Phi}^2 - Q{\Phi}^1=0, \hspace{5mm} {\Psi}^2 \equiv R{\Phi}^1 - P{\Phi}^3=0, \hspace{5mm} {\Psi}^1=Q{\Phi}^3 - R{\Phi}^2=0 \]
Exactly like in the Poincar\'{e} sequence, we finally get the new second order CC:  \\
\[        P{\Psi}^1 + Q{\Psi}^2 + R{\Psi}^3 =0   \]
Writing out only the number of respective equations, we obtain the formally exact differential sequence with vanishing Euler-Poincar\'{e} characteristics: \\
\[       0 \rightarrow \Theta \rightarrow 1 \rightarrow 3 \rightarrow 3 \rightarrow 1 \rightarrow 0   \]
{\it which is not a Janet sequence}. We let the reader check that $h_2$ is {\it not} $2$-acyclic but that $h_3$ is $2$-acyclic and thus $s'=1$ because $2+1+1=4$ and $g_4=0$ is trivially involutive. A similar situation will be met with the conformal Killing equations.  \\

We may finally extend the restriction $D:R_{q+1} \rightarrow T^*\otimes R_q$ of the Spencer operator to:
\[ D:{\wedge}^rT^*\otimes R_{q+1} \rightarrow {\wedge}^{r+1}T^* \otimes R_q: \alpha \otimes f_{q+1} \rightarrow d\alpha \otimes f_q +(-1)^r \alpha \wedge Df_{q+1}  \Rightarrow D^2=D\circ D \equiv 0\]
in order to construct the {\it first Spencer sequence} which is another resolution of $\Theta$ because the kernel of the first $D$ is such that $f_{q+1}\in R_{q+1},Df_{q+1}=0 \Leftrightarrow f_{q+1}=j_{q+1}(f), f\in \Theta$ when $q$ is large enough. This standard notation for the Spencer operator must not be confused with the same notation used in the next section for the ring of differential operators but the distinction will always be pointed out whenever a confusion could exist.  \\

\noindent
{\bf 4) DIFFERENTIAL MODULES}  \\

Let $K$ be a {\it differential field}, that is a field containing $\mathbb{Q}$ with $n$ commuting {\it derivations} $\{{\partial}_1,...,{\partial}_n\}$ with ${\partial}_i{\partial}_j={\partial}_j{\partial}_i={\partial}_{ij}, \forall i,j=1,...,n$ such that ${\partial}_i(a+b)={\partial}_ia+{\partial}_ib, \hspace{2mm} {\partial}_i(ab)=({\partial}_ia)b+a{\partial}_ib, \forall a,b\in K$ and ${\partial}_i(1/a)= - (1/a^2){\partial}_ia, \forall a\in K$. Using an implicit summation on multi-indices, we may introduce the (noncommutative) {\it ring of differential operators} $D=K[d_1,...,d_n]=K[d]$ with elements $P=a^{\mu}d_{\mu}$ such that $\mid \mu\mid<\infty$ and $d_ia=ad_i+{\partial}_ia$. The highest value of ${\mid}\mu {\mid}$ with $a^{\mu}\neq 0$ is called the {\it order} of the {\it operator} $P$ and the ring $D$ with multiplication $(P,Q)\longrightarrow P\circ Q=PQ$ is filtred by the order $q$ of the operators. We have the {\it filtration} $0\subset K=D_0\subset D_1\subset  ... \subset D_q \subset ... \subset D_{\infty}=D$. Moreover, it is clear that $D$, as an algebra, is generated by $K=D_0$ and $T=D_1/D_0$ with $D_1=K\oplus T$ if we identify an element $\xi={\xi}^id_i\in T$ with the vector field $\xi={\xi}^i(x){\partial}_i$ of differential geometry, but with ${\xi}^i\in K$ now. It follows that $D={ }_DD_D$ is a {\it bimodule} over itself, being at the same time a left $D$-module ${ }_DD$ by the composition $P \longrightarrow QP$ and a right $D$-module $D_D$ by the composition $P \longrightarrow PQ$ with $D_rD_s=D_{r+s}, \forall r,s \geq 0$. \\

If we introduce {\it differential indeterminates} $y=(y^1,...,y^m)$, we may extend $d_iy^k_{\mu}=y^k_{\mu+1_i}$ to ${\Phi}^{\tau}\equiv a^{\tau\mu}_ky^k_{\mu}\stackrel{d_i}{\longrightarrow} d_i{\Phi}^{\tau}\equiv a^{\tau\mu}_ky^k_{\mu+1_i}+{\partial}_ia^{\tau\mu}_ky^k_{\mu}$ for $\tau=1,...,p$. Therefore, setting $Dy^1+...+Dy^m=Dy\simeq D^m$ and calling $I=D\Phi\subset Dy$ the {\it differential module of equations}, we obtain by residue the {\it differential module} or $D$-{\it module} $M=Dy/D\Phi$, denoting the residue of $y^k_{\mu}$ by ${\bar{y}}^k_{\mu}$ when there can be a confusion. Introducing the two free differential modules $F_0\simeq D^{m_0}, F_1\simeq D^{m_1}$, we obtain equivalently the {\it free presentation} $F_1\stackrel{d_1}{\longrightarrow} F_0 \rightarrow M \rightarrow 0$ of order $q$ when $m_0=m, m_1=p$ and $d_1={\cal{D}}=\Phi \circ j_q$ with $(P_1, ... ,P_p) \rightarrow (P_1, ... , P_p) \circ {\cal{D}} =(Q_1, ... , Q_m)$. We shall moreover assume that ${\cal{D}}$ provides a {\it strict morphism}, namely that the corresponding system $R_q$ is formally integrable. It follows that $M$ can be endowed with a {\it quotient filtration} obtained from that of $D^m$ which is defined by the order of the jet coordinates $y_q$ in $D_qy$. We have therefore the {\it inductive limit} $0=M_{-1} \subseteq M_0 \subseteq M_1 \subseteq ... \subseteq M_q \subseteq ... \subseteq M_{\infty}=M$ with $d_iM_q\subseteq M_{q+1}$ but it is important to notice that $D_rD_q=D_{q+r} \Rightarrow D_rM_q= M_{q+r}, \forall q,r\geq 0 \Rightarrow M=DM_q, \forall q\geq 0$ {\it in this particular case}. It also follows from noetherian arguments and involution that $D_ rI_q=I_{q+r}, \forall r\geq 0$ though we have in general only $D_rI_s\subseteq I_{r+s}, \forall r\geq 0, \forall s<q$. As $K\subset D$, we may introduce the {\it forgetful functor} $for : mod(D) \rightarrow mod(K): { }_DM \rightarrow { }_KM$. \\

More generally, introducing the successive CC as in the preceding section while changing slightly the numbering of the respective operators, we may finally obtain the {\it free resolution} of $M$, namely the exact sequence $\hspace{5mm} ... \stackrel{d_3}{\longrightarrow} F_2  \stackrel{d_2}{\longrightarrow} F_1 \stackrel{d_1}{\longrightarrow}F_0\longrightarrow M \longrightarrow 0 $. In actual practice, {\it one must never forget that} ${\cal{D}}=\Phi \circ j_q$ {\it acts on the left on column vectors in the operator case and on the right on row vectors in the module case}. Also, with a slight abuse of language, when ${\cal{D}}=\Phi \circ j_q$ is involutive as in section 3 and thus $R_q=ker( \Phi)$ is involutive, one should say that $M$ has an {\it involutive presentation} of order $q$ or that $M_q$ is {\it involutive}. \\

\noindent
{\bf DEFINITION  4.1}: Setting $P=a^{\mu}d_{\mu}\in D  \stackrel{ad}{\longleftrightarrow} ad(P)=(-1)^{\mid\mu\mid}d_{\mu}a^{\mu}   \in D $, we have $ad(ad(P))=P$ and $ad(PQ)=ad(Q)ad(P), \forall P,Q\in D$. Such a definition can be extended to any matrix of operators by using the transposed matrix of adjoint operators and we get:  
\[ <\lambda,{\cal{D}} \xi>=<ad({\cal{D}})\lambda,\xi>+\hspace{1mm} {div}\hspace{1mm} ( ... )  \]
from integration by part, where $\lambda$ is a row vector of test functions and $<  > $ the usual contraction. We quote the useful formulas $[ad(\xi),ad(\eta)]=ad(\xi)ad(\eta)-ad(\eta)ad(\xi)= - ad([\xi, \eta]), \forall \xi, \eta \in T$ ({\it care about the minus sign}) and $rk_D({\cal{D}})=rk_D(ad({\cal{D}}))$ as 
in ([23], p 610-612).\\

\noindent
{\bf REMARK  4.2}: As can be seen from the last two examples of Section 3, when ${\cal{D}}$ is involutive, then $ad({\cal{D}})$ may not be involutive. In the differential framework, we may set $ rk_D ({\cal{D}})=m-{\alpha}^n_q={\beta}^n_q$. Comparing to similar concepts used in {\it differential algebra}, this number is just the maximum number of differentially independent equations to be found in the differential module $I$ of equations. Indeed, pointing out that differential indeterminates in differential algebra are nothing else than jet coordinates in differential geometry and using standard notations, we have $K\{y\}=lim_{q\rightarrow \infty}K[y_q]$. In that case, the differential ideal $I$ {\it automatically} generates a prime differential ideal $\mathfrak{p}\subset K\{y\}$ providing a {\it differential extension} $L/K$ with $L=Q(K\{y\}/\mathfrak{p})$ and {\it differential transcendence degree} $diff\hspace{1mm}trd (L/K)={\alpha}^n_q$, a result explaining the notations ([12],[22]). Now, from the dimension formulas of $R_{q+r}$, we obtain at once $rk_D(M)={\alpha}^n_q$ and thus $rk_D({\cal{D}})=m - rk_D(M)$ in a coherent way with any free presentation of $M$ starting with ${\cal{D}}$. However, ${\cal{D}}$ acts on the left in differential geometry but on the right in the theory of differential modules. For an operator of order zero, we recognize the fact that the rank of a matrix is eqal to the rank of the transposed matrix.\\

\noindent
{\bf PROPOSITION  4.3}: If $f\in aut(X)$ is a local diffeomorphisms on $X$, we may set $ x=f^{-1}(y)=g(y)$ and we have the {\it identity}:
\[   \frac{\partial}{\partial y^k}(\frac{1}{\Delta (g(y))} {\partial}_if^k(g(y))) \equiv 0 \hspace{2mm} \Rightarrow \hspace{2mm} \frac{\partial}{\partial y^k}(\frac{1}{\Delta}\frac{\partial f^k}{\partial x^i}{\cal{A}}^i)=   \frac{1}{\Delta}\frac{\partial f^k}{\partial x^i}\frac{\partial {\cal{A}}^i}{\partial y^k}=\frac{1}{\Delta}{\partial}_i{\cal{A}}^i\]
and the adjoint of the well defined intrinsic operator ${\wedge}^0T^* \stackrel{d}{\longrightarrow} {\wedge}^1T^*=T^*:A \longrightarrow {\partial}_iA$ is ({\it minus}) the well defined intrinsic operator ${\wedge}^nT^* \stackrel{d}{\longleftarrow} {\wedge}^nT^*\otimes T \simeq {\wedge}^{n-1}T^*: {\partial}_i{\cal{A}}^i\longleftarrow {\cal{A}}^i$. Hence, if we have an operator $E\stackrel{\cal{D}}{\longrightarrow} F$, we obtain the {\it formal adjoint} operator ${\wedge}^nT^*\otimes E^*\stackrel{ad(\cal{D})}{\longleftarrow} {\wedge}^nT^*\otimes F^*$. \\

Having in mind that $D$ is a $K$-algebra, that $K$ is a left $D$-module with the standard action $(D,K) \longrightarrow K:(P,a) \longrightarrow P(a):(d_i,a)\longrightarrow {\partial}_ia$ and that $D$ is a bimodule over itself, {\it we have only two possible constructions leading to the following two definitions}:  \\

\noindent
{\bf DEFINITION  4.4}: We may define the {\it inverse system} $R=hom_K(M,K)$ of $M$ and introduce $R_q=hom_K(M_q,K)$ as the {\it inverse system of order} $q$. \\

\noindent
{\bf DEFINITION  4.5}: We may define the right differential module $M^*=hom_D(M,D)$ by using the bimodule structure of $D={ }_DD_D$.  \\

\noindent
{\bf THEOREM  4.6}: When $M$ and $N$ are left $D$-modules, then $hom_K(M,N)$ and $M{\otimes}_KN$ are left $D$-modules. In particular $R=hom_K(M,K)$ is also a left $D$-module for the {\it Spencer operator}.  \\

\noindent
{\it Proof}:  For any $f\in hom_K(M,N)$, let us define:   \\
\[   (af)(m)=af(m)=f(am) \hspace{1cm} \forall a\in K, \forall m\in M\]
\[   (\xi f)(m)=\xi f(m)-f(\xi m)  \hspace{1cm}  \forall \xi ={\xi}^id_i\in T, \forall m\in M  \]
It is easy to check that $\xi a=a \xi+\xi (a)$ in the operator sense and that $\xi\eta -\eta\xi =[\xi,\eta]$ is the standard bracket of vector fields. We have in particular with $d$ in place of any $d_i$: \\
\[  \begin{array}{rcl}
((da)f)(m)=(d(af))(m)=d(af(m))-af(dm)&=&(\partial a)f(m)+ad(f(m))-af(dm)\\
       &=& (a(df))(m)+(\partial a)f(m)  \\
       &=& ((ad+\partial a)f)(m)
       \end{array}  \]
 For any $m\otimes n\in M{\otimes}_KN$ with arbitrary $m\in M$ and $n\in N$, we may then define:   \\
 \[      a(m\otimes n)=am\otimes n=m\otimes an\in M{\otimes}_KN  \]
 \[  \xi (m\otimes n)=\xi m\otimes n + m\otimes \xi n \in M{\otimes}_KN   \]
 and conclude similarly with:   \\
 \[  \begin{array}{rcl}
  (da)(m\otimes n)=d(a(m\otimes n)) & = & d(am\otimes n)\\
                         & = &  d(am)\otimes n+am\otimes dn  \\
                             & = & (\partial a)m\otimes n + a(dm)\otimes n + am\otimes dn  \\
                                & = & (ad+\partial a)(m\otimes n)
                                \end{array}    \]
Using $K$ in place of $N$, we finally get $(d_if)^k_{\mu}=(d_if)(y^k_{\mu})={\partial}_if^k_{\mu}-f^k_{\mu +1_i}$ that is {\it we recognize exactly the Spencer operator} 
with {\it now} $Df=dx^i\otimes d_if$ and thus:\\
\[  (d_i(d_jf))^k_{\mu}={\partial}_{ij}f^k_{\mu}-({\partial}_if^k_{\mu+1_j}+{\partial}_jf^k_{\mu+1_i})+f^k_{\mu+1_i+1_j} \Rightarrow d_i(d_jf)=d_j(d_if)=d_{ij}f \]
In fact, $R$ is the {\it projective limit} of ${\pi}^{q+r}_q:R_{q+r}\rightarrow R_q$ in a coherent way with jet theory ([2],[27],[38]).\\
\hspace*{12cm}  Q.E.D.  \\

\noindent
{\bf COROLLARY  4.7}: If $M$ and $N$ are right $D$-modules, then $hom_K(M,N)$ is a left $D$-module. Moreover, if $M$ is a left $D$-module and $N$ is a right $D$-module, then $M{\otimes}_KN$ is a right $D$-module. \\

\noindent
{\it Proof}: If $M$ and $N$ are right $D$-modules, we just need to set $(\xi f)(m)=f(m\xi)-f(m)\xi, \forall \xi\in T, \forall m\in M$ and conclude as before. Similarly, if $M$ is a left $D$-module and $N$ is a right $D$-module, we just need to set $(m\otimes n)\xi=m\otimes n\xi - \xi m \otimes n$. \\
\hspace*{12cm}  Q.E.D.  \\
 
\noindent
{\bf REMARK  4.8}: When $M={Ê}_DM\in mod(D)$ and $N=N_D$, , then $hom_K(N,M)$ cannot be endowed with any left or right differential structure. When $M=M_D$ and $N=N_D$, then $M{\otimes}_KN$ cannot be endowed with any left or right differential structure (See [2], p 24 for more details).  \\

As $M^*=hom_D(M,D)$ is a right $D$-module, let us define the right $D$-module $N_D$ by the ker/coker long exact sequence $0\longleftarrow N_D \longleftarrow F^*_1 \stackrel{{\cal{D}}^*}{\longleftarrow} F^*_0 \longleftarrow M^* \longleftarrow 0$.\\

\noindent
{\bf COROLLARY  4.9}: We have the {\it side changing} procedure $N_D \rightarrow N={ }_DN=hom_K({\wedge}^nT^*,N_D)$ with inverse $M={Ê}_DM \rightarrow M_D={\wedge}^nT^*{\otimes}_K M$ whenever $M,N \in mod(D)$.  \\

\noindent
{\it Proof}: According to the above Theorem, we just need to prove that ${\wedge}^nT^*$ has a natural right module structure over $D$. For this, if $\alpha=adx^1\wedge ...\wedge dx^n\in T^*$ is a volume form with coefficient $a\in K$, we may set $\alpha.P=ad(P)(a)dx^1\wedge...\wedge dx^n$ when $P\in D$. As $D$ is generated by $K$ and $T$, we just need to check that the above formula has an intrinsic meaning for any $\xi={\xi}^id_i\in T$. In that case, we check at once:
\[  \alpha.\xi=-{\partial}_i(a{\xi}^i)dx^1\wedge...\wedge dx^n=-\cal{L}(\xi)\alpha \]
by introducing the Lie derivative of $\alpha$ with respect to $\xi$, along the intrinsic formula ${\cal{L}}(\xi)=i(\xi)d+di(\xi)$ where $i( )$ is the interior multiplication and $d$ is the exterior derivative of exterior forms. According to well known properties of the Lie derivative, we get :
\[\alpha.(a\xi)=(\alpha.\xi).a-\alpha.\xi(a), \hspace{5mm} \alpha.(\xi\eta-\eta\xi)=-[\cal{L}(\xi),\cal{L}(\eta)]\alpha=-\cal{L}([\xi,\eta])\alpha=\alpha.[\xi,\eta].  \]
\hspace*{12cm}  Q.E.D.  \\

Collecting the previous results, if a differential operator ${\cal{D}}$ is given in the framework of differential geometry, we may keep the same notation ${\cal{D}}$ in the framework of differential modules which are {\it left} modules over the ring $D$ of linear differential operators and apply duality, provided we use the notation ${\cal{D}}^*$ and deal with {\it right} differential modules or use the notation $ad({\cal{D}})$ and deal again with {\it left} differential modules by using the $left \leftrightarrow right$ {\it conversion} procedure.  \\

\noindent
{\bf DEFINITION  4.10}: If an operator $\xi \stackrel{\cal{D}}{\longrightarrow} \eta$ is given, a {\it direct problem} is to look for (generating) {\it compatibility conditions} (CC) as an operator $\eta \stackrel{{\cal{D}}_1}{\longrightarrow} \zeta $ such that ${\cal{D}}\xi=\eta \Rightarrow {\cal{D}}_1\eta=0$. Conversely, given $\eta \stackrel{{\cal{D}}_1}{\longrightarrow} \zeta$, the {\it inverse problem} will be to look for $\xi \stackrel{\cal{D}}{\longrightarrow} \eta$ such that ${\cal{D}}_1$ generates the CC of ${\cal{D}}$ and we shall say that ${\cal{D}}_1$ {\it is parametrized by} ${\cal{D}}$ {\it if such an operator} ${\cal{D}}$ {\it is existing}.  \\

As $ad(ad(P))=P, \forall P \in D$, any operator is the adjoint of a certain operator and we get:  \\

\noindent
{\bf DOUBLE DUALITY CRITERION  4.11}: An operator ${\cal{D}}_1$ can be parametrized by an operator ${\cal{D}}$ if, whenever $ad({\cal{D}})$ generates the CC of $ad({\cal{D}}_1)$, then ${\cal{D}}_1$ generates the CC of ${\cal{D}}$. However, as shown in the example below, many other parametrizations may exist.  \\

Reversing the arrows, we finally obtain:  \\

\noindent
{\bf TORSION-FREE CRITERION  4.12}: A differential module $M$ having a finite free presentation $F_1 \stackrel{{\cal{D}}_1}{\longrightarrow} F_0 \rightarrow M \rightarrow 0$ is {\it torsion-free}, that is to say $t(M)=\{ m\in M\mid \exists 0 \neq P\in D, Pm=0 \}=0$, if and only if there exists a free differential module $E$ and an exact sequence $F_1 \stackrel{{\cal{D}}_1}{\longrightarrow} F_0  \stackrel{{\cal{D}}}{\longrightarrow} E  $ providing the {\it parametrization} $M\subseteq E $. \\

\noindent
{\bf REMARK  4.13}: Of course, solving the direct problem (Janet, Spencer) is {\it necessary} for solving the inverse problem. However, though the direct problem always has a solution, the inverse problem may not have a solution at all and the case of the Einstein operator is one of the best non-trivial PD counterexamples ([24],[30]). It is rather striking to discover that, in the case of OD operators, it took almost 50 years to understand that the possibility to solve the inverse problem was equivalent to the controllability of the corresponding control system ([24],[34]).\\

\noindent
{\bf EXAMPLE  4.14}: ({\it contact transformations}) With $n=3,K=\mathbb{Q}(x^1,x^2,x^3)$, let us consider the Lie pseudogroup of transfomations preserving the first order geometric object $\omega$ like a $1$-form but up to the square root of $\Delta$. The infinitesimal transformations are among the solutions $\Theta$ of the {\it general} system:\\ 
\[  {\Omega}_i\equiv ({\cal{L}}(\xi)\omega)_i\equiv {\omega}_r(x){\partial}_i{\xi}^r-(1/2){\omega}_i(x){\partial}_r{\xi}^r+{\xi}^r{\partial}_r{\omega}_i(x)=0   \]
When $\omega=(1,-x^3,0)$, we obtain the {\it special} involutive system:ÊÊ\\
\[  {\partial}_3{\xi}^3+{\partial}_2{\xi}^2+2x^3{\partial}_1{\xi}^2-{\partial}_1{\xi}^1=0, {\partial}_3{\xi}^1-x^3{\partial}_3{\xi}^2=0, {\partial}_2{\xi}^1-x^3{\partial}_2{\xi}^2+x^3{\partial}_1{\xi}^1-(x^3)^2{\partial}_1{\xi}^2-{\xi}^3=0   \]
with  $2$ equations of class $3$, $1$ equation of class 2 and thus only $1$ first order CC for the second members coming from the linearization of the Vessiot structure equation:\\
\[  {\omega}_1({\partial}_2{\omega}_3-{\partial}_3{\omega}_2)+{\omega}_2({\partial}_3{\omega}_1-{\partial}_1{\omega}_3)+
{\omega}_3({\partial}_1{\omega}_2-{\partial}_2{\omega}_1)=c    \]
involving the only {\it structure constant} $c$. This system can be parametrized by a single potential $\theta$: \\
\[-x^3{\partial}_3\theta + \theta={\xi}^1, -{\partial}_3\theta={\xi}^2, {\partial}_2\theta-x^3{\partial}_1\theta={\xi}^3 \Rightarrow {\xi}^1-x^3{\xi}^2=\theta \]
and we have the formally exact differential sequence  $0 \rightarrow 1 \stackrel{{\cal{D}}_{-1}}{\longrightarrow} 3 \stackrel{{\cal{D}}}{\longrightarrow} 3 \stackrel{{\cal{D}}_1}{\longrightarrow} 1 \rightarrow 0 $. \\

However, we have yet not proved the most difficult result that could not be obtained without homological algebra and the next example will explain ths additional 
difficulty.  \\

\noindent
{\bf EXAMPLE  4.15}:  With ${\partial}_{22}\xi={\eta}^2, {\partial}_{12}\xi={\eta}^1$ for $\cal{D}$, we get  ${\partial}_1{\eta}^2-{\partial}_2{\eta}^1=\zeta$ for ${\cal{D}}_1$. Then $ad({\cal{D}}_1)$ is defined by ${\mu}^2=-{\partial}_1\lambda, {\mu}^1={\partial}_2\lambda$ while $ad(\cal{D})$ is defined by $\nu={\partial}_{12}{\mu}^1+{\partial}_{22}{\mu}^2$ but the CC of $ad({\cal{D}}_1)$ are generated by ${\nu}'={\partial}_1{\mu}^1+{\partial}_2{\mu}^2$. In the operator framework, we have the differential sequences:\\  
\[  \begin{array}{rcccccl}
       &   \xi & \stackrel{\cal{D}}{\longrightarrow} & \eta & \stackrel{{\cal{D}}_1}{\longrightarrow} & \zeta  & \rightarrow 0 \\
 0 \leftarrow & \nu& \stackrel{ad(\cal{D})}{\longleftarrow} & \mu & \stackrel{ad({\cal{D}}_1)}{\longleftarrow} & \lambda &
  \end{array}  \]
where the upper sequence is formally exact at $\eta$ but the lower sequence is not formally exact at $\mu$.  \\
Passing to the module framework, we obtain the sequences:  \\
\[  \begin{array}{rcccccccl}
 0 \rightarrow   &  D & \stackrel{{\cal{D}}_1}{\longrightarrow} & D^2 & \stackrel{\cal{D}}{\longrightarrow} & D & \rightarrow &M &  \rightarrow  0  \\
    &     D& \stackrel{ad({\cal{D}}_1)}{\longleftarrow} & D^2 & \stackrel{ad(\cal{D})}{\longleftarrow} & D & \leftarrow    &  0  &
  \end{array}  \]
where the lower sequence is not exact at $D^2$.  \\

Therefore, we have to find out situations in which $ad({\cal{D}})$ generates the CC of $ad({\cal{D}}_1)$ whenever ${\cal{D}}_1$ generates the CC of ${\cal{D}}$ and conversely. This problem will be studied in Section 5, Part C.  \\

\noindent
{\bf 5) APPLICATIONS}  \\

Though the next pages will only be concerned with a study of the Lie pseudogroups of isometries (A) and conformal isometries (B), the reader must never forget that they can be used similarly for {\it any} arbitrary transitive Lie pseudogroup of transformations ([19],[21],[22],[29]). \\

\noindent
{\bf A) RIEMANN TENSOR}\\

 Let $\omega=({\omega}_{ij}={\omega}_{ji})\in S_2T¬*$be a non-degenerate metric with $det(\omega)\neq 0$. We shall apply the Main Theorem to the first order {\it Killing system} $R_1\subset J_1(T)$ defined by the $n(n+1)/2$ linear equations ${\Omega}_{ij}\equiv {\omega}_{rj}{\xi}^r_i+{\omega}_{ir}{\xi}^r_j +{\xi}^r{\partial}_r{\omega}_{ij}=0$ for any section ${\xi}_1\in R_1$. Its symbol $g_1\subset T^*\otimes T$ is defined by the $n(n+1)/2$ linear equations ${\omega}_{rj}{\xi}^r_i+{\omega}_{ir}{\xi}^r_j=0$ and we obtain at once isomorphisms $g_1\simeq {\wedge}^2T\simeq {\wedge}^2T^*$ by lowering or raising the indices by means of the metric, obtaining for example ${\xi}_{i,j}+{\xi}_{j,i}=0$. As $det(\omega)\neq 0$, we may introduce the well known Chrisoffel symbols $\gamma=({\gamma}^k_{ij}={\gamma}^k_{ji})$ through the standard Ricci/Levi-Civita isomorphism $j_1(\omega)\simeq (\omega,\gamma)$ and obtain by one prolongation the linear second order equations for any section ${\xi}_2\in R_2$: \\
 \[   {\Gamma}^k_{ij} \equiv {\xi}^k_{ij} +{\gamma}^k_{rj}{\xi}^r_i + {\gamma}^k_{ir} - {\gamma}^r_{ij}{\xi}^k_r + {\xi}^r{\partial}_r{\gamma}^k_{ij} =0  \]
 and we have $\Omega \in S_2T^* \Rightarrow \Gamma \in S_2T^* \otimes T $ for the respective linearization/variation of $\omega$ and $\gamma$. As we shall see that $g_1$ is {\it not} $2$-acyclic and $g_2=0$ is defined by the $n^2(n+1)/2$ linear equations ${\xi}^k_{ij}=0$, we may apply the Main Theorem with $q=1,s=1,E=T$, on the condition that $R_1$ should be formally integrable as it is finite type and cannot therefore be involutive. First of all, we have the following commutative and exact diagram  allowing to define $F_0=S_2T^*$:  \\
 \[  \begin{array}{rcccccl}
    &  0  &  &  0  &  &  0  &  \\
    & \downarrow & & \downarrow & & \downarrow &  \\
    0  \rightarrow & g_1 & \longrightarrow &T^*\otimes T& \stackrel{\sigma(\Phi)}{\longrightarrow} & F_0 & \rightarrow 0 \\
    & \downarrow & & \downarrow & & \parallel &  \\
 0  \rightarrow & R_1 & \longrightarrow &J_1(T)& \stackrel{\Phi}{\longrightarrow} & F_0 & \rightarrow 0 \\
& \downarrow & & \downarrow & & \downarrow &  \\
  0  \rightarrow & T & = & T & \longrightarrow & 0 &  \\
  & \downarrow & & \downarrow & & &  \\ 
 &  0  &  &  0  &  &   &  
\end{array}   \]
Now, $R_2 \stackrel{{\pi}^2_1}{\longrightarrow} R_1$ is an isomorphism because $g_2=0$ and $dim(R_2)=dim(R_1)=n(n+1)/2$. Hence, $R_2\subset J_2(T)$ is involutive if and only if $R_3 \stackrel{{\pi}^3_2}{\longrightarrow} R_2$ is also an isomorphism too because $g_2=0 \Rightarrow g_{2+r}=0, \forall r\geq 0$. Such a differential condition for $\omega$ has been shown by L.P. Eisenhart in ([8]) to be equivalent to the {\it Vessiot structure equation with one constant} called {\it constant riemannian curvature} ${\rho}^k_{l,ij}=c({\delta}^k_i{\omega}_{lj} - {\delta}^k_j{\omega}_{li})$ (See [19],[22] and [29] for effective calculations still not acknowledged today). In this formula, $c$ is an arbitrary constant and the {\it Riemann tensor} $({\rho}^k_{l,ij})\in {\wedge}^2T^*\otimes T^*\otimes T$ satisfies the two types of purely algebraic relations: \\
\[  {\omega}_{rl}{\rho}^r_{k,ij} +{\omega}_{kr}{\rho}^r_{l,ij}=0, \hspace{5mm} {\rho}^k_{l,ij} +  {\rho}^k_{i,jl} +  {\rho}^k_{j,li}=0  \]
We shall suppose that $\omega$ is the Euclidean metric if $n=2,3$ and the Minkowskian metric if $n=4$ but any other compatible choice should be convenient. As a next step, we know from the Main Theorem that the generating CC for the operator $Killing=\Phi \circ j_1:T \rightarrow F_0$ are made by an operator $Riemann=\Psi \circ j_2:F_0 \rightarrow F_1$ of order $s+1=2$. We shall define $F_1$ by setting $q=1,r=0,s=1$ in the corresponding diagram in order to get the following commutative diagram:  \\
\[  \begin{array}{rcccccccl}
   &  0 & & 0 & & 0 &  &  &   \\
   & \downarrow & & \downarrow & & \downarrow & & &  \\
0\rightarrow & g_3 & \rightarrow &  S_3T^*\otimes T & \rightarrow & S_2T^*\otimes F_0& \rightarrow & F_1 & \rightarrow 0  \\
   & \hspace{2mm}\downarrow  \delta  & & \hspace{2mm}\downarrow \delta & &\hspace{2mm} \downarrow \delta & & &  \\
0\rightarrow& T^*\otimes g_2&\rightarrow &T^*\otimes S_2T^*\otimes T & \rightarrow &T^*\otimes T^*\otimes F_0 &\rightarrow & 0 &  \\
   &\hspace{2mm} \downarrow \delta &  &\hspace{2mm} \downarrow \delta & &\hspace{2mm}\downarrow \delta &  &  &   \\
0\rightarrow & {\wedge}^2T^*\otimes g_1 & \rightarrow & \underline{{\wedge}^2T^*\otimes T^*\otimes T} & \rightarrow & {\wedge}^2T^*\otimes F_0 & \rightarrow & 0 &  \\
   &\hspace{2mm}\downarrow \delta  &  & \hspace{2mm} \downarrow \delta  &  & \downarrow  & &  &  \\
0\rightarrow & {\wedge}^3T^*\otimes T & =  & {\wedge}^3T^*\otimes T  &\rightarrow   & 0  &  &  &   \\
    &  \downarrow  &  &  \downarrow  &  &  &  &  &  \\
    &  0  &   & 0  & &  &  &  &
\end{array}  \]
where all the rows are exact and all the columns are also exact but the first at ${\wedge}^2T^*\otimes g_1$ with $g_2=0 \Rightarrow g_3=0$. We shall denote by $B^2(g_1)$ the {\it coboundary} as the image of the central $\delta$, by $Z^2(g_1)$ the {\it cocycle} as the kernel of the lower $\delta$ and by $H^2(g_1)=Z^2(g_1)/B^2(g_1)$ the {\it Spencer} $\delta$-{\it  cohomology} at ${\wedge}^2T^*\otimes g_1$ as the quotient. Chasing in the previous diagram, we discover that the {\it Riemann tensor} is a section of the bundle $ F_1=H^2(g_1)=Z^2(g_1)$ with $dim(F_1)= (n^2(n+1)^2/4)-(n^2(n+1)(n+2)/6)=(n^2(n-1)^2/4)-(n^2(n-1)(n-2)/6)=n^2(n^2-1)/12$ by using the top row or the left column. We discover at once the two properties of the (linearized) Riemann tensor through the chase involved, namely $(R^k_{l,ij})\in {\wedge}^2T^*\otimes T^*\otimes T$ is killed by both $\delta$ and ${\sigma}_0(\Phi)$. Similarly, going one step further, we get the (linearized) Bianchi identities by defining the first order operator $Bianchi:F_1 \rightarrow F_2$ where $F_2=H^3(g_1)=Z^3(g_1)$ with $ dim(F_2)=dim({\wedge}^3T^*\otimes g_1)-dim({\wedge}^4T^*\otimes T)=n^2(n-1)^2(n-2)/12 - n^2(n-1)(n-2)(n-3)/24= n^2(n^2-1)(n-2)/24$ may be defined by the following commutative diagram:  \\
\[  \begin{array}{rcccccccccl}
   & 0  &  & 0  & &  0  &  & 0  &  &     \\
   & \downarrow  &  &  \downarrow & & \downarrow & & \downarrow & &  & \\
0 \rightarrow & g_4 & \rightarrow &S_4T^*\otimes T& \rightarrow &S_3T^ *\otimes F_0 &\rightarrow & T^*\otimes F_1&\rightarrow & F_2 & \rightarrow 0 \\
 & \downarrow  &  &  \downarrow & & \downarrow & & \parallel & &  \\
0 \rightarrow & T^*\otimes g_3 & \rightarrow &T^*\otimes S_3T^*\otimes T& \rightarrow &T^*\otimes S_2T^ *\otimes F_0 &\rightarrow & T^*\otimes F_1& \rightarrow & 0 &  \\
 & \downarrow  &  &  \downarrow & & \downarrow & & \downarrow & & &  \\
0 \rightarrow &{\wedge}^2 T^*\otimes g_2 & \rightarrow &{\wedge}^2T^*\otimes S_2T^*\otimes T& \rightarrow &{\wedge}^2T^*\otimes T^ *\otimes F_0 &\rightarrow &
 0&&& \\
& \downarrow  &  &  \downarrow & & \downarrow & &  & & & \\
0 \rightarrow &{\wedge}^3 T^*\otimes g_1 & \rightarrow &\underline{{\wedge}^3T^*\otimes T^*\otimes T}& \rightarrow &{\wedge}^3T^*\otimes F_0 &\rightarrow & 0 & &&  \\
  & \downarrow  &  &  \downarrow & & \downarrow & & & & &  \\
0 \rightarrow &{\wedge}^4 T^*\otimes T & = &{\wedge}^4T^*\otimes  T& \rightarrow &0 & & & &&  \\
 & \downarrow  &  &  \downarrow & &  & & & & & \\
 & 0  &  & 0  & &    &  &   &   &  & 
\end{array}  \]
This approach is relating for the first time the concept of {\it Riemann tensor candidate}, introduced by Lanczos and others, to the Spencer $\delta$-cohomology of the Killing symbols. We obtain therefore the formally exact sequence:  \\
\[  0 \rightarrow \Theta \rightarrow n \stackrel{Killing}{\longrightarrow} n(n+1)/2 \stackrel{Riemann}{\longrightarrow} n^2(n^2-1)/12\stackrel{Bianchi}{\longrightarrow} n^2(n^2-1)(n-2)/24  \rightarrow ...  \]
with operators of successive orders $1,2,1, ...$ and so on. \\
{\it In the present situation}, we have the (split) short exact sequences: \\
\[  0 \rightarrow F_1 \rightarrow {\wedge}^2T^*\otimes g_1 \stackrel{\delta}{\longrightarrow} {\wedge}^3T^*\otimes T \rightarrow 0, \hspace{5mm}
   0 \rightarrow F_2 \rightarrow {\wedge}^3T^*\otimes g_1 \stackrel{\delta}{\rightarrow} {\wedge}^4T^*\otimes T \rightarrow 0  \]
and obtain the operator $ad(Bianchi):{\wedge}^nT^*\otimes F_2^*\rightarrow {\wedge}^nT^*\otimes F_1^*$ with the short exact sequence: \\
\[  0 \leftarrow {\wedge}^nT^*\otimes F_2^* \leftarrow {\wedge}^{n-3}T^*\otimes g_1^* \leftarrow {\wedge}^{n-2}T^*\otimes T^* \leftarrow 0  \]
explaining at once why the {\it Lagrange multipliers} $\lambda \in {\wedge}^nT^*\otimes F_2^*$ can be represented by a section of ${\wedge}^{n-3}T^*\otimes {\wedge}^2T^*$, that is by a Lanczos potential in $T\otimes {\wedge}^2T^*$ when $n=4$. We shall see in part C that $ad(Bianchi)$ is parametrizing $ad(Riemann)$ contrary to the claims of Lanczos. Moreover, we have already pointed out in many books ( [21],[23]) or papers ([28],[32]) that continuum mechanics may be presented through a variational problem with a differential constraints which is {\it shifted by one step backwards in the previous differential sequence} because the infinitesimal deformation tensor $\epsilon=\frac{1}{2}\Omega \in S_2T^*$ must be now  killed by the operator $Riemann$ and the corresponding Lagrange multipliers $\lambda \in {\wedge}^nT^*\otimes F_1^*$ must be used because $ad(Riemann)$ is {\it } parametrizing $ad(Killing)=Cauchy$. Anybody using computations with finite elements also knows that a similar situation is held by electromagnetism too because the EM field is killed by $d:{\wedge}^2T^* \rightarrow {\wedge}^3T^* \Rightarrow ad(d):{\wedge}^1T^*\rightarrow {\wedge}^2T^*$, another fact contradicting Lanczos claims.\\
Finally, the passage to differential modules can be achieved easily by using $K=\mathbb{Q}$ as will be done in the Appendix or $K=\mathbb{Q}<\omega>$ with standard notations because the Lie pseudogroup of isometries is an {\it algebraic Lie pseudogroup} as it can be defined by differential polynomials in the jets of order 
$\geq 1$ (See [12],[20],[22] for details).  \\

\noindent
{\bf B) WEYL TENSOR}\\

If the study of the Riemann tensor/operator has been related to many classical results, the study of the Weyl tensor/operator in this new framework is quite different because these new mathematical tools have not been available before $1975$ and are still not acknowledged today by mathematical physicists. In particular, we may quote the link existing between acyclicity and formal integrability both with the possibility to use the {\it Vessiot structure equations} in order to combine in a unique framework the constant riemannian curvature condition needed for the Killing system, which only depends on one arbitrary constant, with the zero Weyl tensor condition needed for the conformal Killing system, which does not depend on any constant. For this reason, we shall follow as closely as possible the previous part $A$, putting a "$hat$" on the corresponding concepts.\\

 The {\it conformal Killing system} ${\hat{R}}_1\subset J_1(T) $ is defined by eliminating the function $A(x)$ in the system ${\cal{L}}(\xi)\omega=A(x)\omega$. It is also a {\it Lie operator} $\hat{\cal{D}}$ with solutions $\hat{\Theta}\subset T$ satisfying $[\hat{\Theta},\hat{\Theta}]\subset \hat{\Theta}$. Its symbol ${\hat{g}}_1$ is defined by the linear equations ${\omega}_{rj}{\xi}^r_i+{\omega}_{ir}{\xi}^r_j - \frac{2}{n}{\omega}_{ij}{\xi}^r_r=0$ which do not depend on any conformal factor and is finite type because ${\hat{g}}_3=0$ when $n\geq 3$. We have ([19],[20],[32]): \\

\noindent
{\bf LEMMA  5.1}: ${\hat{g}}_2 \subset S_2T^*\otimes T$ is {\it now} $2$-acyclic {\it only when} $n\geq 4$ and $3$-acyclic {\it only when} $n\geq 5$.  \\

 It is known that ${\hat{R}}_2$ and thus ${\hat{R}}_1$ too (by a chase) are formally integrable if and only if $\omega$ has zero {\it Weyl tensor}:  \\
 \[  {\tau}^k_{l,ij}\equiv {\rho}^k_{l,ij} - \frac{1}{(n-2)}({\delta}^k_i{\rho}_{lj} - {\delta}^k_j{\rho}_{li} +{\omega}_{lj}{\omega}^{ks}{\rho}_{si} - {\omega}_{li}{\omega}^{ks}{\rho}_{sj}) + \frac{1}{(n-1)(n-2)}({\delta}^k_i{\omega}_{lj} - {\delta}^k_j{\omega}_{li})\rho=0  \]
If we use the formula $id_M-f\circ u=v\circ g$ of Proposition 2.4 in the {\it split short exact sequence} induced by the inclusions $g_1\subset {\hat{g}}_1, 0=g_2\subset {\hat{g}}_2, g_3={\hat{g}}_3=0$ ([21],[22],[28]):  \\
\[ 0 \longrightarrow Ricci \longrightarrow Riemann \longrightarrow Weyl \longrightarrow  0  \]
according to the Vessiot structure equations, in particular if $\omega$ has constant Riemannian curvature and thus ${\rho}_{ij}={\rho}^r_{i,rj}=c(n-1){\omega}_{ij} \Rightarrow \rho={\omega}^{ij}{\rho}_{ij}=cn(n-1)$ ([19],[21],[30],[31]). Using the same diagrams as before, we get ${\hat{F}}_0=T^*\otimes T/{\hat{g}}_1$ with $dim({\hat{F}}_0)=(n-1)(n+2)/2$ and ${\hat{F}}_1=H^2({\hat{g}}_1)\neq Z^2({\hat{g}}_1)$ for defining any {\it Weyl tensor candidate}. As a byproduct, {\it we could believe} that the linearized operator $Weyl:{\hat{F}}_0 \rightarrow {\hat{F}}_1$ is of order $2$ with a symbol ${\hat{h}}_2\subset S_2T^*\otimes {\hat{F}}_0$ which is {\it not} $2$-acyclic by applying the $\delta$-map to the short exact sequence:  \\
 \[  0 \rightarrow {\hat{g}}_{3+r} \longrightarrow S_{3+r}T^*\otimes T \stackrel{{\sigma}_{2+r}(\Phi)}{\longrightarrow} {\hat{h}}_{2+r} \rightarrow  0  \]
and chasing through the commutative diagram thus obtained with $r=0,1,2$. As ${\hat{h}}_3$ becomes $2$-acyclic after one prolongation of ${\hat{h}}_2$ only, it follows that {\it the generating CC for the Weyl operator are of order} $2$, a result showing that the so-called Bianchi identities for the Weyl tensor are {\it not} CC in the strict sense of the definition as they do not involve only the Weyl tensor. \\

{\it In fact, things are quite different and we have to distinguish three different cases}: \\

\noindent
$\bullet$ $n=3$: According to the last Lemma, ${\hat{g}}_2$ is {\it not} $2$-acyclic but ${\hat{g}}_3=0$ becomes trivially $2$-ayclic and even involutive, that is $s=2$. According to the Main Theorem, the operator $Weyl:{\hat{F}}_0 \rightarrow {\hat{F}}_1$ is {\it third order} because $s+1=3$ (See Appendix) and ${\hat{F}}_1$ is defined by the short exact sequences: \\
\[ 0\rightarrow S_4T^*\otimes T \rightarrow S_3T^*\otimes {\hat{F}}_0 \rightarrow {\hat{F}}_1 \rightarrow 0, \hspace{5mm} 0 \rightarrow {\hat{F}}_1 \rightarrow {\wedge}^2T^*\otimes {\hat{g}}_2 \stackrel{\delta}{\longrightarrow} {\wedge}^3T^*\otimes {\hat{g}}_1 \rightarrow 0 \]
with $dim({\hat{F}}_1)=50-45=9-4=5$. As {\it now} ${\hat{h}}_3\subset S_3T^*\otimes {\hat{F}}_0$, applying the $\delta$-map to the short exact sequence: \\
\[  0 \rightarrow {\hat{g}}_6 \rightarrow S_6T^* \otimes T \rightarrow {\hat{h}}_5 \rightarrow 0\]
and chasing, we discover that ${\hat{h}}_3$ is $2$-acyclic because ${\hat{g}}_3=0$. Accordingly, the operator $Bianchi: {\hat{F}}_1 \rightarrow {\hat{F}}_2$ is {\it first order} and ${\hat{F}}_2$ is defined by the long exact sequence:\\
\[  0 \rightarrow S_5T^*\otimes T \rightarrow S_4T^*\otimes {\hat{F}}_0 \rightarrow T^*\otimes {\hat{F}}_1 \rightarrow {\hat{F}}_2 \rightarrow 0  \]
or by the isomorphism $0 \rightarrow {\hat{F}}_2 \rightarrow {\wedge}^3T^*\otimes {\hat{g}}_2 \rightarrow 0$ giving $dim({\hat{F}}_2)=63-75+15=1\times 3=3$.\\
Recapitulating, when $n=3$ we have the formally exact differential sequence with $3-5+5-3=0$:\\
\[ 0 \rightarrow {\hat{\Theta}} \rightarrow 3 \stackrel{CKilling}{\longrightarrow} 5 \stackrel{Weyl}{\longrightarrow} 5 \stackrel{Bianchi}{\longrightarrow} 3 \rightarrow 0  \]
In actual practice, introducing the new geometric objects ${\hat{\gamma}}^k_{ij}={\gamma}^k_{ij}-\frac{1}{n}({\delta}^k_i{\gamma}^r_{rj}+{\delta}^k_j{\gamma}^r_{ri}-{\omega}_{ij}{\omega}^{ks}{\gamma}^r_{rs})$, linearizing and using the cyclic summation ${\cal{C}}(1,2,3)$, we get for example:  \\
\[  {\cal{C}}(1,2,3)[d_{23}({\hat{\Gamma}}^2_{12}+{\hat{\Gamma}}^1_{33}-{\hat{\Gamma}}^1_{11} - {\hat{\Gamma}}^1_{22})] =0  \]

\noindent
$\bullet$ $n=4$: {\it This situation is even more striking} because ${\hat{g}}_2$ is $2$ acyclic but {\it not} $3$-acyclic and thus $s=1$. As before, we have $dim({\hat{F}}_0)=(n-1)(n+2)/2=9$ but, according to the Main Theorem, the operator $Weyl:{\hat{F}}_0 \rightarrow {\hat{F}}_1$ is of order $s+1=2$ and ${\hat{F}}_1$ is defined by the short exact sequence:\\
\[  0 \rightarrow S_3T^*\otimes T \rightarrow S_2T^*\otimes {\hat{F}}_0 \rightarrow {\hat{F}}_1 \rightarrow 0 \hspace{5mm}\Rightarrow \hspace{5mm} 
dim({\hat{F}}_1)=90-80=10\]
or by ${\hat{F}}_1=H^2({\hat{g}}_1)=Z^2({\hat{g}}_1)/B^2({\hat{g}}_1)$ with exact sequences:\\
\[ 0 \rightarrow T^*\otimes {\hat{g}}_2 \stackrel{\delta}{\longrightarrow} B^2({\hat{g}}_1) \rightarrow 0, \hspace{5mm} 0Ê\rightarrow Z^2({\hat{g}}_1) \rightarrow {\wedge}^2T^*\otimes {\hat{g}}_1 \stackrel{\delta}{\longrightarrow} {\wedge}^3T^*\otimes T \rightarrow 0  \]
 providing again $dim({\hat{F}}_1)=26-16=10$. The main problem is that,{\it now}, ${\hat{g}}_2$ is {\it not} $3$-acyclic and thus ${\hat{h}}_2$ is {\it not} $2$-acyclic according to a chase in the commutative diagram:  \\
\[ \begin{array}{rcccccl}
   & 0                    &  & 0                      &  & 0                       & \\
  & \downarrow &  &  \downarrow  &  &  \downarrow  &  \\
0    \rightarrow & {\hat{g}}_5 & \rightarrow & S_5T^*\otimes T & \rightarrow & {\hat{h}}_4 &  \rightarrow 0  \\  
   & \downarrow &  &  \downarrow  &    & \downarrow &  \\                                                                                   
0    \rightarrow & T^*\otimes {\hat{g}}_4 & \rightarrow & T^*\otimes S_4T^*\otimes T & \rightarrow & T^*\otimes {\hat{h}}_3 & \rightarrow 0 \\  
   & \downarrow &  &  \downarrow  &  & \downarrow &  \\
0    \rightarrow & {\wedge}^2T^*\otimes {\hat{g}}_3 & \rightarrow & {\wedge}^2T^*\otimes S_3T^*\otimes T&\rightarrow&
 {\wedge}^2T^*\otimes {\hat{h}}_2 &\rightarrow 0   \\  
& \downarrow &  &  \downarrow  &  &  \downarrow &     \\
0    \rightarrow & {\wedge}^3T^*\otimes {\hat{g}}_2& \rightarrow & {\wedge}^3T^*\otimes S_2T^*\otimes T & \rightarrow & {\wedge}^3T^*\otimes T^*\otimes {\hat{F}}_0&
\rightarrow 0  \\
  & \downarrow &  &  \downarrow  &  & \downarrow &  \\
0    \rightarrow & {\wedge}^4T^*\otimes {\hat{g}}_1 & \rightarrow & {\wedge}^4T^*\otimes T^*\otimes T&\rightarrow&
 {\wedge}^4T^* \otimes {\hat{F}}_0 &\rightarrow 0   \\  
  & \downarrow &  &  \downarrow  &  &  \downarrow  &  \\
   & 0                    &  & 0                      &  & 0                       & 
\end{array}  \]
but ${\hat{h}}_3$ {\it becomes} $2$-acyclic by chasing in the next diagram:  \\
\[ \begin{array}{rcccccl}
   & 0                    &  & 0                      &  & 0                       & \\
  & \downarrow &  &  \downarrow  &  &  \downarrow  &  \\
0    \rightarrow & {\hat{g}}_6 & \rightarrow & S_6T^*\otimes T & \rightarrow & {\hat{h}}_5 &  \rightarrow 0  \\  
   & \downarrow &  &  \downarrow  &    & \downarrow &  \\                                                                                   
0    \rightarrow & T^*\otimes {\hat{g}}_5 & \rightarrow & T^*\otimes S_5T^*\otimes T & \rightarrow & T^*\otimes {\hat{h}}_4 & \rightarrow 0 \\  
   & \downarrow &  &  \downarrow  &  & \downarrow &  \\
0    \rightarrow & {\wedge}^2T^*\otimes {\hat{g}}_4 & \rightarrow & {\wedge}^2T^*\otimes S_4T^*\otimes T&\rightarrow&
 {\wedge}^2T^*\otimes {\hat{h}}_3 &\rightarrow 0   \\  
& \downarrow &  &  \downarrow  &  &  \downarrow &     \\
0    \rightarrow & {\wedge}^3T^*\otimes {\hat{g}}_3& \rightarrow & {\wedge}^3T^*\otimes S_3T^*\otimes T & \rightarrow & {\wedge}^3T^*\otimes {\hat{h}}_2 &
\rightarrow 0  \\
  & \downarrow &  &  \downarrow  &  & \downarrow &  \\
0    \rightarrow & {\wedge}^4T^*\otimes {\hat{g}}_2 & \rightarrow & {\wedge}^4T^*\otimes S_2T^*\otimes T&\rightarrow&
 {\wedge}^4T^* \otimes T^* \otimes {\hat{F}}_0 &\rightarrow 0   \\  
  & \downarrow &  &  \downarrow  &  &  \downarrow  &  \\
   & 0                    &  & 0                      &  & 0                       & 
\end{array}  \]
and we have $s'=1$. Accordingly, the operator $Bianchi: {\hat{F}}_1 \rightarrow {\hat{F}}_2$ is of order $s'+1=2$ and ${\hat{F}}_2$ is defined by following commutative diagram where all the rows are exact and all the columns are exact but the first: Ê\\

\[   \begin{array}{rcccccccl}
   & 0  &  & 0  & &  0  &  & 0  &       \\
   & \downarrow  &  &  \downarrow & & \downarrow & & \downarrow &   \\
0 \rightarrow & {\hat{g}}_5 & \rightarrow &S_5T^*\otimes T& \rightarrow &S_4T^ *\otimes {\hat{F}}_0 &\rightarrow & S_2T^*\otimes {\hat{F}}_1& \rightarrow {\hat{F}}_2 \rightarrow 0 \\
 & \downarrow  &  &  \downarrow & & \downarrow & & \downarrow &    \\
0 \rightarrow & T^*\otimes {\hat{g}}_4 & \rightarrow &T^*\otimes S_4T^*\otimes T& \rightarrow &T^*\otimes S_3T^ *\otimes {\hat{F}}_0 &\rightarrow & T^*\otimes T^*\otimes {\hat{F}}_1& \rightarrow \hspace{1mm} 0   \\
 & \downarrow  &  &  \downarrow & & \downarrow & & \downarrow &    \\
0 \rightarrow &{\wedge}^2 T^*\otimes {\hat{g}}_3 & \rightarrow &{\wedge}^2T^*\otimes S_3T^*\otimes T& \rightarrow &{\wedge}^2T^*\otimes S_2T^ *\otimes {\hat{F}}_0 &\rightarrow & {\wedge}^2T^*\otimes {\hat{F}}_1& \rightarrow  \hspace{1mm} 0   \\
& \downarrow  &  &  \downarrow & & \downarrow & & \downarrow &   \\
0 \rightarrow &{\wedge}^3 T^*\otimes {\hat{g}}_2 & \rightarrow &{\wedge}^3T^*\otimes S_2T^*\otimes T& \rightarrow &{\wedge}^3T^*\otimes T^ *\otimes {\hat{F}}_0 &\rightarrow &0& \\
  & \downarrow  &  &  \downarrow & & \downarrow & & &    \\
0 \rightarrow &{\wedge}^4 T^*\otimes {\hat{g}}_1 & \rightarrow &{\wedge}^4T^*\otimes T^*\otimes T& \rightarrow &{\wedge}^4T^*\otimes {\hat{F}}_0 &\rightarrow & 0 &  \\
 & \downarrow  &  &  \downarrow & & \downarrow & & &    \\
 & 0  &  & 0  & &  0  &  &   &        
\end{array}  \]

\[\begin{array}{rcccccccccl}
0 \rightarrow & {\hat{g}}_5 & \rightarrow &S_5T^*\otimes T& \rightarrow &S_4T^ *\otimes {\hat{F}}_0 &\rightarrow & S_2T^*\otimes {\hat{F}}_1& \rightarrow &{\hat{F}}_2& \rightarrow 0 \\
                       &      0 & \rightarrow & 224                      & \rightarrow & 315                            & \rightarrow & 100                           & \rightarrow & 9  & \rightarrow 0  
     \end{array}   \]

\[  \begin{array}{rcccccccl}
0 \rightarrow &{\hat{g}}_4 & \rightarrow &S_4T^*\otimes T& \rightarrow &S_3T^ *\otimes {\hat{F}}_0 &\rightarrow & T^*\otimes {\hat{F}}_1& \rightarrow 0 \\
                        &   0    & \rightarrow & 140                      & \rightarrow &   180                         & \rightarrow  &      40  &                \rightarrow 0  
                        \end{array}   \]
We could define similarly the {\it first order} generating CC ${\hat{F}}_2 \rightarrow {\hat{F}}_3$ where ${\hat{F}}_3$ is defined by the following long exact sequence: \\
\[ \begin{array}{rccccccccccccl}
 0 \rightarrow& {\hat{g}}_6& \rightarrow &S_6T^*\otimes T& \rightarrow &S_5T^*\otimes {\hat{F}}_0& \rightarrow &S_3T^*\otimes {\hat{F}}_1& \rightarrow &T^*\otimes {\hat{F}}_2 &\rightarrow &{\hat{F}}_3 &\rightarrow 0 \\
 &0& \rightarrow &336 &\rightarrow & 504 &\rightarrow & 200 &\rightarrow & 36 & \rightarrow & 4 & \rightarrow 0  
\end{array}   \]
and obtain finally the formally exact differential sequence:\\
\[   0 \rightarrow \hat{\Theta} \rightarrow 4 \stackrel{CKilling}{\longrightarrow} 9 \stackrel{Weyl}{\longrightarrow} 10 \stackrel{Bianchi}{\longrightarrow} 9 \longrightarrow 4 \rightarrow 0  \]
with vanishing Euler-Poincar\'{e} characteristic $4-9+10-9+4=0$. We conclude the study of $n=4$ by exhibiting the short exact sequence:  \\ 
\[  \begin{array}{rcccccl}
0 \rightarrow & {\wedge}^2T^*\otimes {\hat{g}}_2 & \rightarrow & {\wedge}^3T^*\otimes {\hat{g}}_1 & \rightarrow & {\wedge}^4T^* \otimes T & \rightarrow 0  \\
0 \rightarrow  &                        24                    & \rightarrow &                 28                           & \rightarrow  &                       4                    & \rightarrow 0   
\end{array}   \]
a result showing that $H^3({\hat{g}}_1)=0$ for $n=4$, {\it contrary to what will happen when} $n\geq 5$. \\

\noindent
$\bullet$ $n\geq 5$: We still have $dim({\hat{F}}_0)=(n-1)(n+2)/2$ but now ${\hat{g}}_2$ is $2$-acyclic {\it and} $3$-acyclic with $s=1$ again but {\it now} with $s'=0$. Hence, according to the Main Theorem and its Corollary, the operator $Weyl$ is {\it second order} while the operator $Bianchi$ is {\it first order}. We may therefore use the same diagrams already introduced in part A but with a "{\it hat} " symbol and $n\geq 4$. In particular we get:  \\
\[  dim(Z^3({\hat{g}}_1))=dim ({\wedge}^3T^*\otimes {\hat{g}}_1) - dim ({\wedge}^4T^*\otimes T)=n(n-1)(n-2)(n^2+n+4)/24 \]
\[  dim(B^3({\hat{g}}_1))=n^2(n-1)/2  \hspace{5mm} \Rightarrow \hspace{5mm} dim(H^3({\hat{g}}_1))=n(n^2-1)(n+2)(n-4)/24  \]
We find again $H^3({\hat{g}}_1)=0$ when $n=4$ but $H^3({\hat{g}}_1)\neq 0$ when $n\geq 5$, a key step that no classical technique can even imagine. We have therefore a {\it first order} operator $Bianchi:{\hat{F}}_1 \rightarrow {\hat{F}}_2$ with ${\hat{F}}_1=H^2({\hat{g}}_1)$ and ${\hat{F}}_2=H^3({\hat{g}}_1)$ when $n\geq 5$. We obtain therefore a formally exact differential sequence: \\
\[  0 \rightarrow \hat{\Theta} \rightarrow n \stackrel{CKilling}{\longrightarrow} (n-1)(n+2)/2 \stackrel{Weyl}{\longrightarrow} n(n+1)(n+2)(n-3)/12 \stackrel{Bianchi}{\longrightarrow} n(n^2-1)(n+2)(n-4)/24  \]
In order to convince the reader about the powerfulness of the previous methods, let us prove that the generating CC of $Bianchi$ is a second order operator with $14$ equations when $n=5$. First, we ask the reader to prove, as an exercise, that $H^4({\hat{g}}_1)=0$ through the exact $\delta$-sequence $0 \rightarrow {\wedge}^3T^*\otimes {\hat{g}}_2\rightarrow {\wedge}^4T^*\otimes {\hat{g}}_1 \rightarrow {\wedge}^5T^*\otimes T \rightarrow 0$ (Hint: $50=55-5$). Conclude that the generating CC cannot be of first order. Then, prove that $dim(H^4({\hat{g}}_2))=25-11=14$. Finally, prove that the generating CC of these CC is first order with $5$ equations (Hint: Use the vanishing of the Euler-Poincar\'{e} characteristic with $5-14+35-35+14-5=0$) (See Appendix for confirmation).  \\

\noindent
{\bf C) PARAMETRIZATIONS}\\

Let $E$ and $F$ be two vector bundles with respective fiber dimensions $dim(E)=m$ and $dim(F)=p$. Starting with a differential operator ${\cal{D}}:E \rightarrow F$ of order $q$ with solutions $\Theta \subset E$ and such that the corresponding system $R_q\subset J_q(E)$ is formally integrable, we have explained in Section 3 how to construct a formally exact differential sequence: \\
\[ 0 \rightarrow \Theta \rightarrow E \stackrel{\cal{D}}{\longrightarrow} F_0 \stackrel{{\cal{D}}_1}{\longrightarrow} F_1 \stackrel{{\cal{D}}_2}{\longrightarrow}F_2 \stackrel{{\cal{D}}_3}{\longrightarrow} ...   \]
where $F_0=F$ and each operator ${\cal{D}}_i$ of order $q_i$ generates the CC of the previous one. In particular, if the starting operator ${\cal{D}}$ is involutive, then $q_1=...=q_n=1$ and each ${\cal{D}}_i$ is involutive in the resulting Janet sequence finishing at $F_n$ in the sense that ${\cal{D}}_n$ is formally surjective. Equivalently, it is possible to pass to the framework of differential module and look for a free resolution of a differential module $M$ starting with a free finite presentation $D^p \stackrel{\cal{D}}{\longrightarrow} D^m \rightarrow M \rightarrow 0$ with an operator acting on the right. In general, we have proved in the last two parts A and B that the succession of the orders $q,q_1,q_2, q_3,... $ can be nevertheless quite strange. Meanwhile, we have proved through examples that many possible finite length such sequences can be exhibited and the purpose of {\it Homological Algebra} is to study formal properties that should not depend on the sequence used. At the end od Section 4, we have pointed out the fact that, whenever an operator ${\cal{D}}_1$ generates all the CC of an operator $\cal{D}$, this does not imply in general that the operator $ad(\cal{D})$ generates all the CC of $ad({\cal{D}}_1)$. The following (quite difficult) theorem is a main result of homological algebra, adapted to differential systems and differential modules ([2],[3],[9],[17],[23],[37]):  \\

\noindent
{\bf THEOREM  5.2}: The fact that a formally exact differential sequence considered as a resolution of $\Theta$ has the property that the adjoint sequence is also formally exact does not depend on the sequence but only on $\Theta$. \\

\noindent
{\bf COROLLARY  5.3}: The fact that a free resolution of a differential module $M$ has the property that the adjoint sequence is also a free resolution does not depend on the sequence but only on $M$.  \\

Our problem in this part C is to describe a sufficiently general situation in such a way that {\it all} the results of the parts A and B can fit together, in the sense that we shall no longer need to use a $hat$ in order to distinguish them. For this, with $E=T$, let us say that an operator ${\cal{D}}=\Phi \circ j_q:T \rightarrow F$ is a {\it Lie operator} if ${\cal{D}} \xi=0, {\cal{D}}\eta=0 \Rightarrow {\cal{D}} [\xi,\eta]=0$ that is to say $[\Theta,\Theta]\subseteq \Theta$. The corresponding system $R_q=ker(\Phi)$ is called a {\it system of infinitesimal Lie equations} and one can define a "{\it bracket on sections}" satisfying $[R_q,R_q] \subseteq R_q$ in order to check formally the previous definition ([19],[29],[30]). It has been found by E. Vessiot, as early as in ... 1903 ([19],[40]), that the condition of formal integrability of $R_q$ can be described by the {\it Vessiot structure equations}, a set of (non-linear in general) differential conditions depending on a certain number of constants, for one or a family of {\it geometric objects} that can be vectors, forms, tensors or even higher order objects. The idea has been to look for "{\it general}" systems or symbols having the same dimensions as for a model object called "{\it special}", for example the euclidean metric when $n=2,3$ or the minkowskian metric when $n=4$. The case of part A has been a metric $\omega$ with $det(\omega)\neq 0$ and constant riemannian curvature with one constant while the case of part B has been a metric density $\hat{\omega}=\omega/(\mid det(\omega)\mid ^{\frac{1}{n}})$ with zero Weyl tensor and no constant involved. The following results will be local.\\

Let us suppose that we have a Lie group of transformations of $X$, namely a Lie group $G$ and an {\it action} $X \times G \rightarrow X:(x,a)\rightarrow y=ax=f(x,a)$ or, better, its graph $X\times G \rightarrow X\times X:(x,a) \rightarrow (x,y=ax=f(x,a)$. Differentiating enough times, we may eliminate the parameters $a$ among the equations $y_q=j_q(f)(x,a)$ for $q$ large enough and get a (non-linear in general) {\it system of finite Lie equations}. Linearising this system for $a$ close to the identity $e\in G$, that is for $y$ close to $x$, provides the system $R_q\subset J_q(T)$ and the corresponding Lie operator of finite type. Equivalently, the three theorems of Sophus Lie assert that there exists a {\it finite number of infinitesimal generators} $\{{\theta}_{\tau}\}$ of the action that should be linearly independent over the constants and satisfy $[{\theta}_{\rho},{\theta}_{\sigma}]=c^{\tau}_{\rho \sigma} {\theta}_{\tau}$ where the {\it structure constants} $c$ define a Lie algebra ${\cal{G}}=T_e(G)$. We have therefore $\xi \in \Theta \Leftrightarrow \xi = {\lambda}^{\tau}{\theta}_{\tau}$ with ${\lambda}^{\tau}=cst$. Hence, we may replace locally the system of infinitesimal Lie equations by the system ${\partial}_i{\lambda}^{\tau}=0$, getting therefore the differential sequence:  \\
\[  0 \rightarrow \Theta \rightarrow {\wedge}^0T^*\otimes {\cal{G}} \stackrel{d}{\longrightarrow} {\wedge}^1T^*\otimes {\cal{G}} \stackrel{d}{\longrightarrow} ... \stackrel{d}{\longrightarrow} {\wedge}^nT^* \otimes {\cal{G}} \rightarrow 0  \]
 which is the tensor product of the Poincar\'{e} sequence by ${\cal{G}}$. Finally, we are in a position to apply the previous Theorem and Corollary because the Poincar\'{e} sequence is self adjoint (up to sign), that is $ad(d)$ generates the CC of $ad(d)$ at any position, exactly like $d$ generates the CC of $d$ at any position. We invite the reader to compare with the situation of the Maxwell equations in electromagnetisme. However, we have proved in ([21],[22],[30],[31],[32]) why neither the Janet sequence nor the Poincar\'{e} sequence can be used in physics and must be replaced by the {\it Spencer sequence} which is another resolution of $\Theta$. We provide a few additional details on the motivations for such a procedure.\\
 
 For this, if $q$ is large enough in such a way that $g_q=0$ and thus $R_{q+1}\simeq R_q$, let us define locally a section ${\xi}_{q+1}\in R_{q+1}$ by the formula ${\xi}^k_{\mu+1_i}={\lambda}^{\tau}(x){\partial}_{\mu+1_i}{\theta}^k_{\tau}(x)$ and apply the Spencer operator $D$. We obtain at once ${(D{\xi}_{q+1})}^k_{\mu,i}= {\partial}_i{\xi}^k_{\mu}-{\xi}^k_{\mu+1_i}={\partial}_i{\lambda}^{\tau}{\partial}_{\mu}{\theta}^k_{\tau}$, a result proving that the previous sequence is (locally) isomorphic to the Spencer sequence:  \\
 \[  0 \rightarrow \Theta \rightarrow {\wedge}^0T^*\otimes R_q \stackrel{D}{\longrightarrow} {\wedge}^1T^*\otimes R_q \stackrel{D}{\longrightarrow} ... \stackrel{D}{\longrightarrow} {\wedge}^nT^* \otimes R_q \rightarrow 0  \]
In the present paper, we had $dim({\cal{G}})=n(n+1)/2$ in part A and $dim(\hat{\cal{G}})=(n+1)(n+2)/2$ in part B. Moreover, whatever is the part concerned, $ad({\cal{D}})$ generates the CC of $ad({\cal{D}}_1)$  because ${\cal{D}}_1$ generates the CC of ${\cal{D}}$ while, similarly, $ad({\cal{D}}_1)$ generates the CC of $ad({\cal{D}}_2)$  because ${\cal{D}}_2$ generates the CC of ${\cal{D}}_1$ and so on. We conclude with the following comments.  \\

1) Coming back to the first differential sequence constructed in part A, the Riemann tensor is a section of $F_1$ which is in the image of the operator $Riemann$ or in the kernel of the operator $Bianchi$. Indeed, Lanczos has been considering the Riemann tensor as a section of $F_1$ killed by the operator $Bianchi$ considered as a differential constraint. Accordingly, he has used the action of $ad(Bianchi)$ on the corresponding Lagrange multipliers. However, this operator parametrizes $ad(Riemann)$ as we saw and {\it cannot be used in order to parametrize the Riemann tensor by means of the Lanczos potential}. \\

2) Coming back to the second differential sequence constructed (independently) in part B, the Weyl tensor is a section of ${\hat{F}}_1$ which is in the image of the operator $Weyl$ or in the kernel of the operator $Bianchi$. As most of the results presented are unknown, in particular the fact that both operators $Weyl$ and $Bianchi$ are second order in dimension $n=4$, we believe that even the proper concept of a Weyl tensor candidate must be revisited within this new framework.  \\ 

3) In parts A and B, only linear differential operators have been used. However, it is known from {\it the formal theory of Lie pseudogroups} that non-linear differential sequences can be similarly constructed ([13],[19],[22]). {\it As a matter of fact}, if non-linear analogues of ${\cal{D}}$ and ${\cal{D}}_1$ may be exhibited, {\it this is not possible for} ${\cal{D}}_2$. Moreover, the only important problem is to compare the image of ${\cal{D}}$ with the kernel of ${\cal{D}}_1$ in the finite/infinitesimal {\it equivalence problem} (See [19], p 333 for a nice counterexample). We believe that {\it this shift by one step backwards in the interpretation of a differential sequence} will become important for future physics. It is commonly done in the variational calculus using finite elements where the free energy brings together the deformation tensor and the EM field on equal footing, {\it quite contrary to the ideas of Lanczos}.  \\

\noindent
{\bf 6) CONCLUSION}\\

In most textbooks, the Weyl tensor is always presented today by comparison with the Riemann tensor after eliminating a conformal factor and its derivatives. We have exhibited new methods in order to introduce both the Riemann and the Weyl tensor {\it independently} by using the formal theory of systems of partial differential equations (Spencer cohomology) in the study of arbitrary Lie pseudogroups while using the {\it Vessiot structure equations} for the Killing and conformal Killing systems {\it separately}. In particular, we have revisited, in both cases, the proper concept of Bianchi identities by means of homological algebra and diagram chasing, obtaining explicit numbers and orders for each dimension. These striking results are confirmed by means of computer algebra in the Appendix. They prove that the work of Lanczos and followers must be revisited within this new framework. \\

\noindent
{\bf REFERENCES FOR THE PAPER} \\

\noindent
[1] G.B. AIRY: On the Strains in the Interior of Beams, Phil. Trans. Roy. Soc.London, 153, 1863, 49-80. \\
 \noindent
[2] J.E. BJORK: Analytic D-Modules and Applications, Kluwer, 1993.\\
\noindent
[3] N. BOURBAKI: El\'{e}ments de Math\'{e}matiques, Alg\`{e}bre, Ch. 10, Alg\`{e}bre Homologique, Masson, Paris, 1980.   \\
\noindent 
[4] E. CARTAN: Les Syst\`{e}mes Diff\'{e}rentiels Ext\'{e}rieurs et Leurs Applications G\'{e}om\'{e}triques, Hermann, Paris, 1945.  \\
\noindent
[5] S.B. EDGAR: On Effective Constraints for the Riemann-Lanczos Systems of Equations. J. Math. Phys., 44, 2003, 5375-5385. \\ 
 http://arxiv.org/abs/gr-qc/0302014   \\
\noindent
[6] S.B. EDGAR, A. H\"{O}GLUND: The Lanczos potential for Weyl-Candidate Tensors Exists only in Four Dimension, General Relativity and Gravitation, 
 32, 12, 2000, 2307. \\
http://rspa.royalsocietypublishing.org/content/royprsa/453/1959/835.full.pdf   \\ 
 \noindent
[7] S.B. EDGAR, J.M.M. SENOVILLA: A Local Potential for the Weyl tensor in all dimensions, Classical and Quantum Gravity, 21, 2004, L133.\\
 http://arxiv.org/abs/gr-qc/0408071     \\
\noindent
[8] L.P. EISENHART: Riemannian Geometry, Princeton University Press, 1926.  \\
\noindent
[9] S.-T. HU: Introduction to Homological Algebra, Holden-Day, 1968.  \\
\noindent
[10] M. JANET: Sur les Syst\`{e}mes aux D\'{e}riv\'{e}es Partielles, Journal de Math., 8(3),1920, 65-151.  \\
\noindent
[11] M. KASHIWARA: Algebraic Study of Systems of Partial Differential Equations, M\'emoires de la Soci\'et\'e Math\'ematique de France 63, 1995, 
(Transl. from Japanese of his 1970 Master's Thesis).\\
\noindent
[12] E.R. KOLCHIN: Differential Algebra and Algebraic Groups, Academic Press, 1973.  \\
\noindent
[13] A. KUMPERA, D.C. SPENCER, Lie Equations, Ann. Math. Studies 73, Princeton University Press, Princeton, 1972.\\
\noindent
[14] C. LANCZOS: Lagrange Multiplier and Riemannian Spaces, Reviews of Modern Physics, 21, 1949, 497-502.  \\ 
\noindent
[15] C. LANCZOS: The Splitting of the Riemann Tensor, Rev. Mod. Phys. 34, 1962, 379-389.  \\
\noindent
[16] F. S. MACAULAY, The Algebraic Theory of Modular Systems, Cambridge Tracts 19, Cambridge University Press, London, 1916; 
Reprinted by Stechert-Hafner Service Agency, New York, 1964.\\
\noindent
[17] D.G. NORTHCOTT: An Introduction to Homological Algebra, Cambridge University Press, 1966.  \\
\noindent
[18] P. O'DONNELL, H. PYE: A Brief Historical Review of the Important Developments in Lanczos Potential Theory, EJTP, 24, 2010, 327-350.  \\
noindent
[19] J.-F. POMMARET: Systems of Partial Differential Equations and Lie Pseudogroups, Gordon and Breach, New York, 1978 
(Russian translation by MIR, Moscow, 1983) \\
\noindent
[20] J.-F. POMMARET: Differential Galois Theory, Gordon and Breach, New York, 1983. \\
\noindent
[21]  J.-F. POMMARET: Lie Pseudogroups and Mechanics, Gordon and Breach, New York, 1988.\\
\noindent
[22] J.-F. POMMARET: Partial Differential Equations and Group Theory,New Perspectives for Applications, Mathematics and its Applications 293, Kluwer, 1994.\\
http://dx.doi.org/10.1007/978-94-017-2539-2   \\
\noindent
[23] J.-F. POMMARET: Partial Differential Control Theory, Kluwer, 2001, 957 pp.\\
\noindent
[24] J.-F. POMMARET: Algebraic Analysis of Control Systems Defined by Partial Differential Equations, in Advanced Topics in Control Systems Theory, Lecture Notes in Control and Information Sciences 311, Chapter 5, Springer, 2005, 155-223.\\
\noindent
[25] J.-F. POMMARET: Gr\"{o}bner Bases in Algebraic Analysis: New perspectives for applications, Radon Series Comp. Appl. Math 2, 1-21, de Gruyter, 2007.\\
\noindent
[26] J.-F. POMMARET: Parametrization of Cosserat Equations, Acta Mechanica, 215, 2010, 43-55.\\
\noindent
[27] J.-F. POMMARET: Macaulay Inverse Systems Revisited, Journal of Symbolic Computation, 46, 2011, 1049-1069.  \\
http://dx.doi.org/10.1016/j.jsc.2011.05.007    \\
\noindent
[28] J.-F. POMMARET: Spencer Operator and Applications: From Continuum Mechanics to Mathematical Physics, in "Continuum Mechanics-Progress in Fundamentals and Engineering Applications", Dr. Yong Gan (Ed.), ISBN: 978-953-51-0447--6, InTech, 2012, Chapter 1, Available from: \\
http://www.intechopen.com/books/continuum-mechanics-progress-in-fundamentals-and-engineering-applications   \\
\noindent
[29] J.-F. POMMARET: Deformation Cohomology of Algebraic and Geometric Structures, arXiv:1207.1964.\\
http://arxiv.org/abs/1207.1964  \\
\noindent
[30] J.-F. POMMARET: The Mathematical Foundations of General Relativity Revisited, Journal of Modern Physics, 4, 2013, 223-239.\\
http://dx.doi.org/10.4236/jmp.2013.48A022   \\
\noindent
[31] J.-F. POMMARET: The Mathematical Foundations of Gauge Theory Revisited, Journal of Modern Physics, 2014, 5, 157-170.  \\
http://dx.doi.org/10.4236/jmp.2014.55026    \\
\noindent
[32] J.-F. POMMARET: From Thermodynamics to Gauge Theory: the Virial Theorem Revisited, in " Gauge Theories and Differential geometry ", NOVA Science Publishers, 2015, Chapter 1, 1-44.  \\
http://arxiv.org/abs/1504.04118  \\
\noindent
[33] J.-F. POMMARET: Airy,  Beltrami, Maxwell, Einstein and Lanczos potentials Revisited.  \\
http://arxiv.org/abs/1512.05982   \\
\noindent
[34] J.-F. POMMARET, A. QUADRAT: Algebraic Analysis of Linear Multidimensional  Control Systems, IMA Journal of Mathematical Control and Informations, 
16, 1999, 275-297.\\
\noindent
[35] A. QUADRAT, D. ROBERTZ: A Constructive Study of the Module Structure of Rings of Partial Differential Operators, Acta Applicandae Mathematicae, 
133, 2014, 187-234. \\
http://hal-supelec.archives-ouvertes.fr/hal-00925533  \\
\noindent
[36] M.D. ROBERTS: The Physical Interpretation of the Lanczos Tensor, Il Nuovo Cimento, B110, 1996, 1165-1176.  \\ 
http://arxiv.org/abs/gr-qc/9904006   \\   
\noindent
[37] J.J. ROTMAN: An Introduction to Homological Algebra, Pure and Applied Mathematics, Academic Press, 1979.\\
\noindent
[38] J.-P. SCHNEIDERS: An Introduction to D-Modules, Bull. Soc. Roy. Sci. Li\`{e}ge, 63, 1994, 223-295.  \\
\noindent
[39] D.C. SPENCER: Overdetermined Systems of Partial Differential Equations, Bull. Amer. Math. Soc., 75, 1965, 1-114.\\
\noindent
[40] E. VESSIOT:  Sur la Th\'{e}orie des Groupes Continus, Annales Scientifiques de l'Ecole Normale Sup\'{e}rieure, Vol. 20, 1903, 411-451.  \\
(Can be obtained from  \hspace{1cm}       http://www.numdam.org)\\

\vspace{8cm}

\noindent
{\bf 7) APPENDIX} by Alban QUADRAT (INRIA, Lille;  alban.quadrat@inria.fr):  \\  

\noindent
{\bf 7.1) Conformal Killing system}\\

Using an Euclidean or Minkowskian metric, the system of conformal Killing equations is:\\
\begin{equation}
\omega_{rj} \, \xi_i^r+\omega_{ir} \, \xi_j^r-\frac{2}{n} \, \omega_{ij} \, \xi_r^r=0.
\end{equation}
but any other choice could be convenient too as we have explained in the paper. \\

\noindent
{\bf 7.2) Weyl tensor in dimension 2}\\

Let us consider $n=2$ and the euclidean metric $\omega_{ij}=\delta_j^i$ for $1 \leq i < j \leq 2$. If $D=\Q[d_1, d_2]$ denotes the commutative ring of PD operators in $d_1=\frac{\partial}{\partial x_1}$ and $d_2=\frac{\partial}{\partial x_2}$, then the system can be rewritten as $W_2 \, \eta=0$, where $\eta=(\xi^1 \quad \xi^2)^T$ and:
$$W_2=\left( \begin {array}{cc} d_{{1}}&-d_{{2}}\\ \noalign{\medskip} \, d_{{2}
}&d_{{1}}\\ \noalign{\medskip}-d_{{1}}&d_{{2}}\end {array} \right) $$

If we note $\lambda_2=\left( \begin {array}{ccc} 1&0&1\end {array} \right)$, then we have $\lambda_2 \, W_2=0$, which shows that the last row of $W_2$ is minus the first one. Hence, the $D$-module ${\rm im}_D(W_2)=\{\lambda \, W_2 \; | \; \lambda \in D^{1 \times 2}\}$ generated by the rows of $W_2$ can be generated by the first two rows of $W_2$, i.e., we have ${\rm im}_D(W_2)={\rm im}_D(R_2)$, where $R_2$ is the following matrix:
$$R_2=\left( \begin {array}{cc} d_{{1}}&-d_{{2}}\\ \noalign{\medskip} \, d_{{2}
}&d_{{1}}\end {array} \right) $$

We can check that the $D$-module $\ker_D(R_2)=\{ \mu \in D^{1 \times 2} \; | \; \mu \, R_2=0\}$, called the {\em second syzygy module} ([4]) of the $D$-module $M=D^2/(D^2 \,R_2)$, is reduced to 0, i.e., $R_2$ has full row rank. Hence, we obtain the following {\em finite free resolution} of $M$:
\[0 \rightarrow D^2 \stackrel{{R_2}}{\longrightarrow}  D^2 \rightarrow  M \rightarrow 0 \]

\noindent
{\bf 7.3) Weyl tensor in dimension 3}\\

Let us consider $n=3$ and the euclidean metric $\omega_{ij}=\delta_j^i$ for $1 \leq i < j \leq 3$. If $D=\Q[d_1, d_2, d_3]$ denotes the commutative ring of PD operators in $d_i=\frac{\partial}{\partial x_i}$ for $i=1, \, 2,  \, 3$, then the system can be rewritten as $W_3 \, \eta=0$, where $\eta=(\xi^1 \quad \xi^2 \quad \xi^3)^T$ and:
$$W_3 =\left( \begin {array}{ccc} \frac{4}{3} \,d_{{1}}&-\frac{2}{3}\,d_{{2}}& -\frac{2}{3} \,d_{{3}}
\\ \noalign{\medskip} \, d_{{2}}&d_{{1}}&0\\ \noalign{\medskip} \, d_{{3}}&0&d
_{{1}}\\ \noalign{\medskip}-\frac{2}{3}\,d_{{1}}& \frac{4}{3} \,d_{{2}}&-\frac{2}{3}\,d_{{3}}
\\ \noalign{\medskip}0&d_{{3}}&d_{{2}}\\ \noalign{\medskip}-\frac{2}{3}\,d_{{1
}}& -\frac{2}{3}\,d_{{2}}& \frac{4}{3}\,d_{{3}}\end {array} \right)$$ 

If we note
$\lambda_3=\left( \begin {array}{cccccc} 1&0&0&1&0&1\end {array} \right) \in D^{1 \times 6}$, 
then we have $\lambda_3 \, W_3=0$ and shows that the last row of $W_3$ is a linear combination of the first and fourth rows of $W_3$. Hence, the $D$-module ${\rm im}_D(W_3)$, generated by the rows of $W_3$, can be generated by the first five rows of $W_3$, i.e., we have ${\rm im}_D(W_3)={\rm im}_D(R_3)$, where $R_3$ is the following matrix:
$$R_3=\left( \begin {array}{ccc} \frac{4}{3}\,d_{{1}}&-\frac{2}{3}\,d_{{2}}&-\frac{2}{3}\,d_{{3}}
\\ \noalign{\medskip} \, d_{{2}}&d_{{1}}&0\\ \noalign{\medskip} \, d_{{3}}&0&d
_{{1}}\\ \noalign{\medskip}-\frac{2}{3}\,d_{{1}}&\frac{4}{3}\,d_{{2}}&-\frac{2}{3}\,d_{{3}}
\\ \noalign{\medskip}0&d_{{3}}&d_{{2}}\end {array} \right) $$

Using the {\sc OreModules} package, we can show that the left kernel of $R_3$, i.e., the $D$-module $\ker_D(R_3)=\{\lambda \in D^{1 \times 3} \; | \; \lambda \, R_3=0\}$, is generated by the rows of the following matrix $S_3$: \\
$$
\begin{array}{c}
\left( \begin {array}{ccccc} -d_{{2}} \, d_{{3}} \, d_{{1}}& d_{{2}}^{2} \, d_{{
3}}+d_{{3}}^{3}&-d_{{2}}^{3}-d_{{2}} \, d_{{3}}^{2}&-2\,d_{{2}} \, d_{{3
}} \, d_{{1}}&d_{{1}} \, d_{{2}}^{2}-d_{{1}} \, d_{{3}}^{2}
\\ \noalign{\medskip}-d_{{1}}^{2} \, d_{{3}}+2\,d_{{2}}^{2} \, d_{{3}}+d_
{{3}}^{3}&0&-2\,d_{{1}} \, d_{{2}}^{2}-2\,d_{{1}} \, d_{{3}}^{2}&-2\, d_{
{1}}^{2} \, d_{{3}}+d_{{2}}^{2} \, d_{{3}}-d_{{3}}^{3}&2\,d_{{1}}^{2} \, d_
{{2}}+2\,d_{{2}} \, d_{{3}}^{2}\\ \noalign{\medskip} d_{{1}}^{2} \, d_{{2}}
+2\, d_{{2}}^{3}+d_{{2}} \, d_{{3}}^{2}&-2\,d_{{1}} \, d_{{2}}^{2}-2\,d_{
{1}} \, d_{{3}}^{2}&0&2\, d_{{1}}^{2} \, d_{{2}}+d_{{2}}^{3}-d_{{2}} \, d_{{
3}}^{2}&2\, d_{{1}}^{2} \, d_{{3}}+2\, d_{{2}}^{2} \, d_{{3}}
\\ \noalign{\medskip}2\,d_{{2}} \, d_{{3}} \, d_{{1}}&-d_{{1}}^{2} \, d_{{3}}-d
_{{3}}^{3}&-d_{{1}}^{2} \, d_{{2}}+d_{{2}} \, d_{{3}}^{2}&d_{{2}} \, d_{{3}} \, d
_{{1}}& d_{{1}}^{3}+d_{{1}} \, d_{{3}}^{2}\\ \noalign{\medskip} d_{{1}}
^{3}+2\,d_{{1}} \, d_{{2}}^{2}-d_{{1}} \, d_{{3}}^{2}&-2\, d_{{1}}^{2} \, d_
{{2}}-2\,d_{{2}} \, d_{{3}}^{2}&2\, d_{{1}}^{2} \, d_{{3}}+2\, d_{{2}}^{2}
d_{{3}}&2\, d_{{1}}^{3}+d_{{1}} \, d_{{2}}^{2}+d_{{1}} \, d_{{3}}^{2}&0
\end {array} \right)
\end{array}$$
i.e., we have $\ker_D(R_3)={\rm im}_D(S_3)=\{\mu \in D^{1 \times 5} \; | \; \mu \, S_3 \}$. The computation of $S_3$ takes 0.026 CPU seconds with Maple~18 on Mac OS~10.10.5 equipped with 2.8~GHz Intel Core~i7 and 16~Go. 
For more details on the left kernel/syzygy computation, see Algorithm~1 on page 330 of ([1]) (see also [3]). Similarly, the left kernel of $T_3$ is generated by the rows of the following matrix
$$T_3=\left( \begin {array}{ccccc} -2\,d_{{3}}&0&d_{{1}}&-2\,d_{{3}}&-d_{{2
}}\\ \noalign{\medskip}2\,d_{{1}}&-d_{{2}}&d_{{3}}&0&0
\\ \noalign{\medskip}0&d_{{1}}&0&-2\,d_{{2}}&d_{{3}}\end {array}
 \right) $$
i.e., we have $\ker_D(S_3)={\rm im}_D(T_3)$. We can check that that the matrix $T_3$ has full row rank, i.e., $\ker_D(T_3)=0$, or equivalently that the rows of $T_3$ are $D$-linearly independent. Finally, if we denote by $M=D^3/(D^5 \, R_3)=D^3/(D^6 \, W_3)$ the $D$-module {\em finitely presented} by $R_3$, associated with the system, then we obtain the following finite free resolution:\\
\[  0 \rightarrow D^3 \stackrel{T_3}{\longrightarrow}  D^5 \stackrel{S_3}{\longrightarrow}D^5 \stackrel{R_3}{\longrightarrow} D^3 \rightarrow M \rightarrow 0\]

\noindent
{\bf 7.4) Weyl tensor in dimension 4}\\

Let us consider $n=4$ and the Minkowski metric $\omega=(1, \; 1, \; 1, \; -1)$. If $D=\Q[d_1, d_2, d_3, d_4]$ denotes the commutative ring of PD operators in $d_i=\frac{\partial}{\partial x_i}$ for $i=1, \ldots, 4$, then the system can be rewritten as $W_4 \, \eta=0$, where $\eta=(\xi^1 \quad \xi^2 \quad \xi^3 \quad \xi^4)^T$ and:
$$W_4= \left( \begin {array}{cccc} \frac{3}{2}\,d_{{1}}&-\frac{1}{2}\,d_{{2}}&-\frac{1}{2}\,d_{{3}}&
-\frac{1}{2}\,d_{{4}}\\ \noalign{\medskip}d_{{2}}&d_{{1}}&0&0
\\ \noalign{\medskip}d_{{3}}&0&d_{{1}}&0\\ \noalign{\medskip}d_{{4}}&0
&0&-d_{{1}}\\ \noalign{\medskip}-\frac{1}{2}\,d_{{1}}&\frac{3}{2}\,d_{{2}}&-\frac{1}{2}\,d_{{3}
}&-\frac{1}{2}\,d_{{4}}\\ \noalign{\medskip}0&d_{{3}}&d_{{2}}&0
\\ \noalign{\medskip}0&d_{{4}}&0&-d_{{2}}\\ \noalign{\medskip}-\frac{1}{2}\,d_{
{1}}&-\frac{1}{2}\,d_{{2}}&\frac{3}{2}\,d_{{3}}&-\frac{1}{2}\,d_{{4}}\\ \noalign{\medskip}0&0&
d_{{4}}&-d_{{3}}\\ \noalign{\medskip}\frac{1}{2}\,d_{{1}}&\frac{1}{2}\,d_{{2}}&\frac{1}{2}\,
d_{{3}}&-\frac{3}{2}\,d_{{4}}\end {array} \right) $$

If we note $\lambda_4=\left( \begin {array}{cccccccccc} 1&0&0&0&1&0&0&1&0&-1\end {array}
 \right)$,  then we have $\lambda_4 \, W_4=0$, which shows that the last row of $W_4$ is a linear combination of the first, fourth, and eighth rows. Hence, the $D$-module ${\rm im}_D(W_4)$, generated by the rows of $W_4$, can be generated by the first nine rows of $W_4$, i.e., we have ${\rm im}_D(W_4)={\rm im}_D(R_4)$, where $R_4$ is the following matrix:
$$R_4=\left( \begin {array}{cccc} \frac{3}{2}\,d_{{1}}&-\frac{1}{2}\,d_{{2}}&-\frac{1}{2}\,d_{{3}}&
-\frac{1}{2}\,d_{{4}}\\ \noalign{\medskip}d_{{2}}&d_{{1}}&0&0
\\ \noalign{\medskip}d_{{3}}&0&d_{{1}}&0\\ \noalign{\medskip}d_{{4}}&0
&0&-d_{{1}}\\ \noalign{\medskip}-\frac{1}{2}\,d_{{1}}&\frac{3}{2}\,d_{{2}}&-\frac{1}{2}\,d_{{3}
}&-\frac{1}{2}\,d_{{4}}\\ \noalign{\medskip}0&d_{{3}}&d_{{2}}&0
\\ \noalign{\medskip}0&d_{{4}}&0&-d_{{2}}\\ \noalign{\medskip}-\frac{1}{2}\,d_{
{1}}&-\frac{1}{2}\,d_{{2}}&\frac{3}{2}\,d_{{3}}&-\frac{1}{2}\,d_{{4}}\\ \noalign{\medskip}0&0&
d_{{4}}&-d_{{3}}\end {array} \right) \in D^{9 \times 4}.$$

We can show that the left kernel of $R_4$, i.e., the $D$-module $\ker_D(R_4)$, is generated by the rows of the following matrix $S_4\in D^{10\times 9}$: \\
{\footnotesize $$\hspace{-2.5cm} 
\begin{array}{c} \left( \begin {array}{ccccccccc} -d_{{1}} \, d_{{3}}&d_{{2}} \,d_{{3}}&-d_{
{2}}^{2}-d_{{4}}^{2}&d_{{3}} \, d_{{4}}&-2\,d_{{1}} \, d_{{3}}&d_{{1}} \, d_{{
2}}&0&-d_{{1}} \, d_{{3}}&d_{{4}} \, d_{{1}}\\ \noalign{\medskip}d_{{1}} \, d_{{2
}}&d_{{3}}^{2}+d_{{4}}^{2}&-d_{{2}} \, d_{{3}}&-d_{{4}} \, d_{{2}}&d_{{1}} \, d
_{{2}}&-d_{{1}}\, d_{{3}}&-d_{{4}} \, d_{{1}}&2\,d_{{1}} \, d_{{2}}&0
\\ \noalign{\medskip}d_{{3}} \, d_{{4}}&0&-d_{{4}} \, d_{{1}}&-d_{{1}} \, d_{{3}}&
-d_{{3}} \, d_{{4}}&d_{{4}} \, d_{{2}}&d_{{2}} \, d_{{3}}&0&d_{{1}}^{2}-d_{{2}}
^{2}\\ \noalign{\medskip}d_{{4}} \, d_{{2}}&-d_{{4}} \, d_{{1}}&0&-d_{{1}} \, d_{
{2}}&0&d_{{3}} \, d_{{4}}& d_{{1}}^{2}-d_{{3}}^{2}&-d_{{4}} \, d_{{2}}&d_{{
2}} \, d_{{3}}\\ \noalign{\medskip}2\,d_{{2}} \, d_{{3}}&-d_{{1}} \, d_{{3}}&-d_{{
1}} \, d_{{2}}&0&d_{{2}} \, d_{{3}}& d_{{1}}^{2}+d_{{4}}^{2}&-d_{{3}} \, d_{{4}}
&d_{{2}} \, d_{{3}}&-d_{{4}} \, d_{{2}}\\ \noalign{\medskip}2\, d_{{2}}^{2}-2
\,d_{{3}}^{2}&-2\,d_{{1}} \, d_{{2}}&2\,d_{{1}} \, d_{{3}}&0& d_{{1}}^{2}+
d_{{2}}^{2}-d_{{3}}^{2}+d_{{4}}^{2}&0&-2\,d_{{4}} \, d_{{2}}&-d_{{1}}
^{2}+d_{{2}}^{2}-d_{{3}}^{2}-d_{{4}}^{2}&2\,d_{{3}} \, d_{{4}}
\\ \noalign{\medskip}d_{{1}}^{2}-3\,d_{{1}}^{2}+3\, d_{{3}}^{2}+
d_{{4}}^{2}&2\,d_{{1}} \, d_{{2}}&-4\,d_{{1}} \, d_{{3}}&-2\,d_{{4}} \, d_{{1}}&-2
\,d_{{2}}^{2}-2\, d_{{4}}^{2}&2\,d_{{2}} \, d_{{3}}&4\,d_{{4}} \, d_{{2}}&
3\, d_{{1}}^{2}-3\, d_{{2}}^{2}+d_{{3}}^{2}+d_{{4}}^{2}&-2\,d_{{3
}} \, d_{{4}}\\ \noalign{\medskip}0&d_{{3}} \, d_{{4}}&0&-d_{{2}} \, d_{{3}}&0&-d_
{{4}} \, d_{{1}}&0&0&d_{{1}} \, d_{{2}}\\ \noalign{\medskip}0&d_{{4}} \, d_{{2}}&-
d_{{3}} \, d_{{4}}&-d_{{2}}^{2}+d_{{3}}^{2}&-d_{{4}} \, d_{{1}}&0&d_{{1}} \, d
_{{2}}&d_{{4}} \, d_{{1}}&-d_{{1}} \, d_{{3}}\\ \noalign{\medskip}0&0&d_{{4}} \, d
_{{2}}&-d_{{2}} \, d_{{3}}&0&-d_{{4}} \, d_{{1}}&d_{{1}} \, d_{{3}}&0&0
\end {array} \right)
\end{array}
$$}

i.e., we have $\ker_D(R_4)={\rm im}_D(S_4)$. In fact, using the {\sc OreModules} package, we can prove that $\ker_D(R_4)={\rm im}_D(\overline{S}_4)$, where $\overline{S}_4=(S_4^T \quad {S_4'}^T)^T$, where the matrix $S_4' \in D^{6 \times 9}$ is defined by $S'_4\in D^{6\times 9}$:  \\
 {\tiny $$
\begin{array}{c}
\hspace{-2cm} 
\left( \begin {array}{ccccccccc} -d_{{2}}\, d_{{3}} \, d_{{4}}&0&0&0&-d_{{2}
} \, d_{{3}} \, d_{{4}}&d_{{3}}^{2} \, d_{{4}}-d_{{4}}^{3}&-d_{{3}}^{3}+d_{{
3}} \, d_{{4}}^{2}&-2\,d_{{2}} \, d_{{3}} \, d_{{4}}&d_{{2}} \, d_{{3}} ^{2}+d_{{2}
} \, d_{{4}}^{2}\\ \noalign{\medskip}-d_{{1}} \, d_{{3}} \, d_{{4}}&0& d_{{3}}^
{2} \, d_{{4}}-d_{{4}}^{3}&-d_{{3}}^{3}+d_{{3}} \, d_{{4}}^{2}&-d_{{1}} \, d
_{{3}} \, d_{{4}}&0&0&-2\,d_{{1}} \, d_{{3}} \, d_{{4}}&d_{{1}} \, d_{{3}}^{2}+d_{{1
}} \, d_{{4}}^{2}\\ \noalign{\medskip}-d_{{2}}^{2} \, d_{{4}}+d_{{3}}^{2
} \, d_{{4}}&0&0&0&-d_{{2}}^{2} \, d_{{4}}+2\, d_{{3}}^{2} \, d_{{4}}-d_{{4}}
^{3}&0&-2\,d_{{2}} \, d_{{3}}^{2}+2\,d_{{2}} \, d_{{4}}^{2}&-2\, d_{{2}}^
{2} \, d_{{4}}+d_{{3}}^{2} \, d_{{4}}+d_{{4}}^{3}&2\, d_{{2}}^{2} \, d_{{3}}-
2\,d_{{3}} \, d_{{4}}^{2}\\ \noalign{\medskip} d_{{2}}^{2} \, d_{{3}}+d_{{
3}}^{3}&0&0&0& d_{{2}}^{2} \, d_{{3}}+2\, d_{{3}}^{3}-d_{{3}} \, d_{{4}}^
{2}&-2\,d_{{2}} \, d_{{3}}^{2}+2\,d_{{2}} \, d_{{4}}^{2}&0&2\, d_{{2}}^{2
} \, d_{{3}}+d_{{3}}^{3}+d_{{3}} \, d_{{4}}^{2}&-2\, d_{{2}}^{2} \, d_{{4}}-2
\, d_{{3}}^{2} \, d_{{4}}\\ \noalign{\medskip}d_{{2}} \, d_{{3}} \, d_{{4}}&0&0&0
&2\,d_{{2}} \, d_{{3}} \, d_{{4}}&-d_{{2}}^{2} \, d_{{4}}+d_{{4}}^{3}&-d_{{2}
}^{2} \, d_{{3}}-d_{{3}} \, d_{{4}}^{2}&d_{{2}} \, d_{{3}} \, d_{{4}}& d_{{2}}^{3}
-d_{{2}} \, d_{{4}}^{2}\\ \noalign{\medskip} d_{{2}}^{3}+d_{{2}} \, d_{{3}
}^{2}&0&0&0& d_{{2}}^{3}+2\,d_{{2}} \, d_{{3}}^{2}+d_{{2}} \, d_{{4}}^{2
}&-2\,d_{{2}}^{2} \, d_{{3}}+2\,d_{{3}} \, d_{{4}}^{2}&-2\, d_{{2}}^{2} \, d_{
{4}}-2\, d_{{3}}^{2} \, d_{{4}}&2\, d_{{2}}^{3}+d_{{2}} \, d_{{3}}^{2}-d_{
{2}} \, d_{{4}}^{2}&0\end {array} \right).
\end{array}$$} 
But, we have $S_4'=F_4 \, S_4$, where
$$F_4=\left( \begin {array}{cccccccccc} 0&0&0&d_{{3}}&-d_{{4}}&0&0&0&0&-d_{
{1}}\\ \noalign{\medskip}d_{{4}}&0&0&0&0&0&0&0&-d_{{3}}&d_{{2}}
\\ \noalign{\medskip}0&0&-d_{{3}}&d_{{2}}&0&-d_{{4}}&0&0&-d_{{1}}&0
\\ \noalign{\medskip}-d_{{1}}&0&d_{{4}}&0&d_{{2}}&-2\,d_{{3}}&-d_{{3}
}&0&0&0\\ \noalign{\medskip}0&0&-d_{{2}}&0&d_{{4}}&0&0&d_{{1}}&0&0
\\ \noalign{\medskip}0&d_{{1}}&0&d_{{4}}&d_{{3}}&-d_{{2}}&-d_{{2}}&0&0
&0\end {array} \right) \in D^{6 \times 10},$$
which shows that ${\rm im}_D(\overline{S}_4)={\rm im}_D(S_4)$. Similarly, we get that $\ker_D(S_4)={\rm im}_D(T_4)$, where the matrix $T_4\in D^{9\times 10}$ is defined by
{\footnotesize $$\hspace{-1cm} \left( \begin {array}{cccccccccc} d_{{4}} \, d_{{2}}&d_{{3}}  \, d_{{4}}&0&0&0
&0&0&-d_{{3}}^{2}-d_{{4}}^{2}&-d_{{2}}  \, d_{{3}}& d_{{2}}^{2}+d_{{4
}}^{2}\\ \noalign{\medskip}d_{{4}}  \, d_{{1}}&0&-d_{{3}}^{2}-d_{{4}}^
{2}&-d_{{2}}  \, d_{{3}}&d_{{4}}  \, d_{{2}}&d_{{3}}  \, d_{{4}}&d_{{3}}  \, d_{{4}}&0&-d_
{{1}}  \, d_{{3}}&2\,d_{{1}}  \, d_{{2}}\\ \noalign{\medskip}d_{{1}}  \, d_{{2}}&d_{{
1}}  \, d_{{3}}&d_{{4}}  \, d_{{2}}&-d_{{3}}  \, d_{{4}}&-d_{{2}}^{2}+d_{{3}}^{2}
&d_{{2}}  \, d_{{3}}&0&-2\,d_{{4}}  \, d_{{1}}&0&2\,d_{{4}}  \, d_{{1}}
\\ \noalign{\medskip}-d_{{3}}  \, d_{{4}}&-d_{{4}}  \, d_{{2}}&d_{{1}}  \, d_{{3}}&-d
_{{1}}  \, d_{{2}}&0&d_{{4}}  \, d_{{1}}&0&2\,d_{{2}}  \, d_{{3}}& d_{{1}}^{2}+d_{{
4}}^{2}&-2\,d_{{2}}  \, d_{{3}}\\ \noalign{\medskip}d_{{2}}  \, d_{{3}}& d_{{1}
}^{2}-d_{{2}}^{2}&0&d_{{4}}  \, d_{{1}}&d_{{1}}  \, d_{{3}}&-d_{{1}}  \, d_{{2}}&
-d_{{1}}  \, d_{{2}}&0&d_{{4}}  \, d_{{2}}&2\,d_{{3}}  \, d_{{4}}
\\ \noalign{\medskip} d_{{1}}^{2}-d_{{3}}^{2}&d_{{2}}  \, d_{{3}}&-d_{{4}
}  \, d_{{1}}&0&-d_{{1}}  \, d_{{2}}&2\,d_{{1}}  \, d_{{3}}&d_{{1}}  \, d_{{3}}&-2\,d_{{4}}  \, 
d_{{2}}&d_{{3}}  \, d_{{4}}&0\\ \noalign{\medskip}0&-d_{{4}}  \, d_{{1}}&-d_{{2
}}  \, d_{{3}}&-d_{{2}}^{2}-d_{{4}}^{2}&d_{{3}}  \, d_{{4}}&2\,d_{{4}}  \, d_{{2}
}&d_{{4}}  \, d_{{2}}&2\,d_{{1}}  \, d_{{3}}&d_{{1}}  \, d_{{2}}&0
\\ \noalign{\medskip}0&-d_{{3}}  \, d_{{4}}&0&-d_{{1}}  \, d_{{3}}&d_{{4}}  \, d_{{1}
}&0&0& d_{{3}}^{2}+d_{{4}}^{2}&0& d_{{1}}^{2}-d_{{3}}^{2}
\\ \noalign{\medskip}0&-d_{{3}}  \, d_{{4}}&-d_{{1}}  \, d_{{2}}&0&d_{{4}}  \, d_{{1}
}&0&0& d_{{1}}^{2}-d_{{2}}^{2}+d_{{3}}^{2}+d_{{4}}^{2}&d_{{2}}  \, d
_{{3}}&0\end {array} \right),$$}
and $\ker_D(T_4)={\rm im}_D(U_4)$, where the matrix $U_4 \in D^{4\times 9}$ is defined by:
$$U_4=\left( \begin {array}{ccccccccc} -2\,d_{{3}}&0&0&-d_{{2}}&d_{{4}}&0&d
_{{1}}&0&-2\,d_{{3}}\\ \noalign{\medskip}2\,d_{{1}}&-d_{{2}}&-d_{{4}}&0
&0&0&d_{{3}}&0&0\\ \noalign{\medskip}0&d_{{1}}&0&d_{{3}}&0&-d_{{4}}&0&
-2\,d_{{2}}&0\\ \noalign{\medskip}0&0&d_{{1}}&0&-d_{{3}}&-d_{{2}}&0&-2
\,d_{{4}}&2\,d_{{4}}\end {array} \right).$$
Finally, we can check that $U_4$ has full row rank, i.e., $\ker_D(U_4)=0$, which shows that the $D$-module $M=D^4/(D^9 \, R_4)=D^4/(D^{10} \, W_4)$, associated with the system, admits the following finite free resolution:\\
\[ 0 \rightarrow D^4 \stackrel{U_4}{\longrightarrow} D^9 \stackrel{T_4}{\longrightarrow} D^{10} \stackrel{S_4}{\longrightarrow} D^9 \stackrel{R_4}{\longrightarrow} D^4 \rightarrow  M \rightarrow 0  \]
Finally, we note that a free resolution of $M$ can be computed in 0.616 CPU seconds when $n=4$. \\

\noindent
{\bf 7.5) Weyl tensor in dimension 5}\\

Let us consider $n=5$ and the euclidean metric $\omega_{ij}=\delta_j^i$ for $1 \leq i < j \leq 5$. If $D=\Q[d_1, d_2, d_3, d_4,d_5]$ denotes the commutative ring of PD operators in $d_i=\frac{\partial}{\partial x_i}$ for $i=1, \ldots, 5$, then the system can be rewritten as $W_5 \, \eta=0$, where $\eta=(\xi^1 \quad \ldots  \quad \xi^5)^T$ and:
$$W_5=\left( \begin {array}{ccccc} \frac{8}{5}\,d_{{1}}&-\frac{2}{5}\,d_{{2}}&-\frac{2}{5}\,d_{{3}}
&-\frac{2}{5}\,d_{{4}}&-\frac{2}{5}\,d_{{5}}\\ \noalign{\medskip}d_{{2}}&d_{{1}}&0&0&0
\\ \noalign{\medskip}d_{{3}}&0&d_{{1}}&0&0\\ \noalign{\medskip}d_{{4}}
&0&0&d_{{1}}&0\\ \noalign{\medskip}d_{{5}}&0&0&0&d_{{1}}
\\ \noalign{\medskip}-\frac{2}{5}\,d_{{1}}&\frac{8}{5}\,d_{{2}}&-\frac{2}{5}\,d_{{3}}&-\frac{2}{5}\,d_
{{4}}&-\frac{2}{5}\,d_{{5}}\\ \noalign{\medskip}0&d_{{3}}&d_{{2}}&0&0
\\ \noalign{\medskip}0&d_{{4}}&0&d_{{2}}&0\\ \noalign{\medskip}0&d_{{5
}}&0&0&d_{{2}}\\ \noalign{\medskip}-\frac{2}{5}\,d_{{1}}&-\frac{2}{5}\,d_{{2}}&\frac{8}{5}\,d_
{{3}}&-\frac{2}{5}\,d_{{4}}&-\frac{2}{5}\,d_{{5}}\\ \noalign{\medskip}0&0&d_{{4}}&d_{{
3}}&0\\ \noalign{\medskip}0&0&d_{{5}}&0&d_{{3}}\\ \noalign{\medskip}-2
/5\,d_{{1}}&-\frac{2}{5}\,d_{{2}}&-\frac{2}{5}\,d_{{3}}&\frac{8}{5}\,d_{{4}}&-\frac{2}{5}\,d_{{5}}
\\ \noalign{\medskip}0&0&0&d_{{5}}&d_{{4}}\\ \noalign{\medskip}-\frac{2}{5}\,d
_{{1}}&-\frac{2}{5}\,d_{{2}}&-\frac{2}{5}\,d_{{3}}&-\frac{2}{5}\,d_{{4}}&\frac{8}{5}\,d_{{5}}
\end {array} \right) \in D^{15 \times 5}.$$

If we note $\lambda_5=\left( \begin {array}{ccccccccccccccc} 1&0&0&0&0&1&0&0&0&1&0&0&1&0&1
\end {array} \right)$, then we have $\lambda_5 \, W_5=0$, which shows that the $D$-module ${\rm im}_D(W_5)$, generated by the rows of $W_5$, can be generated by the first fourteen rows of $W_5$, i.e., we have ${\rm im}_D(W_5)={\rm im}_D(R_5)$, where $R_5$ is the following matrix:
$$R_5= \left( \begin {array}{ccccc} \frac{8}{5}\,d_{{1}}&-\frac{2}{5}\,d_{{2}}&-\frac{2}{5}\,d_{{3}}
&-\frac{2}{5}\,d_{{4}}&-\frac{2}{5}\,d_{{5}}\\ \noalign{\medskip}d_{{2}}&d_{{1}}&0&0&0
\\ \noalign{\medskip}d_{{3}}&0&d_{{1}}&0&0\\ \noalign{\medskip}d_{{4}}
&0&0&d_{{1}}&0\\ \noalign{\medskip}d_{{5}}&0&0&0&d_{{1}}
\\ \noalign{\medskip}-\frac{2}{5}\,d_{{1}}&\frac{8}{5}\,d_{{2}}&-\frac{2}{5}\,d_{{3}}&-\frac{2}{5}\,d_
{{4}}&-\frac{2}{5}\,d_{{5}}\\ \noalign{\medskip}0&d_{{3}}&d_{{2}}&0&0
\\ \noalign{\medskip}0&d_{{4}}&0&d_{{2}}&0\\ \noalign{\medskip}0&d_{{5
}}&0&0&d_{{2}}\\ \noalign{\medskip}-\frac{2}{5}\,d_{{1}}&-\frac{2}{5}\,d_{{2}}&\frac{8}{5}\,d_
{{3}}&-\frac{2}{5}\,d_{{4}}&-\frac{2}{5}\,d_{{5}}\\ \noalign{\medskip}0&0&d_{{4}}&d_{{
3}}&0\\ \noalign{\medskip}0&0&d_{{5}}&0&d_{{3}}\\ \noalign{\medskip}-2
/5\,d_{{1}}&-\frac{2}{5}\,d_{{2}}&-\frac{2}{5}\,d_{{3}}&\frac{8}{5}\,d_{{4}}&-\frac{2}{5}\,d_{{5}}
\\ \noalign{\medskip}0&0&0&d_{{5}}&d_{{4}}\end {array} \right) \in D^{14 \times 5}.$$

Using the {\sc OreModules} package, we can show that the left kernel of $R_5$, i.e., the $D$-module $\ker_D(R_5)$, is generated by the rows of a matrix 
$S_5=(S_A \quad S_B) \in D^{35 \times 14}$, where the matrices $S_A \in D^{35 \times 7}$ and $S_B \in D^{35 \times 7}$ are respectively defined by:
{\tiny $$\hspace{-1cm}
S_A=\left( \begin {array}{ccccccc} d_{{4}}^{2}-d_{{4}} \, d_{{2}}&0&0&0&0&-d_{{4}} \, d_{{2}
}&d_{{3}} \, d_{{4}}\\ \noalign{\medskip}d_{{2}} \, d_{{3}}&0&0&0&0&d_{{2}} \, d_{
{3}}&d_{{4}}^{2}-d_{{5}}^{2}\\ \noalign{\medskip}-d_{{4}} \, d_{{1}}&0
&d_{{3}} \, d_{{4}}&-d_{{3}}^{2}+d_{{5}}^{2}&-d_{{5}} \, d_{{4}}&-d_{{4}}d
_{{1}}&0\\ \noalign{\medskip}d_{{1}} \, d_{{3}}&0&d_{{4}}^{2}-{d_{{5}}}^
{2}&-d_{{3}} \, d_{{4}}&d_{{3}} \, d_{{5}}&d_{{1}} \, d_{{3}}&0
\\ \noalign{\medskip}d_{{3}} \, d_{{4}}&0&0&0&0&2\,d_{{3}} \, d_{{4}}&-d_{{4}}
d_{{2}}\\ \noalign{\medskip}d_{{3}}^{2}-d_{{4}}^{2}&0&0&0&0&2\,{d_
{{3}}}^{2}-2\,d_{{4}}^{2}&-2\,d_{{2}} \, d_{{3}}\\ \noalign{\medskip}{d_
{{2}}}^{2}-2\,d_{{3}}^{2}+d_{{4}}^{2}&0&0&0&0&d_{{2}}^{2}-3\,{d_
{{3}}}^{2}+3\,d_{{4}}^{2}-d_{{5}}^{2}&2\,d_{{2}} \, d_{{3}}
\\ \noalign{\medskip}-d_{{4}} \, d_{{1}}&d_{{4}} \, d_{{2}}&0&-d_{{2}}^{2}+{
d_{{5}}}^{2}&-d_{{5}} \, d_{{4}}&-2\,d_{{4}} \, d_{{1}}&0\\ \noalign{\medskip}
d_{{1}} \, d_{{2}}&-d_{{3}}^{2}+2\,d_{{4}}^{2}-d_{{5}}^{2}&d_{{2}} \, d_
{{3}}&-2\,d_{{4}} \, d_{{2}}&d_{{2}} \, d_{{5}}&d_{{1}} \, d_{{2}}&d_{{1}} \, d_{{3}}
\\ \noalign{\medskip}d_{{5}} \, d_{{4}}&0&0&-d_{{1}} \, d_{{5}}&-d_{{4}} \, d_{{1}
}&0&0\\ \noalign{\medskip}d_{{3}} \, d_{{5}}&0&-d_{{1}} \, d_{{5}}&0&-d_{{1}}d
_{{3}}&0&0\\ \noalign{\medskip}2\,d_{{3}} \, d_{{4}}&0&-d_{{4}} \, d_{{1}}&-d_
{{1}} \, d_{{3}}&0&d_{{3}} \, d_{{4}}&0\\ \noalign{\medskip}2\,d_{{3}}^{2}-2
\,d_{{4}}^{2}&0&-2\,d_{{1}} \, d_{{3}}&2\,d_{{4}} \, d_{{1}}&0&d_{{3}}^{2}
-d_{{4}}^{2}&0\\ \noalign{\medskip}d_{{2}} \, d_{{5}}&-d_{{1}} \, d_{{5}}&0&0
&-d_{{1}} \, d_{{2}}&0&0\\ \noalign{\medskip}2\,d_{{4}} \, d_{{2}}&-d_{{4}} \, d_{
{1}}&0&-d_{{1}} \, d_{{2}}&0&d_{{4}} \, d_{{2}}&0\\ \noalign{\medskip}3\,{d_{{
3}}}^{2}-3\,d_{{4}}^{2}&-2\,d_{{1}} \, d_{{2}}&0&2\,d_{{4}} \, d_{{1}}&0&{d_
{{1}}}^{2}-d_{{2}}^{2}+4\,d_{{3}}^{2}-5\,d_{{4}}^{2}+{d_{{5}}}^{
2}&-2\,d_{{2}} \, d_{{3}}\\ \noalign{\medskip}d_{{1}}^{2}-6\,{d_{{3}}}^{
2}+6\,d_{{4}}^{2}-d_{{5}}^{2}&2\,d_{{1}} \, d_{{2}}&2\,d_{{1}} \, d_{{3}}&
-6\,d_{{4}} \, d_{{1}}&2\,d_{{1}} \, d_{{5}}&d_{{2}}^{2}-6\,d_{{3}}^{2}+6
\,d_{{4}}^{2}-d_{{5}}^{2}&2\,d_{{2}} \, d_{{3}}\\ \noalign{\medskip}0&
d_{{5}} \, d_{{4}}&0&0&-d_{{4}} \, d_{{2}}&0&0\\ \noalign{\medskip}0&d_{{3}} \, d_
{{5}}&0&0&-d_{{2}} \, d_{{3}}&0&-d_{{1}} \, d_{{5}}\\ \noalign{\medskip}0&d_{{
3}} \, d_{{4}}&0&-d_{{2}} \, d_{{3}}&0&0&-d_{{4}} \, d_{{1}}\\ \noalign{\medskip}0
&d_{{3}}^{2}-d_{{4}}^{2}&-d_{{2}} \, d_{{3}}&d_{{4}} \, d_{{2}}&0&0&-d_{{1
}} \, d_{{3}}\\ \noalign{\medskip}0&d_{{2}} \, d_{{5}}&0&-d_{{5}} \, d_{{4}}&-{d_{
{2}}}^{2}+d_{{4}}^{2}&-d_{{1}} \, d_{{5}}&0\\ \noalign{\medskip}0&d_{{2}
} \, d_{{3}}&-d_{{2}}^{2}+d_{{4}}^{2}&-d_{{3}} \, d_{{4}}&0&-d_{{1}} \, d_{{3}
}&d_{{1}} \, d_{{2}}\\ \noalign{\medskip}0&-d_{{1}} \, d_{{3}}&-d_{{1}} \, d_{{2}}
&0&0&-d_{{2}} \, d_{{3}}&d_{{1}}^{2}-2\,d_{{4}}^{2}+d_{{5}}^{2}
\\ \noalign{\medskip}0&0&d_{{5}} \, d_{{4}}&0&-d_{{3}} \, d_{{4}}&0&0
\\ \noalign{\medskip}0&0&d_{{3}} \, d_{{5}}&-d_{{5}} \, d_{{4}}&-d_{{3}}^{2}
+d_{{4}}^{2}&0&0\\ \noalign{\medskip}0&0&d_{{2}} \, d_{{5}}&0&-d_{{2}} \, d_
{{3}}&0&-d_{{1}} \, d_{{5}}\\ \noalign{\medskip}0&0&d_{{4}} \, d_{{2}}&-d_{{2}
} \, d_{{3}}&0&0&-d_{{4}} \, d_{{1}}\\ \noalign{\medskip}0&0&0&d_{{3}} \, d_{{5}}&
-d_{{3}} \, d_{{4}}&0&0\\ \noalign{\medskip}0&0&0&d_{{2}} \, d_{{5}}&-d_{{4}}d
_{{2}}&0&0\\ \noalign{\medskip}0&0&0&0&0&d_{{5}} \, d_{{4}}&0
\\ \noalign{\medskip}0&0&0&0&0&d_{{3}} \, d_{{5}}&-d_{{2}} \, d_{{5}}
\\ \noalign{\medskip}0&0&0&0&0&0&d_{{5}} \, d_{{4}}\\ \noalign{\medskip}0&0
&0&0&0&0&d_{{3}} \, d_{{5}}\\ \noalign{\medskip}0&0&0&0&0&0&0\end {array}
 \right),$$}
{\footnotesize $$\hspace{-1cm} 
S_B=\left( \begin {array}{ccccccc} -d_{{3}}^{2}+d_{{5}}^{2}&-d_{{5}}d
_{{4}}&-2\,d_{{4}} \, d_{{2}}&d_{{2}} \, d_{{3}}&0&-d_{{4}} \, d_{{2}}&-d_{{2}} \, d_{
{5}}\\ \noalign{\medskip}-d_{{3}} \, d_{{4}}&d_{{3}} \, d_{{5}}&d_{{2}} \, d_{{3}}
&-d_{{4}} \, d_{{2}}&d_{{2}} \, d_{{5}}&2\,d_{{2}} \, d_{{3}}&0
\\ \noalign{\medskip}0&0&-2\,d_{{4}} \, d_{{1}}&d_{{1}} \, d_{{3}}&0&-d_{{4}}d
_{{1}}&-d_{{1}} \, d_{{5}}\\ \noalign{\medskip}0&0&d_{{1}} \, d_{{3}}&-d_{{4}}
d_{{1}}&d_{{1}} \, d_{{5}}&2\,d_{{1}} \, d_{{3}}&0\\ \noalign{\medskip}-d_{{2}
} \, d_{{3}}&0&d_{{3}} \, d_{{4}}&d_{{2}}^{2}-d_{{5}}^{2}&d_{{5}} \, d_{{4}}&d
_{{3}} \, d_{{4}}&d_{{3}} \, d_{{5}}\\ \noalign{\medskip}2\,d_{{4}} \, d_{{2}}&0&{
d_{{2}}}^{2}+d_{{3}}^{2}-d_{{4}}^{2}-d_{{5}}^{2}&0&2\,d_{{3}} \, d_{
{5}}&-d_{{2}}^{2}+d_{{3}}^{2}-d_{{4}}^{2}+d_{{5}}^{2}&-2\,d_{{
5}} \, d_{{4}}\\ \noalign{\medskip}-4\,d_{{4}} \, d_{{2}}&2\,d_{{2}} \, d_{{5}}&-2
\,d_{{3}}^{2}+2\,d_{{5}}^{2}&2\,d_{{3}} \, d_{{4}}&-4\,d_{{3}} \, d_{{5}}&
3\,d_{{2}}^{2}-3\,d_{{3}}^{2}+d_{{4}}^{2}-d_{{5}}^{2}&2\,d_{{5
}} \, d_{{4}}\\ \noalign{\medskip} d_{{1}} \, d_{{2}}&0&-d_{{4}} \, d_{{1}}&0&0&-d_
{{4}} \, d_{{1}}&-d_{{1}} \, d_{{5}}\\ \noalign{\medskip}-2\,d_{{4}} \, d_{{1}}&d_
{{1}} \, d_{{5}}&0&0&0&3\,d_{{1}} \, d_{{2}}&0\\ \noalign{\medskip}0&0&-d_{{5}
} \, d_{{4}}&d_{{3}} \, d_{{5}}&d_{{3}} \, d_{{4}}&0&d_{{1}}^{2}-d_{{3}}^{2}
\\ \noalign{\medskip}0&0&0&d_{{5}} \, d_{{4}}&d_{{1}}^{2}-d_{{4}}^{2}&
-d_{{3}} \, d_{{5}}&d_{{3}} \, d_{{4}}\\ \noalign{\medskip}0&0&d_{{3}} \, d_{{4}}&
d_{{1}}^{2}-d_{{5}}^{2}&d_{{5}} \, d_{{4}}&d_{{3}} \, d_{{4}}&d_{{3}} \, d_{{5
}}\\ \noalign{\medskip}0&0&d_{{1}}^{2}+d_{{3}}^{2}-d_{{4}}^{2}-{
d_{{5}}}^{2}&0&2\,d_{{3}} \, d_{{5}}&-d_{{1}}^{2}+d_{{3}}^{2}-{d_{{4}}
}^{2}+d_{{5}}^{2}&-2\,d_{{5}} \, d_{{4}}\\ \noalign{\medskip}d_{{5}} \, d_{{
4}}&d_{{1}}^{2}-d_{{4}}^{2}&0&0&0&-d_{{2}} \, d_{{5}}&d_{{4}} \, d_{{2}}
\\ \noalign{\medskip}d_{{1}}^{2}-d_{{5}}^{2}&d_{{5}} \, d_{{4}}&d_{{4}
} \, d_{{2}}&0&0&d_{{4}} \, d_{{2}}&d_{{2}} \, d_{{5}}\\ \noalign{\medskip}6\,d_{{
4}} \, d_{{2}}&-2\,d_{{2}} \, d_{{5}}&3\,d_{{3}}^{2}-3\,d_{{5}}^{2}&-4\,d_
{{3}} \, d_{{4}}&6\,d_{{3}} \, d_{{5}}&-d_{{1}}^{2}-4\,d_{{2}}^{2}+5\,{d_{
{3}}}^{2}-2\,d_{{4}}^{2}+2\,d_{{5}}^{2}&-4\,d_{{5}} \, d_{{4}}
\\ \noalign{\medskip}-6\,d_{{4}} \, d_{{2}}&2\,d_{{2}} \, d_{{5}}&-5\,{d_{{3}}
}^{2}+5\,d_{{5}}^{2}&6\,d_{{3}} \, d_{{4}}&-10\,d_{{3}} \, d_{{5}}&4\,{d_{{1
}}}^{2}+4\,d_{{2}}^{2}-8\,d_{{3}}^{2}+3\,d_{{4}}^{2}-3\,{d_{{5}}
}^{2}&6\,d_{{5}} \, d_{{4}}\\ \noalign{\medskip}-d_{{1}} \, d_{{5}}&0&0&0&0&0&
d_{{1}} \, d_{{2}}\\ \noalign{\medskip}0&0&0&0&d_{{1}} \, d_{{2}}&0&0
\\ \noalign{\medskip}0&0&0&d_{{1}} \, d_{{2}}&0&0&0\\ \noalign{\medskip}d_
{{4}} \, d_{{1}}&0&d_{{1}} \, d_{{2}}&0&0&-d_{{1}} \, d_{{2}}&0
\\ \noalign{\medskip}0&d_{{1}} \, d_{{2}}&0&0&0&d_{{1}} \, d_{{5}}&-d_{{4}} \, d_{
{1}}\\ \noalign{\medskip}0&0&0&-d_{{4}} \, d_{{1}}&0&d_{{1}} \, d_{{3}}&0
\\ \noalign{\medskip}2\,d_{{3}} \, d_{{4}}&-d_{{3}} \, d_{{5}}&-d_{{2}} \, d_{{3}}
&2\,d_{{4}} \, d_{{2}}&-d_{{2}} \, d_{{5}}&-3\,d_{{2}} \, d_{{3}}&0
\\ \noalign{\medskip}0&0&0&-d_{{1}} \, d_{{5}}&0&0&d_{{1}} \, d_{{3}}
\\ \noalign{\medskip}0&0&-d_{{1}} \, d_{{5}}&0&d_{{1}} \, d_{{3}}&d_{{1}} \, d_{{5
}}&-d_{{4}} \, d_{{1}}\\ \noalign{\medskip}0&d_{{1}} \, d_{{3}}&0&0&0&0&0
\\ \noalign{\medskip}d_{{1}} \, d_{{3}}&0&0&0&0&0&0\\ \noalign{\medskip}0&0
&0&-d_{{1}} \, d_{{5}}&d_{{4}} \, d_{{1}}&0&0\\ \noalign{\medskip}-d_{{1}} \, d_{{
5}}&d_{{4}} \, d_{{1}}&0&0&0&0&0\\ \noalign{\medskip}-d_{{2}} \, d_{{5}}&-d_{{
4}} \, d_{{2}}&-d_{{5}} \, d_{{4}}&d_{{3}} \, d_{{5}}&d_{{3}} \, d_{{4}}&0&{d_{{2}}}^{
2}-d_{{3}}^{2}\\ \noalign{\medskip}0&-d_{{2}} \, d_{{3}}&0&d_{{5}} \, d_{{4}
}&d_{{2}}^{2}-d_{{4}}^{2}&-d_{{3}} \, d_{{5}}&d_{{3}} \, d_{{4}}
\\ \noalign{\medskip}0&-d_{{3}} \, d_{{4}}&0&-d_{{2}} \, d_{{5}}&0&0&d_{{2}} \, d_
{{3}}\\ \noalign{\medskip}-d_{{5}} \, d_{{4}}&-d_{{3}}^{2}+d_{{5}}^{2}
&-d_{{2}} \, d_{{5}}&0&d_{{2}} \, d_{{3}}&d_{{2}} \, d_{{5}}&-d_{{4}} \, d_{{2}}
\\ \noalign{\medskip}d_{{3}} \, d_{{5}}&-d_{{3}} \, d_{{4}}&0&-d_{{2}} \, d_{{5}}&
d_{{4}} \, d_{{2}}&0&0\end {array} \right).$$}

Similarly, we can prove that we have $\ker_D(S_5)={\rm im}_D(T_5)$, where the matrix $T_5=(T_A \quad T_B)  \in D^{35 \times 35}$ is defined by:

\newpage

\begin{sideways}
{\tiny $\begin{array}{c}
 T_A=\left( \begin {array}{cccccccccccccccccccc} d_{{5}}&0&0&0&0&0&0&0&0&0
&0&0&0&-d_{{4}}&d_{{5}}&0&0&0&0&0\\ \noalign{\medskip}d_{{5}}&0&0&0&0&0
&0&0&0&-d_{{2}}&0&0&0&0&d_{{5}}&0&0&d_{{1}}&0&0\\ \noalign{\medskip}d_
{{4}}&-3\,d_{{3}}&0&0&0&0&0&0&d_{{1}}&0&0&0&0&-d_{{5}}&2\,d_{{4}}&-d_{
{2}}&-d_{{2}}&0&0&0\\ \noalign{\medskip}-d_{{2}}&0&0&0&d_{{3}}&-d_{{4}
}&d_{{4}}&d_{{1}}&0&d_{{5}}&0&0&d_{{4}}&0&-d_{{2}}&2\,d_{{4}}&d_{{4}}&0
&0&0\\ \noalign{\medskip}-d_{{2}}&0&d_{{1}}&0&d_{{3}}&-2\,d_{{4}}&-d_{
{4}}&0&0&d_{{5}}&0&-d_{{3}}&2\,d_{{4}}&0&0&d_{{4}}&d_{{4}}&0&0&0
\\ \noalign{\medskip}d_{{1}}&0&0&0&0&0&0&0&d_{{4}}&0&0&0&0&0&0&0&0&d_{
{5}}&0&-d_{{3}}\\ \noalign{\medskip}0&2\,d_{{4}}&0&0&0&0&0&0&0&0&0&0&0
&0&-d_{{3}}&0&0&0&0&0\\ \noalign{\medskip}0&d_{{5}}&0&0&0&0&0&0&0&0&0&0
&0&-d_{{3}}&0&0&0&0&0&0\\ \noalign{\medskip}0&-2\,d_{{2}}&0&0&0&d_{{3}
}&2\,d_{{3}}&0&0&0&0&d_{{4}}&0&0&0&d_{{3}}&0&0&0&0
\\ \noalign{\medskip}0&2\,d_{{3}}&0&0&0&0&0&0&0&0&0&0&-d_{{2}}&0&-d_{{
4}}&0&0&0&0&0\\ \noalign{\medskip}0&2\,d_{{4}}&0&0&0&0&0&0&0&0&0&-d_{{
2}}&0&0&0&0&0&0&0&d_{{1}}\\ \noalign{\medskip}0&d_{{5}}&0&0&0&0&0&0&0&0
&-d_{{2}}&0&0&0&0&0&0&0&d_{{1}}&0\\ \noalign{\medskip}0&-d_{{2}}&0&d_{
{1}}&-d_{{4}}&d_{{3}}&d_{{3}}&0&0&0&-d_{{5}}&d_{{4}}&-d_{{3}}&0&0&-d_{
{3}}&-d_{{3}}&0&0&0\\ \noalign{\medskip}0&d_{{1}}&0&0&0&0&0&0&-d_{{3}}
&0&0&0&0&0&0&0&0&0&-d_{{5}}&d_{{4}}\\ \noalign{\medskip}0&0&-d_{{5}}&0
&0&0&0&d_{{5}}&0&0&0&0&0&0&0&0&0&0&0&0\\ \noalign{\medskip}0&0&d_{{2}}
&0&0&0&0&0&d_{{4}}&0&0&0&0&0&0&0&0&d_{{5}}&0&-d_{{3}}
\\ \noalign{\medskip}0&0&-d_{{4}}&2\,d_{{3}}&0&0&d_{{1}}&2\,d_{{4}}&-d
_{{2}}&0&0&0&0&0&0&0&0&0&0&0\\ \noalign{\medskip}0&0&0&-d_{{4}}&0&0&0&
-d_{{3}}&0&0&0&0&0&0&0&0&0&0&0&0\\ \noalign{\medskip}0&0&0&d_{{2}}&0&0
&0&0&-d_{{3}}&0&0&0&0&0&0&0&0&0&-d_{{5}}&d_{{4}}\\ \noalign{\medskip}0
&0&0&-d_{{3}}&0&d_{{1}}&0&-d_{{4}}&0&0&0&0&0&0&0&0&0&0&0&0
\\ \noalign{\medskip}0&0&0&0&-d_{{5}}&0&0&0&0&0&-d_{{4}}&d_{{5}}&0&0&0
&0&0&0&0&0\\ \noalign{\medskip}0&0&0&0&-d_{{5}}&0&0&0&0&-d_{{3}}&0&d_{
{5}}&0&0&0&0&0&0&0&0\\ \noalign{\medskip}0&0&0&0&d_{{1}}&0&0&d_{{3}}&0
&0&0&0&0&0&0&0&0&0&0&-d_{{2}}\\ \noalign{\medskip}0&0&0&0&0&-d_{{5}}&0
&0&0&d_{{4}}&-d_{{3}}&0&d_{{5}}&0&0&0&0&0&0&0\\ \noalign{\medskip}0&0&0
&0&0&-d_{{5}}&d_{{5}}&0&0&d_{{4}}&0&0&0&-d_{{2}}&0&d_{{5}}&0&0&0&0
\\ \noalign{\medskip}0&0&0&0&0&0&0&0&0&0&0&0&0&0&0&0&0&d_{{4}}&-d_{{3}
}&0\\ \noalign{\medskip}0&0&0&0&0&0&0&0&0&0&0&0&0&0&0&0&0&-d_{{3}}&0&d
_{{5}}\\ \noalign{\medskip}0&0&0&0&0&0&0&0&0&0&0&0&0&0&0&0&0&d_{{4}}&-
d_{{3}}&0\\ \noalign{\medskip}0&0&0&0&0&0&0&0&0&0&0&0&0&0&0&0&0&-d_{{3
}}&0&d_{{5}}\\ \noalign{\medskip}0&0&0&0&0&0&0&0&0&0&0&0&0&0&0&0&0&-d_
{{2}}&0&0\\ \noalign{\medskip}0&0&0&0&0&0&0&0&0&0&0&0&0&0&0&0&0&0&-d_{
{4}}&d_{{5}}\\ \noalign{\medskip}0&0&0&0&0&0&0&0&0&0&0&0&0&0&0&0&0&0&-
d_{{4}}&d_{{5}}\\ \noalign{\medskip}0&0&0&0&0&0&0&0&0&0&0&0&0&0&0&0&0&0
&-d_{{2}}&0\\ \noalign{\medskip}0&0&0&0&0&0&0&0&0&0&0&0&0&0&0&0&0&0&0&0
\\ \noalign{\medskip}0&0&0&0&0&0&0&0&0&0&0&0&0&0&0&0&0&0&0&0
\end {array} \right)
\end{array}$}
\end{sideways}

{\footnotesize
$T_B=\left( \begin {array}{ccccccccccccccc} 0&0&0&0&0&0&0&0&0&d_{{1}}&0&0&0
&-d_{{4}}&d_{{3}}\\ \noalign{\medskip}0&0&0&0&0&0&0&0&0&0&0&0&-d_{{3}}
&0&d_{{3}}\\ \noalign{\medskip}0&0&0&-d_{{3}}&0&0&0&0&0&0&0&0&0&-2\,d_
{{5}}&0\\ \noalign{\medskip}0&0&0&0&0&0&0&0&0&0&0&0&0&0&0
\\ \noalign{\medskip}0&0&0&0&0&0&0&0&0&0&-d_{{5}}&0&0&0&0
\\ \noalign{\medskip}2\,d_{{4}}&0&0&0&0&0&0&d_{{3}}&0&0&0&0&0&0&0
\\ \noalign{\medskip}0&0&0&d_{{4}}&0&0&0&d_{{1}}&0&0&0&0&d_{{5}}&0&-d_
{{5}}\\ \noalign{\medskip}0&0&0&d_{{5}}&0&0&d_{{1}}&0&0&0&0&0&d_{{4}}&0
&0\\ \noalign{\medskip}0&0&d_{{1}}&-d_{{2}}&0&0&0&0&0&0&0&d_{{5}}&0&0&0
\\ \noalign{\medskip}d_{{1}}&0&0&d_{{3}}&0&0&0&0&0&0&0&0&0&d_{{5}}&0
\\ \noalign{\medskip}0&0&0&d_{{4}}&0&0&0&0&0&0&0&0&d_{{5}}&0&0
\\ \noalign{\medskip}0&0&0&d_{{5}}&0&0&0&0&0&0&0&0&d_{{4}}&0&-d_{{4}}
\\ \noalign{\medskip}0&0&0&0&0&0&0&0&0&0&0&d_{{5}}&0&0&0
\\ \noalign{\medskip}-d_{{3}}&0&0&0&0&0&0&0&0&0&0&0&0&0&0
\\ \noalign{\medskip}0&-d_{{4}}&0&0&0&d_{{4}}&0&0&-d_{{3}}&d_{{2}}&0&0
&0&0&0\\ \noalign{\medskip}2\,d_{{4}}&0&0&0&0&0&0&0&0&-d_{{5}}&0&0&0&0
&0\\ \noalign{\medskip}0&-d_{{5}}&-d_{{3}}&0&0&2\,d_{{5}}&0&0&0&0&0&0&0
&0&0\\ \noalign{\medskip}0&0&d_{{4}}&0&-d_{{5}}&0&0&d_{{2}}&d_{{5}}&0&0
&0&0&0&0\\ \noalign{\medskip}-d_{{3}}&0&0&0&0&0&d_{{5}}&-d_{{4}}&0&0&0
&0&0&0&0\\ \noalign{\medskip}-d_{{2}}&0&d_{{3}}&0&0&-d_{{5}}&0&0&0&0&0
&0&0&0&0\\ \noalign{\medskip}0&0&0&0&0&0&0&0&d_{{1}}&0&0&d_{{4}}&0&0&-
d_{{2}}\\ \noalign{\medskip}0&0&0&0&d_{{1}}&0&0&0&0&0&d_{{3}}&0&-d_{{2
}}&0&0\\ \noalign{\medskip}0&0&0&0&0&0&0&0&-d_{{5}}&0&0&0&0&0&0
\\ \noalign{\medskip}0&0&0&0&0&d_{{1}}&0&0&0&0&-d_{{4}}&d_{{3}}&0&-d_{
{2}}&0\\ \noalign{\medskip}0&d_{{1}}&0&0&0&0&0&0&0&0&0&d_{{3}}&0&-d_{{
2}}&0\\ \noalign{\medskip}d_{{5}}&0&0&0&0&d_{{2}}&0&0&0&0&0&0&0&0&0
\\ \noalign{\medskip}0&0&0&0&d_{{2}}&0&0&-d_{{5}}&0&0&0&0&0&0&0
\\ \noalign{\medskip}d_{{5}}&0&0&0&0&0&d_{{3}}&0&0&-d_{{4}}&0&0&0&d_{{
1}}&0\\ \noalign{\medskip}0&0&0&0&0&0&d_{{4}}&-d_{{5}}&0&0&0&0&d_{{1}}
&0&0\\ \noalign{\medskip}0&d_{{4}}&0&0&d_{{3}}&-d_{{4}}&0&0&0&0&d_{{1}
}&0&0&0&0\\ \noalign{\medskip}0&0&0&0&0&0&0&0&d_{{2}}&0&0&0&0&0&0
\\ \noalign{\medskip}0&0&0&0&0&0&d_{{4}}&-d_{{5}}&0&0&0&0&0&0&d_{{1}}
\\ \noalign{\medskip}0&d_{{3}}&0&0&0&0&0&0&d_{{4}}&0&0&d_{{1}}&0&0&0
\\ \noalign{\medskip}0&-d_{{3}}&d_{{5}}&0&-d_{{4}}&0&d_{{2}}&0&0&0&0&0
&0&0&0\\ \noalign{\medskip}0&0&0&0&0&0&-d_{{4}}&d_{{5}}&0&d_{{3}}&0&0&0
&0&0\end {array} \right)$.}

\bigskip

We have $\ker_D(T_5)={\rm im}_D(U_5)$, where $U_5=(U_A \quad U_B \quad U_C) \in D^{14 \times 35}$ and the matrices $U_A$, $U_B$ and $U_C$ are respectively defined by:
\newpage
\begin{sideways}
{\tiny  $\hspace*{2cm} 
\begin{array}{c}
U_A=\left( 
\begin {array}{ccccccccccccccc} -d_{{2}} \, d_{{3}}&0&0&0&-d_{{3}}
d_{{5}}&0&-d_{{2}} \, d_{{5}}&d_{{4}} \, d_{{2}}&0&0&0&0&-d_{{5}} \, d_{{4}}&0&-d_
{{1}} \, d_{{3}}\\ \noalign{\medskip}-d_{{3}} \, d_{{5}}&d_{{3}} \, d_{{5}}&0&d_{{
2}} \, d_{{3}}&-d_{{2}} \, d_{{3}}&0&-{d_{{2}}}^{2}+{d_{{4}}}^{2}&d_{{5}} \, d_{{4
}}&-d_{{4}} \, d_{{2}}&-d_{{3}} \, d_{{4}}&{d_{{3}}}^{2}-{d_{{4}}}^{2}&-d_{{5}
} \, d_{{4}}&0&0&0\\ \noalign{\medskip}-{d_{{2}}}^{2}+{d_{{5}}}^{2}&{d_{{4
}}}^{2}-{d_{{5}}}^{2}&-d_{{5}} \, d_{{4}}&-d_{{2}} \, d_{{5}}&0&0&0&0&0&-d_{{5
}} \, d_{{4}}&-d_{{3}} \, d_{{5}}&d_{{3}} \, d_{{4}}&0&0&d_{{1}} \, d_{{2}}
\\ \noalign{\medskip}-d_{{1}} \, d_{{3}}&0&0&0&0&d_{{3}} \, d_{{5}}&-d_{{1}} \, d_
{{5}}&d_{{4}} \, d_{{1}}&0&0&0&0&0&d_{{5}} \, d_{{4}}&0\\ \noalign{\medskip}d_
{{1}} \, d_{{5}}&-d_{{1}} \, d_{{5}}&-d_{{4}} \, d_{{1}}&0&-d_{{1}} \, d_{{2}}&-{d_{{2
}}}^{2}+{d_{{4}}}^{2}&0&0&0&-2\,d_{{4}} \, d_{{1}}&d_{{1}} \, d_{{3}}&0&0&-d_{
{3}} \, d_{{4}}&0\\ \noalign{\medskip}-d_{{1}} \, d_{{2}}&0&0&-d_{{1}} \, d_{{5}}&
d_{{1}} \, d_{{5}}&d_{{2}} \, d_{{5}}&0&0&0&0&0&0&0&0&{d_{{1}}}^{2}-{d_{{4}}}^
{2}\\ \noalign{\medskip}0&d_{{3}} \, d_{{4}}&-d_{{3}} \, d_{{5}}&0&0&0&d_{{5}}
d_{{4}}&-{d_{{2}}}^{2}+{d_{{5}}}^{2}&-d_{{2}} \, d_{{5}}&-2\,d_{{3}} \, d_{{5}
}&-d_{{5}} \, d_{{4}}&{d_{{3}}}^{2}-{d_{{5}}}^{2}&d_{{2}} \, d_{{5}}&0&0
\\ \noalign{\medskip}0&d_{{1}} \, d_{{3}}&0&0&0&-d_{{3}} \, d_{{5}}&d_{{1}} \, d_{
{5}}&0&0&0&-d_{{1}} \, d_{{5}}&0&0&0&0\\ \noalign{\medskip}0&-d_{{4}} \, d_{{1
}}&0&0&0&d_{{5}} \, d_{{4}}&0&0&0&-d_{{1}} \, d_{{5}}&0&d_{{1}} \, d_{{3}}&0&d_{{3
}} \, d_{{5}}&-d_{{4}} \, d_{{2}}\\ \noalign{\medskip}0&0&d_{{1}} \, d_{{3}}&0&0&-
d_{{3}} \, d_{{4}}&d_{{4}} \, d_{{1}}&-d_{{1}} \, d_{{5}}&0&d_{{1}} \, d_{{3}}&-d_{{4}
} \, d_{{1}}&d_{{1}} \, d_{{5}}&-d_{{1}} \, d_{{2}}&-{d_{{2}}}^{2}+{d_{{3}}}^{2}&0
\\ \noalign{\medskip}0&0&0&d_{{1}} \, d_{{3}}&0&d_{{2}} \, d_{{3}}&-d_{{1}} \, d_{
{2}}&0&-d_{{4}} \, d_{{1}}&0&0&0&d_{{4}} \, d_{{1}}&d_{{4}} \, d_{{2}}&0
\\ \noalign{\medskip}0&0&0&0&0&0&0&-d_{{1}} \, d_{{2}}&-d_{{1}} \, d_{{5}}&0&0
&0&0&-d_{{2}} \, d_{{5}}&-d_{{3}} \, d_{{4}}\\ \noalign{\medskip}0&0&0&0&0&0&0
&0&0&0&-d_{{1}} \, d_{{5}}&d_{{4}} \, d_{{1}}&0&d_{{5}} \, d_{{4}}&0
\\ \noalign{\medskip}0&0&0&0&0&0&0&0&0&0&0&0&0&0&-d_{{2}} \, d_{{3}}
\end {array} 
\right), \vspace{1cm} \\
U_B=\left(  \begin {array}{cccccccccccccc} 0&0&-d_{{1}} \, d_{{5}}&0&0&-{d_{{3
}}}^{2}+{d_{{5}}}^{2}&{d_{{4}}}^{2}-{d_{{5}}}^{2}&0 &d_{{3}} \, d_{{4}}&0&0
&0&0&0\\ \noalign{\medskip} \, d_{{1}} \, d_{{3}}&0&d_{{1}} \, d_{{2}}&d_{{4}} \, d_{{
1}}&0&d_{{2}} \, d_{{5}}&-d_{{2}} \, d_{{5}}&0&0&0&0&2\,d_{{1}} \, d_{{5}}&0&0
\\ \noalign{\medskip}d_{{1}} \, d_{{5}}&0&0&0&0&-d_{{2}} \, d_{{3}}&0&0&0& d_{{
4}} \, d_{{2}}&-d_{{4}} \, d_{{1}}&0&0&0\\ \noalign{\medskip}0&0&0&0&0&0&0&0&0
&0&0&0&-d_{{3}} \, d_{{4}}&-{d_{{4}}}^{2}+{d_{{5}}}^{2}
\\ \noalign{\medskip}{d_{{1}}}^{2}-{d_{{4}}}^{2}&-d_{{4}} \, d_{{2}}&d_{{2
}} \, d_{{3}}&d_{{3}} \, d_{{4}}&-2\,d_{{4}} \, d_{{2}}&0&0&d_{{2}} \, d_{{3}}&0&0&-d_
{{5}} \, d_{{4}}&2\,d_{{3}} \, d_{{5}}&d_{{5}} \, d_{{4}}&-2\,d_{{3}} \, d_{{5}}
\\ \noalign{\medskip}0&d_{{5}} \, d_{{4}}&d_{{3}} \, d_{{5}}&0&d_{{5}} \, d_{{4}}&
d_{{1}} \, d_{{3}}&0&d_{{3}} \, d_{{5}}&-d_{{4}} \, d_{{1}}&d_{{4}} \, d_{{1}}&0&0&-d_
{{4}} \, d_{{2}}&0\\ \noalign{\medskip}0&0&0&-d_{{1}} \, d_{{5}}&0&0&-d_{{4}}d
_{{2}}&0&-d_{{2}} \, d_{{3}}&d_{{2}} \, d_{{3}}&d_{{1}} \, d_{{3}}&2\,d_{{4}} \, d_{{1
}}&0&0\\ \noalign{\medskip}0&0&-d_{{2}} \, d_{{5}}&-d_{{5}} \, d_{{4}}&0&0&-d_
{{1}} \, d_{{2}}&-d_{{2}} \, d_{{5}}&0&0&d_{{3}} \, d_{{4}}&{d_{{1}}}^{2}-{d_{{3}}
}^{2}+{d_{{4}}}^{2}-{d_{{5}}}^{2}&0&-{d_{{2}}}^{2}+{d_{{3}}}^{2}
\\ \noalign{\medskip}-d_{{5}} \, d_{{4}}&0&0&-d_{{3}} \, d_{{5}}&-d_{{2}} \, d_{{5
}}&0&0&0&-d_{{1}} \, d_{{2}}&0&{d_{{1}}}^{2}-{d_{{5}}}^{2}&2\,d_{{3}} \, d_{{4
}}&-{d_{{2}}}^{2}+{d_{{5}}}^{2}&-2\,d_{{3}} \, d_{{4}}
\\ \noalign{\medskip}d_{{3}} \, d_{{4}}&d_{{2}} \, d_{{3}}&0&{d_{{1}}}^{2}-{d_
{{3}}}^{2}&d_{{2}} \, d_{{3}}&0&0&-d_{{4}} \, d_{{2}}&0&0&-d_{{3}} \, d_{{5}}&0&d_
{{3}} \, d_{{5}}&0\\ \noalign{\medskip}0&0&{d_{{1}}}^{2}-{d_{{3}}}^{2}&0&d
_{{3}} \, d_{{4}}&-d_{{1}} \, d_{{5}}&d_{{1}} \, d_{{5}}&-{d_{{3}}}^{2}+{d_{{4}}}^
{2}&0&0&0&0&0&2\,d_{{2}} \, d_{{5}}\\ \noalign{\medskip}0&d_{{3}} \, d_{{5}}&0
&0&2\,d_{{3}} \, d_{{5}}&0&d_{{4}} \, d_{{1}}&d_{{5}} \, d_{{4}}&0&d_{{1}} \, d_{{3}}&0
&0&d_{{2}} \, d_{{3}}&2\,d_{{4}} \, d_{{2}}\\ \noalign{\medskip}0&0&-d_{{2}} \, d_
{{5}}&-d_{{5}} \, d_{{4}}&0&-d_{{1}} \, d_{{2}}&0&-d_{{2}} \, d_{{5}}&0&0&0&{d_{{4
}}}^{2}-{d_{{5}}}^{2}&0&-{d_{{4}}}^{2}+{d_{{5}}}^{2}
\\ \noalign{\medskip}-d_{{3}} \, d_{{5}}&0&-d_{{2}} \, d_{{5}}&-d_{{5}} \, d_{{4}}
&0&0&0&0&0&0&d_{{3}} \, d_{{4}}&{d_{{4}}}^{2}-{d_{{5}}}^{2}&0&0
\end {array} \right), 
\end{array}$
} 
\end{sideways}
and:
$$U_C=\left( \begin {array}{cccccc} 0&0&0&0&d_{{4}} \, d_{{1}}&2\,d_{{1}} \, d_{{2}
}\\ \noalign{\medskip}0&-2\,d_{{1}} \, d_{{5}}&0&0&0&2\,d_{{1}} \, d_{{5}}
\\ \noalign{\medskip}0&2\,d_{{1}} \, d_{{3}}&0&0&0&0\\ \noalign{\medskip}0
&0&{d_{{3}}}^{2}-{d_{{5}}}^{2}&0&0&{d_{{1}}}^{2}-{d_{{4}}}^{2}
\\ \noalign{\medskip}-d_{{2}} \, d_{{5}}&0&0&0&0&0\\ \noalign{\medskip}-{d
_{{4}}}^{2}+{d_{{5}}}^{2}&0&2\,d_{{2}} \, d_{{3}}&-d_{{3}} \, d_{{4}}&-d_{{3}}
d_{{4}}&0\\ \noalign{\medskip}0&0&0&0&d_{{1}} \, d_{{2}}&0
\\ \noalign{\medskip}d_{{2}} \, d_{{3}}&0&-{d_{{3}}}^{2}+{d_{{5}}}^{2}&0&d
_{{4}} \, d_{{2}}&0\\ \noalign{\medskip}-d_{{4}} \, d_{{2}}&-2\,d_{{3}} \, d_{{4}}
&2\,d_{{3}} \, d_{{4}}&d_{{2}} \, d_{{3}}&d_{{2}} \, d_{{3}}&0
\\ \noalign{\medskip}0&-2\,d_{{5}} \, d_{{4}}&2\,d_{{5}} \, d_{{4}}&d_{{2}} \, d_{
{5}}&0&2\,d_{{5}} \, d_{{4}}\\ \noalign{\medskip}-d_{{3}} \, d_{{5}}&0&-2\,d_{
{2}} \, d_{{5}}&d_{{5}} \, d_{{4}}&0&0\\ \noalign{\medskip}-d_{{3}} \, d_{{4}}&0&0
&-{d_{{3}}}^{2}+{d_{{5}}}^{2}&{d_{{1}}}^{2}-{d_{{3}}}^{2}&2\,d_{{4}} \, d_
{{2}}\\ \noalign{\medskip}0&{d_{{1}}}^{2}-{d_{{4}}}^{2}&-{d_{{2}}}^{2}
+{d_{{4}}}^{2}&d_{{4}} \, d_{{2}}&d_{{4}} \, d_{{2}}&0\\ \noalign{\medskip}0&-
{d_{{3}}}^{2}+{d_{{5}}}^{2}&0&0&d_{{4}} \, d_{{2}}&{d_{{2}}}^{2}-{d_{{5}}}
^{2}\end {array} \right).$$

Moreover, we have $\ker_D(U_5)={\rm im}_D(V_5)$, where $V_5 \in D^{5 \times 14}$ is the full row rank matrix defined by:
$$V_5=\left( \begin {array}{cccccccccccccc} -d_{{2}}&d_{{5}}&d_{{3}}&0&0&0&
-d_{{4}}&0&0&0&0&0&0&2\,d_{{1}}\\ \noalign{\medskip}d_{{1}}&0&0&-2\,d_
{{2}}&0&d_{{3}}&0&0&0&0&d_{{5}}&-d_{{4}}&0&0\\ \noalign{\medskip}0&d_{
{1}}&0&-2\,d_{{5}}&-d_{{3}}&0&0&-2\,d_{{5}}&0&-d_{{4}}&-d_{{2}}&0&2\,d
_{{5}}&0\\ \noalign{\medskip}0&0&d_{{1}}&0&-d_{{5}}&-d_{{2}}&0&0&d_{{4
}}&0&0&0&-2\,d_{{3}}&0\\ \noalign{\medskip}0&0&0&0&0&0&d_{{1}}&-2\,d_{
{4}}&-d_{{3}}&d_{{5}}&0&-d_{{2}}&0&2\,d_{{4}}\end {array} \right).$$ 

Hence, the $D$-module $M=D^5/(D^{14} \, R_5)=D^5 /(D^{15} \, W_5)$ admits the following finite free resolution:\\
\[ 0 \rightarrow  D^5 \stackrel{V_5}{\longrightarrow} D^{14} \stackrel{U_5}{\longrightarrow} D^{35} \stackrel{T_5}{\longrightarrow} D^{35} \stackrel{S_5}{\longrightarrow}
 D^{14} \stackrel{R_5}{\longrightarrow} D^5 \rightarrow M_5 \rightarrow 0 \]

\noindent
{\bf REFERENCES FOR THE APPENDIX}\\

\noindent
[1] F. CHYZAK, A. QUADRAT and D. ROBERTZ:  Effective algorithms for parametrizing linear control systems over Ore algebras,
Appl. Algebra Engrg. Comm. Comput., 16, 319-376, 2005. \\
\noindent
[2] F. CHYZAK, A. QUADRAT and D. ROBERTZ: {\sc OreModules}: A symbolic package for the study of multidimensional linear systems,
Springer, Lecture Notes in Control and Inform. Sci., 352, 233-264, 2007.\\
http://wwwb.math.rwth-aachen.de/OreModules  \\
\noindent
[3] A. QUADRAT: An Introduction to Constructive Algebraic Analysis and its Applications, 
Les cours du CIRM, Journ\'ees Nationales de Calcul Formel (2010), 1(2), 281-471, 2010.\\
\noindent
[5] J. J. ROTMAN: An Introduction to Homological Algebra, Springer, 2009.  \\

\end{document}